\documentclass[preprint]{elsarticle}
\usepackage{lineno,hyperref}
\modulolinenumbers[5]
\journal{Journal of \LaTeX\ Templates}
\usepackage{graphics}
\usepackage{epsfig}
\usepackage{amsmath}
\usepackage{amsthm}

\usepackage{tikz}
\usetikzlibrary{
    positioning,
    arrows.meta,
    fit,
    calc,
    shapes.geometric,
    backgrounds
}
\usepackage{float}
\usepackage{tabularx}      
\usepackage{array}         
\usepackage{placeins}  
\usepackage{booktabs, multirow}
\newtheorem{problem}{Problem}
\usepackage{amsfonts}
\usepackage{cases}
\usepackage{{color}}
\usepackage{graphicx}
\usepackage{caption}
\usepackage{subcaption}
\usepackage{natbib}
\usepackage{placeins}
\usepackage{makecell}

\usepackage[ruled,vlined]{algorithm2e}

\catcode`@=11
\@addtoreset{equation}{section}
\renewcommand\theequation{\thesection.\@arabic\c@equation}
\catcode`@=12

\begin{document}
	
\begin{frontmatter}
\title{\bf A Legendre-Based Modified Extreme Learning Machine for Convection-Diffusion Equations: Steady-State, Continuous-Time, and Discrete-Time Formulations \tnoteref{mytitlenote}}
		
\author[1]{Gadade Navnath Ankush}
\ead{navnathagadade@gmail.com}
\author[2]{Ram Shiromani\corref{cor1}}
\ead{ram.panday.786@gmail.com}	
\cortext[cor1]{Corresponding author}
\author[4]{Tom\'as Gandara}
\ead{tomas.gandara.imilqueo@upc.edu}	
\author[4,5]{Ramon Codina}
\ead{ramon.codina@upc.edu}	
\author[1,2,3]{Sivaram Ambikasaran}
\ead{sivaambi@dsai.iitm.ac.in}

\address[1]{Department of Mathematics, IIT Madras, Chennai, India}
\address[2]{Wadhwani School of Data Science and Artificial Intelligence, IIT Madras, Chennai, India}
 
\address[3]{Department of Data Science and Artificial Intelligence, IIT Madras, Chennai, India}

\address[4]{Universitat Polit\`ecnica de Catalunya, Jordi Girona 1-3, Edifici C1, 08034 Barcelona, Spain}

\address[5]{Centre Internacional de M\`etodes Num\`erics en Enginyeria (CIMNE), Edifici C1, Campus Nord UPC, Gran Capit\`a S/N, 08034 Barcelona, Spain}

\begin{abstract}
Convection--diffusion equations are fundamental models for transport phenomena, but their numerical solution becomes challenging in convection-dominated regimes due to steep gradients and sharp boundary layers. This motivates the development of accurate approaches capable of effectively resolving such features while avoiding costly iterative optimization. In this study, a modified extreme learning machine (MELM) framework based on shifted Legendre polynomials is developed for convection--diffusion equations in steady-state, continuous-time (CT-MELM), and discrete-time (DT-MELM) formulations. The proposed approach is evaluated on steady and unsteady benchmark problems with relatively high Péclet numbers, including continuous and discontinuous boundary conditions. Finite element method (FEM) solutions are used as independent solution references. By exploiting the linearity of the governing equations, MELM determines the output weights through a single least-squares solve, avoiding the iterative gradient-based training required by physics-informed neural networks (PINNs). The shifted Legendre polynomial basis provides an effective spectral representation for resolving steep gradients and sharp boundary layers. For the DT-MELM formulation, the coefficient matrix remains unchanged across successive time levels, allowing its LU or QR factorization to be computed once and reused. A systematic comparison with PINNs that incorporate advanced training strategies demonstrates that MELM provides higher accuracy and shorter solution times for steady-state and DT-MELM formulations across the benchmark problems considered. In contrast, CT-MELM exhibits scalability limitations as the size of the space--time basis increases. Overall, the results demonstrate the effectiveness of MELM for steady-state and discrete-time convection--diffusion problems, while also showing that PINNs provide better scalability for large space--time formulations where global CT-MELM becomes computationally demanding.

\end{abstract}
		
\end{frontmatter}
	
\section {Introduction}
The convection--diffusion equation describes the transport of a scalar quantity, such as heat, mass, or chemical concentration, subject to the simultaneous action of advection by a flow field and diffusive smoothing. It constitutes a fundamental model in a wide range of applications, including the dispersion of pollutants in air and water \cite{stockie2011mathematics}, heat transfer in solids and fluids \cite{patankar1980numerical}, and charge transport in semiconductor devices \cite{selberherr1984analysis}. The qualitative behaviour of the solution is determined by the relative strength of the two mechanisms rather than by either in isolation. In the convection-dominated regime, diffusive smoothing is insufficient to balance the transport induced by the flow, and the solution, although smooth over the greater part of the domain, exhibits a sharp gradient localised within a narrow region adjacent to an outflow boundary or along an internal front \cite{morton1996numerical,roos2008robust}. The resolution of such sharp gradients constitutes the principal difficulty in the numerical treatment of the problem, since discretisations that fail to capture them are prone to spurious oscillations that degrade the solution throughout the computational domain.

{

Standard Galerkin finite element and central finite difference methods may exhibit spurious oscillations in convection-dominated regimes unless the discretization adequately resolves the associated layers \cite{quarteroni2017diffusion,onate2016fic}. Established approaches to improve robustness include stabilization of the numerical formulation and adaptation of the computational mesh. Stabilized methods, such as SUPG, GLS, variational multiscale, and FIC formulations, introduce additional control of the convective contribution \cite{CODINA1998185,onate2016fic,quarteroni2017diffusion}. Alternatively, Shishkin- and Bakhvalov-type meshes concentrate computational points in layer regions and can provide parameter-uniform convergence~\cite{roos2008robust}, while adaptive $h$- and $hp$-refinement strategies identify regions requiring additional resolution through a posteriori error estimation \cite{mitchell2013collection,demkowicz2002hp,schmidt2000posteriori}. These approaches provide effective and well-established remedies for convection-dominated problems, but their performance depends on the selected stabilization or mesh-refinement strategy. This motivates the investigation of neural-network-based approximations as an alternative, mesh-free solution framework.

}

{
Over the past few years, neural-network-based approaches have emerged as an alternative framework for the numerical solution of partial differential equations (PDEs). Early neural-network formulations for differential equations were introduced by Lagaris \emph{et al.}~\cite{lagaris1998ann,lagaris2000ann}, while the modern physics-informed neural network (PINN) framework was established by Raissi \emph{et al.}~\cite{Raissi2019}. PINNs approximate the solution field using a neural network while incorporating the governing differential equation and the associated initial and boundary conditions into the training objective. The required derivatives are computed using automatic differentiation~\cite{Baydin2018}, yielding a mesh-free formulation that does not require labelled solution data. In its standard form, commonly referred to as a vanilla PINN, the solution is represented by a fully connected neural network whose parameters are determined through gradient-based optimization.

Despite their flexibility, vanilla PINNs may encounter significant training difficulties for stiff, multiscale, and convection-dominated PDEs, where steep gradients and localized solution features can exacerbate spectral bias, imbalanced loss contributions, and unfavorable optimization dynamics~\cite{krishnapriyan2021failure,wang2021gradient,Cuomo2022,Luo2025}. These challenges have motivated developments in network architectures, parameterizations, loss construction, collocation strategies, and optimization procedures~\cite{Cuomo2022,Luo2025}. Convection--diffusion problems provide a particularly relevant setting for evaluating these developments, as large Péclet numbers give rise to thin layers and steep solution gradients. Recent studies have specifically investigated PINNs for convection-dominated convection--diffusion problems~\cite{Frerichs-Mihov2026}. In the present work, PINNs are employed as a representative mesh-free solver for comparison with the proposed approach.

For transient PDEs, PINNs can be formulated using either continuous-time (CT) or discrete-time (DT) approaches. In the CT formulation, a single neural network represents the solution over the entire space--time domain, whereas DT formulations introduce the temporal evolution through a time-stepping scheme. The original DT-PINN formulation of Raissi \emph{et al.}~\cite{Raissi2019} employs implicit Runge--Kutta time integration, with the neural network representing the intermediate stages and the solution at the next time level. Subsequent work has extended DT-PINNs to other step-by-step and time-integration strategies~\cite{https://doi.org/10.1002/gamm.202100006,VALENTINO2025541}. In the present work, the temporal derivative is discretized externally using backward differentiation formula (BDF), while a spatial PINN solves the resulting boundary-value problem at each discrete time level. The same BDF formulation and time-step size are used for the proposed method, providing a common temporal discretization for comparing the two spatial solution strategies.

}

The challenges associated with iterative training in PINNs have motivated the development of alternative neural-network approaches that retain their mesh-free character while reducing the computational cost of optimization. The extreme learning machine (ELM), introduced by Huang \emph{et al.}~\cite{Huang2006}, avoids iterative training by fixing the hidden-layer parameters and determining the output weights directly through a linear least-squares problem.

ELM-based approaches have subsequently been applied to the numerical solution of differential equations. Calabr\`o \emph{et al.}~\cite{calabro2021elm} developed an ELM-based collocation method for elliptic PDEs and demonstrated its effectiveness for problems with boundary layers. Their subsequent work~\cite{calabro2023time} investigated time-discretization strategies for solving parabolic PDEs using artificial neural networks. Dwivedi and Srinivasan~\cite{Dwivedi2019} introduced the physics-informed extreme learning machine (PIELM), which incorporates the governing equations and boundary conditions (BCs) into the ELM framework. However, their results also indicated that globally supported PIELM may become less effective for problems containing steep gradients and boundary layers. To improve the representation of localized solution features, distributed and local ELM formulations have been developed using domain decomposition~\cite{Dwivedi2021,dong2021local}, while kernel-adaptive PIELM has been proposed using adaptive radial basis functions~\cite{dwivedi2025kernel}. These developments improve the flexibility of ELM-based approximations but introduce additional mechanisms for resolving localized features.

The hidden-layer representation has also been enhanced through Legendre polynomial-based ELM formulations. Mall and Chakraverty~\cite{Mall2016} employed a Legendre neural network for solving ordinary differential equations, while Yang \emph{et al.}~\cite{Yang2018} developed an improved ELM based on a Legendre neural network for ODEs and subsequently extended the approach to elliptic PDEs~\cite{Yang2020}. These studies demonstrate the potential of Legendre-based representations to improve the approximation capability of ELMs while retaining their non-iterative solution structure. However, their applications have primarily focused on ODEs and elliptic PDEs, leaving the treatment of convection--diffusion problems with high Péclet numbers and sharp boundary or interior layers largely unexplored.

Motivated by these developments, we propose a Legendre polynomial-based modified extreme learning machine (MELM) for convection--diffusion equations with sharp gradients. The proposed MELM retains the fundamental structure of ELM and provides a coordinate-wise input representation for multidimensional problems.  It is positioned at the intersection of deterministic ELM-based function approximation and Legendre spectral approximation, retaining the ELM structure while using shifted Legendre polynomials as a deterministic hidden-layer representation.  The proposed framework includes steady-state, CT-MELM, and DT-MELM formulations, providing a unified approach for steady and transient problems. The DT-MELM formulation employs BDF time discretization and solves the resulting spatial problem through a linear least-squares system at each time level. A key advantage of the DT-MELM formulation is that its coefficient matrix remains unchanged throughout the time integration, allowing its matrix factorization to be computed once and efficiently reused at subsequent time levels. This leads to a non-iterative and computationally efficient framework for time-dependent convection--diffusion problems with sharp gradients and layers.

The remainder of this article is organized as follows. Section~\ref{sec:sec2} presents the mathematical formulation of the convection--diffusion problems. Section~\ref{sec:sec3} presents the learning-based solution framework, including the PINN and proposed MELM approaches. Section~\ref{sec:sec4} describes their numerical implementation. Section~\ref{sec:sec5} presents and discusses the numerical results. Finally, Section~\ref{sec:sec6} concludes the article and outlines potential directions for future work.

\section{Mathematical Formulation}
\label{sec:sec2}

This section presents the mathematical formulation of the convection--diffusion problem, including its governing equation, initial and boundary conditions, and non-dimensional form.

\subsection{ Governing Equations and Problem Setup}

 We shall restrict our attention to 2D problems, with 1D cases considered where appropriate, although the extension to 3D of the ideas presented is straightforward. Let $\Omega\subset\mathbb{R}^{2}$ be a bounded computational domain with
boundary $\partial\Omega$, and let $T>0$ denote the final time. Let
\begin{equation*}
\nabla=\left(\frac{\partial}{\partial x},\frac{\partial}{\partial y}\right),
\qquad
\nabla^2=\frac{\partial^2}{\partial x^2}+\frac{\partial^2}{\partial y^2},
\end{equation*}
denote the gradient and Laplacian operators, respectively, and let
$\mathbf{a}=(a_x,a_y)^T$ denote the constant convection velocity vector.
Consider the two-dimensional transient convection--diffusion equation
\begin{equation}
\frac{\partial \phi}{\partial t}
+\mathbf{a}\cdot\nabla\phi
=
k\,\nabla^2\phi
+f(x,y,t),
\qquad
(x,y,t)\in\Omega\times(0,T],
\label{eq:cde}
\end{equation}
where $\phi(x,y,t)$ denotes the unknown scalar field, $k>0$ is the
constant isotropic diffusivity, and $f(x,y,t)$ denotes the prescribed
source term. Equation~\eqref{eq:cde} is supplemented with the
initial and Dirichlet BCs
\begin{equation}
\phi(x,y,t)=
\begin{cases}
\phi_0(x,y), & (x,y)\in\Omega,\; t=0,\\[4pt]
g(x,y,t), & (x,y,t)\in\partial\Omega\times(0,T].
\end{cases}
\label{eq:icbc}
\end{equation}

\subsection{ Non-dimensionalisation and Scaling Analysis}

{ To express the governing equation in dimensionless form and characterize the relative scaling of the convective and diffusive contributions, characteristic scales are introduced for length, velocity, and the solution field. Let $L$ denote the characteristic length, $U$ the characteristic convection velocity, and $\phi_c$ the characteristic solution value. The corresponding dimensionless variables are defined as}
\begin{equation*}
x'=\frac{x}{L},\qquad
y'=\frac{y}{L},\qquad
t'=\frac{Ut}{L},\qquad
a_x'=\frac{a_x}{U},\qquad
a_y'=\frac{a_y}{U},\qquad
\phi'=\frac{\phi}{\phi_c},\qquad
f'=\frac{Lf}{U\phi_c}.
\end{equation*}
The dimensionless gradient and Laplacian operators are given by
\begin{equation*}
\nabla'=
\left(
\frac{\partial}{\partial x'},
\frac{\partial}{\partial y'}
\right),
\qquad
\nabla'^2=
\frac{\partial^2}{\partial x'^2}
+
\frac{\partial^2}{\partial y'^2}.
\end{equation*}
Substituting these variables into Equation~\eqref{eq:cde} gives
\begin{equation}
\frac{\partial \phi'}{\partial t'}
+\mathbf{a}'\cdot\nabla'\phi'
=
\frac{1}{Pe}\nabla'^2\phi'
+f',
\qquad
(x',y',t')\in\Omega'\times(0,T'],
\label{eq:nondim}
\end{equation}
where
\begin{equation*}
T'=\frac{UT}{L},
\qquad
\mathbf{a}'=(a_x',a_y')^T,
\qquad
Pe=\frac{UL}{k}.
\end{equation*}

{ The Péclet number $Pe$ characterizes the relative importance of convective and diffusive transport and therefore provides a measure of the degree of convection dominance of the problem. The characteristic scales thus provide a common physical scaling for the terms appearing in the governing equation and define the dimensionless formulation used throughout the numerical study. }

\section{Numerical Solution Framework}
\label{sec:sec3}

{
This section presents the numerical solution frameworks used in this work, including the temporal discretization, finite-element reference solver, PINN formulation, and proposed MELM approach.

}

\subsection{Temporal discretization}
\label{sec:temporal_discretization}

For transient problems, the temporal derivative in
Equation~\eqref{eq:cde} is discretized using the BDF. Let $\phi^{(n)}(\mathbf{x})$ denote the numerical approximation to $\phi(\mathbf{x},t_n)$ at time level $t_n$, the first- and second-order BDF schemes (BDF1 and BDF2) approximate the
time derivative at $t=t_{n+1}$ as
\begin{equation}
\left.\frac{\partial\phi}{\partial t}\right|_{t=t_{n+1}}
\approx
\frac{\phi^{(n+1)}-\phi^{(n)}}{\Delta t},
\label{eq:bdf1}
\end{equation}
and
\begin{equation}
\left.\frac{\partial\phi}{\partial t}\right|_{t=t_{n+1}}
\approx
\frac{3\phi^{(n+1)}-4\phi^{(n)}+\phi^{(n-1)}}{2\Delta t},
\label{eq:bdf2}
\end{equation}
respectively. Substituting either approximation into
Equation~\eqref{eq:cde} transforms the transient problem into a
sequence of time-discrete boundary value problems. { BDF2 is adopted as the primary temporal discretization for the DT formulations and for the transient finite element method (FEM) reference calculations.} Since the BDF2 scheme
requires solutions at two previous time levels, the first time step is
computed using BDF1, while BDF2 is employed for all subsequent time
steps. 

{

\subsection{Finite-Element reference solutions}

The FEM reference solutions are obtained using a stabilized finite-element solver developed within the research group. The spatial discretization employs the standard Galerkin method within the Variational Multiscale (VMS) framework, with residual-based stabilization provided by the Algebraic Sub-Grid Scale (ASGS)~\cite{HUGHES19983,CODINA200061}. First-order Lagrange elements are employed throughout, using linear ($P_1$) elements for one-dimensional problems and bilinear quadrilateral ($Q_1$) elements for two-dimensional problems. Stabilization is required for the convection-dominated problems considered here, as the standard Galerkin formulation can produce spurious oscillations in such regimes. The ASGS formulation mitigates these oscillations while preserving the consistency of the finite-element discretization. However, it does not guarantee monotonicity, and small localized oscillations may persist even on refined meshes.

The stabilization parameters are evaluated element-wise and depend on the local mesh size, diffusivity, and convective velocity, following the formulation in~\cite{CODINA200061,Codina2018}. The same ASGS formulation is used consistently across all convection--diffusion problems considered in this work.

For transient problems, temporal discretization follows the BDF formulation described in Section~\ref{sec:temporal_discretization}. At each time level, the resulting linear system is solved. The spatial meshes, time-step sizes, and other discretization parameters are specified for each benchmark in Section~\ref{sec:sec5}.

Mesh resolutions are selected so as to adequately resolve steep gradients and boundary layers, guided by the characteristic length scales of each benchmark. The benchmark-dependent convergence tolerances are used consistently within each class of simulations. Where analytical solutions are available, FEM results serve as an independent numerical comparison. For the remaining benchmarks, the FEM solutions obtained on the adopted meshes are used as reference solutions for evaluating the PINN and MELM approximations.

\subsection{Physics-Informed Neural Networks}\label{subsec:PINNs}

PINNs, introduced by Raissi \emph{et al.}~\cite{Raissi2019}, provide a mesh-free framework for solving PDEs by approximating the unknown solution field with a neural network and enforcing the governing physics through a composite loss function. Unlike conventional supervised learning approaches that require labeled data, PINNs incorporate the differential equation, initial conditions, and BCs directly into the training objective. The required derivatives are computed via automatic differentiation~\cite{Baydin2018}, enabling the network to satisfy the physics while minimizing the physics-informed loss. In this work, PINNs are employed exclusively as physics-based solvers, without observational data.

Consider the convection--diffusion equation~\eqref{eq:cde} subject to the initial and Dirichlet BCs~\eqref{eq:icbc}. Let $\hat{\phi}_\theta(\mathbf{x},t)$ denote a neural network approximation of the solution, parameterized by $\theta$ (weights and biases). Here and hereafter, $\mathbf{x}=(x,y)^{T}$ denotes
the spatial coordinate vector. The PDE residual is defined as 
\begin{equation}
\mathcal{R}[\hat{\phi}_\theta](\mathbf{x},t)
:=
\frac{\partial \hat{\phi}_\theta}{\partial t}
+\mathbf{a}\cdot\nabla\hat{\phi}_\theta
-
k\nabla^2\hat{\phi}_\theta
-
f(\mathbf{x},t).
\label{eq:pinn_residual}
\end{equation}

The network parameters are determined by minimizing the following loss function
\begin{equation}
\mathcal{L}(\theta)
=
\lambda_{\mathrm{PDE}}\mathcal{L}_{\mathrm{PDE}}(\theta)
+
\lambda_{\mathrm{BC}}\mathcal{L}_{\mathrm{BC}}(\theta)
+
\lambda_{\mathrm{IC}}\mathcal{L}_{\mathrm{IC}}(\theta),
\label{eq:pinn_loss}
\end{equation}
where the individual loss terms are defined as
\begin{align}
\mathcal{L}_{\mathrm{PDE}}(\theta)
&=
\frac{1}{N_{\mathrm{PDE}}}
\sum_{i=1}^{N_{\mathrm{PDE}}}
\left|
\mathcal{R}[\hat{\phi}_\theta](\mathbf{x}_i^r, t_i^r)
\right|^2,
\label{eq:pinn_L_pde} \\
\mathcal{L}_{\mathrm{BC}}(\theta)
&=
\frac{1}{N_{\mathrm{BC}}}
\sum_{i=1}^{N_{\mathrm{BC}}}
\left|
\hat{\phi}_\theta(\mathbf{x}_i^{bc}, t_i^{bc})
-
g(\mathbf{x}_i^{bc}, t_i^{bc})
\right|^2,
\label{eq:pinn_L_bc} \\
\mathcal{L}_{\mathrm{IC}}(\theta)
&=
\frac{1}{N_{\mathrm{IC}}}
\sum_{i=1}^{N_{\mathrm{IC}}}
\left|
\hat{\phi}_\theta(\mathbf{x}_i^{ic}, 0)
-
\phi_0(\mathbf{x}_i^{ic})
\right|^2.
\label{eq:pinn_L_ic}
\end{align}
Here, $\{\mathbf{x}_i^r, t_i^r\}_{i=1}^{N_{\mathrm{PDE}}}$, $\{\mathbf{x}_i^{bc}, t_i^{bc}\}_{i=1}^{N_{\mathrm{BC}}}$, and $\{\mathbf{x}_i^{ic}\}_{i=1}^{N_{\mathrm{IC}}}$ denote sets of collocation points sampled from the interior, boundary, and initial domains, respectively. The parameters $\lambda_{\mathrm{PDE}}$, $\lambda_{\mathrm{BC}}$, and $\lambda_{\mathrm{IC}}$ are non-negative weighting coefficients that balance the relative contributions of each loss term.

In its standard form (vanilla PINNs), the neural network is represented by a fully connected multi-layer perceptron (MLP) that maps the input coordinates to the predicted solution. The network parameters are typically initialized using the Glorot scheme~\cite{pmlr-v9-glorot10a} and optimized through gradient-based methods such as Adam~\cite{kingma2014adam}. In the present work, the governing equations are formulated using the characteristic scales introduced in Section~\ref{sec:sec2}, which define the corresponding normalized formulation. Independently of this physical scaling, the spatial and, where applicable, temporal coordinates supplied to the neural network are normalized using the characteristic domain scales. This input normalization is applied only at the network-input level and does not modify the governing equations. The PINN parameters are then optimized using a prescribed learning-rate schedule, with the corresponding initial learning rate and decay settings specified for each benchmark in Section~\ref{sec:sec4}.

Although the vanilla formulation provides a general and effective framework for many PDEs, its training may become challenging for convection-dominated problems, where steep gradients and thin boundary layers can exacerbate spectral bias, imbalanced loss contributions, and unfavorable optimization dynamics. To address these difficulties, the present work considers a set of complementary network architectures, parameterization strategies, loss-balancing mechanisms, and optimization methods described in the following subsection and implemented selectively according to the requirements of each benchmark problem.

\subsubsection{Advanced training strategies for PINNs}\label{sec:pinn_advancedtools}

The numerical experiments employ a set of complementary strategies developed in the PINN literature and implemented through the JAXPI framework (see~\cite{Wang2023} for details). These strategies target different aspects of the neural-network solver, including feature representation, parameterization, loss balancing, temporal training, and optimization. The specific components considered in this work are summarized below.

\begin{itemize}

\item \textbf{Modified MLPs.}
The Modified MLP architecture, proposed by Wang et al.~\cite{wang2021gradient}, improves gradient propagation during PINN training by introducing trainable feature gates that modify the information flow through the hidden layers of the standard MLP architecture. This modification mitigates gradient-flow pathologies associated with stiff PDE losses and improves the representation of sharp solution features.

\item \textbf{Fourier feature embeddings.}
Random Fourier features~\cite{Tancik2020} augment the input representation by projecting the coordinates into a higher-dimensional space with a broader frequency spectrum. This addresses the spectral bias of standard MLPs and improves the representation of higher-frequency spatial features associated with steep gradients and thin boundary layers.

\item \textbf{Random weight factorization.}
Random Weight Factorization (RWF)~\cite{Wang2022RWF} is a reparameterization strategy that improves the optimization stability and reduces the sensitivity of the network to parameter initialization, contributing to more robust training in stiff regimes.

\item \textbf{Adaptive loss balancing.}
The relative contributions of the PDE, boundary, and initial-condition losses can become severely imbalanced during training, particularly in the presence of steep gradients. Several adaptive weighting strategies have been proposed to address this issue, including methods based on the Neural Tangent Kernel (NTK) theory~\cite{Wang2022NTK}. In this work, we employ the NTK-based strategy to dynamically adjust the importance of each loss component, as it has been shown to effectively balance the convergence rates of the different loss terms during PINN training.

\item \textbf{Causal training.}
For transient problems, the temporal structure of the solution must be respected. Causal weighting~\cite{Wang2024} enforces temporal ordering by reducing the influence of later-time residuals before earlier-time dynamics have been adequately learned, preventing the optimization from focusing on future times prematurely. This strategy is employed when collocation points are batched over temporal segments and, in the present work, is applied only to continuous-time (CT) approximations.

\item \textbf{SOAP optimizer.}
The SOAP optimizer~\cite{vyas2025soapimprovingstabilizingshampoo} is a quasi-second-order method that introduces structured preconditioning through Kronecker-factored approximations of the gradient covariance. Wang et al.~\cite{NEURIPS2025_f6557065} demonstrated its effectiveness for PINNs, showing that it mitigates gradient conflicts and yields significant improvements in training efficiency compared with standard first-order methods such as Adam.

\end{itemize}

The specific combination of architectures, sampling procedures, optimization algorithms, and training strategies is problem-dependent and is therefore specified in the implementation tables in Section~\ref{sec:sec4}. These components do not modify the governing PDE formulation, but constitute numerical strategies used to improve the robustness of the PINN baseline for the convection-diffusion problems considered in this work.

\subsection{Legendre-based Modified Extreme Learning Machine}

This section develops the proposed Legendre-based modified extreme
learning machine (Legendre-MELM) framework. First, the construction of
shifted Legendre polynomials on general physical domains is introduced.
The resulting basis is then incorporated into the MELM formulation for
one-dimensional and two-dimensional problems, with multidimensional
inputs treated coordinate-wise and tensor-product basis functions
constructed accordingly. The analytical derivatives of the resulting
approximate solution are subsequently derived for evaluating the
differential equations, initial conditions, and BCs. Hereafter, the proposed method is referred to as the MELM framework.

\subsubsection{Shifted Legendre Polynomials on Intervals and Rectangular Domains}

The standard Legendre polynomials $P_n(\xi)$ are defined on the reference
interval $[-1,1]$, where $\xi$ denotes the reference coordinate. The first few polynomials are
\begin{equation*}
P_0(\xi)=1,\qquad
P_1(\xi)=\xi,\qquad
P_2(\xi)=\frac{1}{2}(3\xi^2-1),\qquad
P_3(\xi)=\frac{1}{2}(5\xi^3-3\xi).
\end{equation*}
These polynomials satisfy the orthogonality relation
\begin{equation*}
\int_{-1}^{1} P_m(\xi)\,P_n(\xi)\,d\xi
=
\frac{2}{2n+1}\delta_{mn},
\end{equation*}
where $m$ and $n$ denote the polynomial degrees, and $\delta_{mn}$ is
the Kronecker delta.

To define these polynomials on a general interval $[a,b]$, the physical
coordinate $x$ is mapped onto the reference interval through the linear
transformation
\begin{equation*}
\xi(x)=\frac{2x-(a+b)}{b-a}, \qquad x\in[a,b], \; \xi(x)\in[-1,1].
\end{equation*}
The shifted Legendre polynomials over the general interval $[a,b]$ are
defined by
\begin{equation*}
\tilde{P}_n(x)
=
P_n(\xi(x))
=
P_n\!\left(\frac{2x-(a+b)}{b-a}\right),
\qquad n=0,1,\ldots,N,
\end{equation*}
where $N$ denotes the maximum polynomial degree. In subsequent sections, the notation $N_x$, $N_y$, and analogous symbols are used to denote the maximum polynomial degrees in the corresponding coordinate directions.
For multidimensional problems, the above mapping is applied
independently to each spatial and temporal variable, yielding
tensor-product shifted Legendre polynomials. According to this, we restrict ourselves to rectangular domains, although the extension to domains that are unions of rectangles is also feasible if possible discontinuities are properly accounted for. Extension to general domains could be based on mapping them to unions of rectangles, if possible, or by using embedding techniques. These extensions will not be studied in this work. 

\subsubsection{One-Dimensional MELM Formulation}

The proposed formulation retains the fundamental ELM principle of a
fixed input-to-hidden representation and a linear output layer without a
bias term. For a one-dimensional input $x$, the input-to-hidden mapping is
written as
\begin{equation*}
z=wx+b,
\end{equation*}
where $w$ and $b$ denote the input weight and hidden-layer bias,
respectively. In the present formulation, these parameters are fixed as
\(w=1, b=0, \)
so that $z=x$. The hidden-layer representation is constructed using
shifted Legendre polynomials as
\begin{equation*}
\psi_n(x)=\tilde{P}_n(x),
\end{equation*}
where $\tilde{P}_n$ denotes the $n$-th degree shifted Legendre
polynomial. These polynomials serve as the hidden-layer activation
functions (basis functions) in place of conventional nonlinear
activation functions and define the approximation space. The approximate
solution is consequently represented by a linear output layer without a
bias term as
\begin{equation*}
\hat{\phi}(x)
=
\sum_{n=0}^{N_x}\beta_n\tilde{P}_n(x),
\end{equation*}
where $\beta_n$ are the output-layer coefficients to be determined and
$N_x$ denotes the maximum polynomial degree.

\subsubsection{Two-Dimensional MELM Formulation}

For a two-dimensional problem with input coordinates $(x,y)$, the
one-dimensional formulation is extended by treating the two spatial
coordinates separately. In the classical ELM formulation, a
multidimensional input
$\mathbf{x}=(x,y)^T$ is transformed into a scalar hidden-layer input
through the affine mapping $\mathbf{w}^{T}\mathbf{x}+b$, where
$\mathbf{w}=(w_1,w_2)^T$ denotes the input-weight vector and $b$ is the
hidden-layer bias. In the present formulation, instead of combining the
two coordinates into a single scalar input, they are retained
individually through the component-wise fixed mapping
\begin{equation*}
\mathbf{z}
=
(w_1x+b,\;w_2y+b),
\end{equation*}
with \(w_1=w_2=1, b=0.\)
Consequently, the transformed input becomes \(\mathbf{z}=(x,y).\)
The two-dimensional hidden-layer features are then constructed as
tensor products of the corresponding one-dimensional shifted Legendre
basis functions:
\begin{equation*}
\psi_{nm}(x,y)
=
\tilde{P}_n(x)\tilde{P}_m(y).
\end{equation*}
The approximate solution is subsequently expressed using the linear
output layer without a bias term as
\begin{equation*}
\hat{\phi}(x,y)
=
\sum_{n=0}^{N_x}
\sum_{m=0}^{N_y}
\beta_{nm}
\tilde{P}_n(x)
\tilde{P}_m(y),
\end{equation*}
where $\beta_{nm}$ are the output-layer coefficients to be determined,
and $N_x$ and $N_y$ denote the maximum polynomial degrees in the
$x$- and $y$-directions, respectively.

Thus, the two-dimensional MELM formulation retains the fundamental ELM
structure of a fixed input-to-hidden representation followed by a
linear output layer, while modifying the multidimensional hidden-layer
representation by treating each coordinate separately and constructing
the corresponding tensor-product Legendre basis. The same
coordinate-wise construction can be extended to three-dimensional and
higher-dimensional problems. For a two-dimensional problem with input variables $(x,y)$, the network architecture is illustrated in Figure~\ref{fig:2d_legendre_melm}.

\begin{figure}[htbp]
\centering

\begin{tikzpicture}[
    font=\small,
    box/.style={
        draw,
        rounded corners=2pt
    },
    basis/.style={
        draw=blue!55!black,
        fill=white,
        ellipse,
        minimum width=4.5cm,
        minimum height=7mm,
        font=\scriptsize,
        align=center
    }
]


\node[
    box,
    draw=green!50!black,
    fill=green!8,
    minimum width=2.6cm,
    minimum height=7.0cm
] (input) {};

\node[
    box,
    draw=blue!50!black,
    fill=blue!5,
    minimum width=7.5cm,
    minimum height=8cm,
    right=6mm of input
] (hidden) {};

\node[
    box,
    draw=violet!50!black,
    fill=violet!7,
    minimum width=2.8cm,
    minimum height=7.0cm,
    right=6mm of hidden
] (output) {};


\node[
    font=\bfseries
] at ([yshift=23mm]input.north)
{Input layer};

\node[
    font=\scriptsize
] at ([yshift=15mm]input.north)
{2D input};

\node[
    draw=green!50!black,
    fill=green!15,
    circle,
    minimum size=13mm
] (xy) at ([yshift=5mm]input.center)
{$(x,y)$};

\node[
    font=\scriptsize
] at ([yshift=-22mm]input.center)
{$\boldsymbol{x}=(x,y)$};


\node[
    font=\bfseries
] at ([yshift=23mm]hidden.north)
{Hidden layer (Basis Functions)};

\node[
    font=\scriptsize
] at ([yshift=15mm]hidden.north)
{Tensor-product shifted Legendre polynomials};


\node[basis] (phi00)
at ([xshift=-7mm,yshift=35mm]hidden.center)
{$
\psi_{00}(x,y)
=
\widetilde P_0(x)\widetilde P_0(y)
$};

\node[basis] (phi01)
at ([xshift=-7mm,yshift=23mm]hidden.center)
{$
\psi_{01}(x,y)
=
\widetilde P_0(x)\widetilde P_1(y)
$};

\node[font=\large]
at ([xshift=-7mm,yshift=17mm]hidden.center)
{$\vdots$};

\node[basis] (phi0N)
at ([xshift=-7mm,yshift=9mm]hidden.center)
{$
\psi_{0N_y}(x,y)
=
\widetilde P_0(x)\widetilde P_{N_y}(y)
$};

\node[basis] (phi10)
at ([xshift=-7mm,yshift=-3mm]hidden.center)
{$
\psi_{10}(x,y)
=
\widetilde P_1(x)\widetilde P_0(y)
$};

\node[font=\large]
at ([xshift=-7mm,yshift=-9mm]hidden.center)
{$\vdots$};

\node[basis] (phi1N)
at ([xshift=-7mm,yshift=-19mm]hidden.center)
{$
\psi_{1N_y}(x,y)
=
\widetilde P_1(x)\widetilde P_{N_y}(y)
$};

\node[font=\large]
at ([xshift=-7mm,yshift=-26mm]hidden.center)
{$\vdots$};

\node[basis] (phiNN)
at ([xshift=-7mm,yshift=-35mm]hidden.center)
{$
\psi_{N_xN_y}(x,y)
=
\widetilde P_{N_x}(x)\widetilde P_{N_y}(y)
$};


\node[
    font=\bfseries
] at ([yshift=23mm]output.north)
{Output layer};

\node[
    font=\scriptsize
] at ([yshift=15mm]output.north)
{Linear combination};

\node[
    draw=violet!60!black,
    fill=violet!5,
    rounded corners=2pt,
    minimum width=2.2cm,
    minimum height=11mm
] (u)
at (output.center)
{$\hat{\phi}(x,y)$};


\draw[->,thick]
(xy.east)
--
(hidden.west |- xy.east);


\draw[->]
(phi00.east)
--
node[above,font=\scriptsize] {$\beta_{00}$}
(u.west);

\draw[->]
(phi01.east)
--
node[above,font=\scriptsize] {$\beta_{01}$}
(u.west);

\draw[->]
(phi0N.east)
--
node[above,font=\scriptsize] {$\beta_{0N_y}$}
(u.west);

\draw[->]
(phi10.east)
--
node[above,font=\scriptsize] {$\beta_{10}$}
(u.west);

\draw[->]
(phi1N.east)
--
node[below,font=\scriptsize] {$\beta_{1N_y}$}
(u.west);

\draw[->]
(phiNN.east)
--
node[below,font=\scriptsize] {$\beta_{N_xN_y}$}
(u.west);


\node[
    draw=blue!50!black,
    fill=blue!8,
    rounded corners=2pt,
    minimum width=14.5cm,
    minimum height=9mm,
    font=\bfseries
] at ([yshift=7mm]hidden.north)
{Two-Dimensional Legendre-MELM Architecture};

\end{tikzpicture}

\vspace{2mm}

\caption{Architecture of the proposed two-dimensional Legendre-MELM.
The two-dimensional input $(x,y)$ is mapped component-wise to the
tensor-product shifted Legendre basis functions, followed by a linear
output layer.}

\label{fig:2d_legendre_melm}

\end{figure}
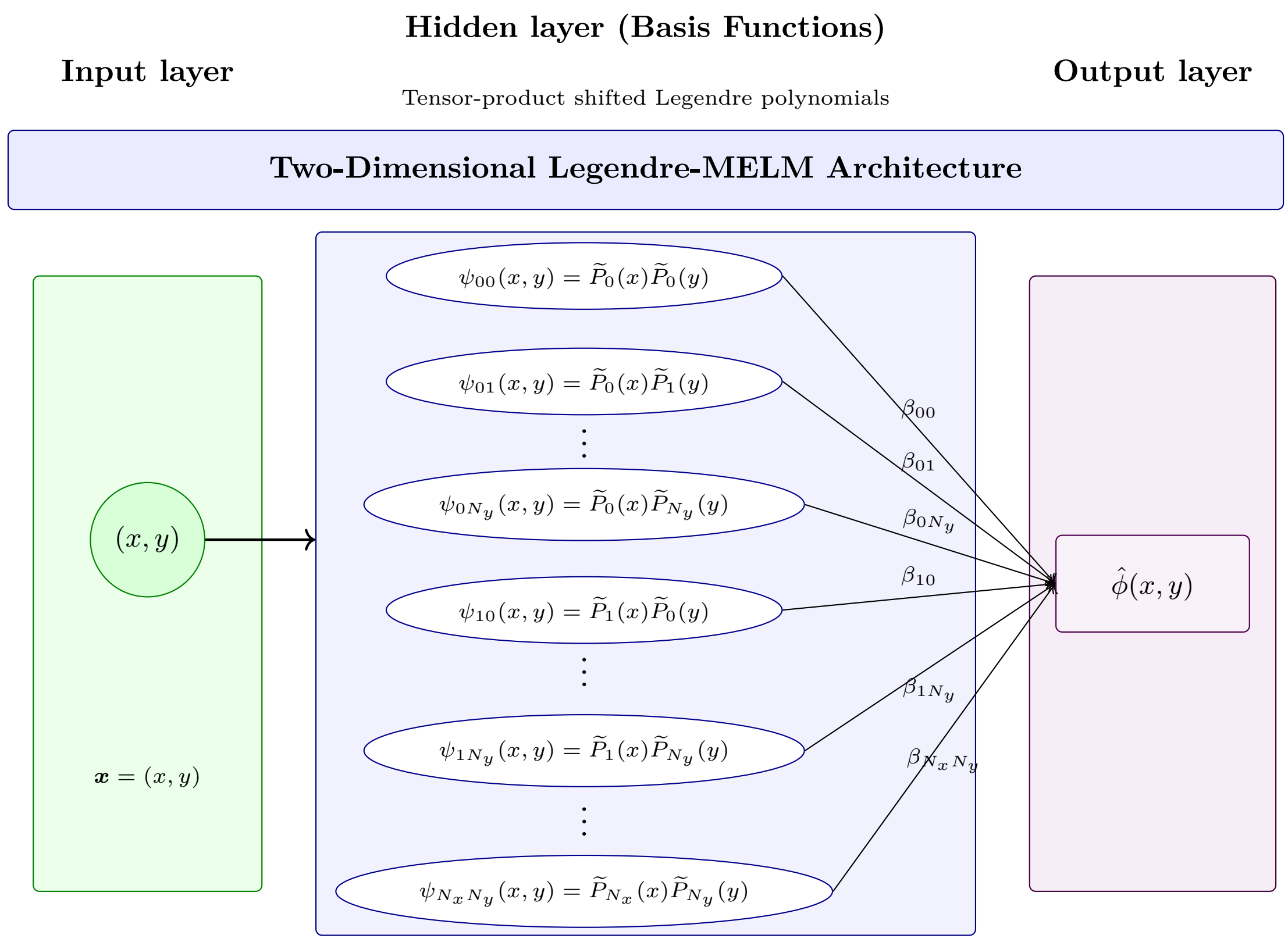

It is worth noting the close mathematical relationship between the proposed MELM framework and standard spectral collocation methods. Both approaches utilize a global polynomial expansion (in this case, shifted Legendre polynomials) and evaluate physical derivatives analytically via coordinate transformations. However, while traditional spectral collocation typically enforces a square linear system at predefined quadrature nodes, MELM views the representation through an extreme learning machine neural network lens. By framing orthogonal polynomials as deterministic hidden-layer activation functions, MELM provides a conceptual bridge to PINN architectures. This perspective opens the door to using PINNs with basis-driven activation functions that are significantly richer and more complex than standard choices (e.g., $\tanh$ or GELU). Furthermore, while the current MELM formulation maintains fixed linear weights without hidden-layer biases or intermediate deep layers, this neural framework creates a natural pathway for future research to selectively reintroduce trainable biases, nonlinear offsets, or multi-layer architectures to further enhance feature representation in complex regimes.

\subsubsection{Derivatives of the Approximate Solution}
\label{sec:derivatives}

The output of the MELM serves as the approximate solution of the DE. Its
derivatives with respect to the spatial variables are required for
evaluating the governing differential equation residuals.

\textbf{One-dimensional case:}
For $x \in [a,b]$,
\begin{equation*}
\hat{\phi}(x) = \sum_{n=0}^{N_x} \beta_n \tilde{P}_n(x),
\qquad
\tilde{P}_n(x)=P_n(\xi(x)),
\qquad
\xi(x)=\frac{2x-(a+b)}{b-a}.
\end{equation*}
The corresponding coordinate-scaling factor is
\begin{equation*}
\kappa_x=\frac{d\xi}{dx}=\frac{2}{b-a}.
\end{equation*}
Using the chain rule,
\begin{equation*}
\frac{d\hat{\phi}}{dx}
=
\frac{d\hat{\phi}}{d\xi}\frac{d\xi}{dx}
=
\kappa_x
\sum_{n=0}^{N_x}\beta_n P_n'(\xi),
\qquad
\frac{d^2\hat{\phi}}{dx^2}
=
\left(\frac{d\xi}{dx}\right)^2
\frac{d^2\hat{\phi}}{d\xi^2}
=
\kappa_x^2
\sum_{n=0}^{N_x}\beta_n P_n''(\xi).
\end{equation*}

\textbf{Two-dimensional case:}
For $(x,y)\in[a,b]\times[c,d]$,
\begin{equation*}
\hat{\phi}(x,y)=
\sum_{n=0}^{N_x}\sum_{m=0}^{N_y}
\beta_{nm}\tilde{P}_n(x)\tilde{P}_m(y),
\end{equation*}
with
\begin{equation*}
\xi(x)=\frac{2x-(a+b)}{b-a},
\qquad
\eta(y)=\frac{2y-(c+d)}{d-c}.
\end{equation*}
The corresponding coordinate-scaling factors are
\begin{equation*}
\kappa_x=\frac{d\xi}{dx}=\frac{2}{b-a},
\qquad
\kappa_y=\frac{d\eta}{dy}=\frac{2}{d-c}.
\end{equation*}
Using the chain rule, the first-order physical derivatives are
\begin{equation*}
\frac{\partial\hat{\phi}}{\partial x}
=
\frac{\partial\hat{\phi}}{\partial\xi}
\frac{\partial\xi}{\partial x}
=
\kappa_x
\sum_{n=0}^{N_x}\sum_{m=0}^{N_y}
\beta_{nm}
P_n'(\xi)P_m(\eta),
\end{equation*}
and
\begin{equation*}
\frac{\partial\hat{\phi}}{\partial y}
=
\frac{\partial\hat{\phi}}{\partial\eta}
\frac{\partial\eta}{\partial y}
=
\kappa_y
\sum_{n=0}^{N_x}\sum_{m=0}^{N_y}
\beta_{nm}
P_n(\xi)P_m'(\eta).
\end{equation*}

The second-order physical derivatives are given by
\begin{equation*}
\frac{\partial^2\hat{\phi}}{\partial x^2}
=
\left(\frac{\partial\xi}{\partial x}\right)^2
\frac{\partial^2\hat{\phi}}{\partial\xi^2}
=
\kappa_x^2
\sum_{n=0}^{N_x}\sum_{m=0}^{N_y}
\beta_{nm}
P_n''(\xi)P_m(\eta),
\end{equation*}
and
\begin{equation*}
\frac{\partial^2\hat{\phi}}{\partial y^2}
=
\left(\frac{\partial\eta}{\partial y}\right)^2
\frac{\partial^2\hat{\phi}}{\partial\eta^2}
=
\kappa_y^2
\sum_{n=0}^{N_x}\sum_{m=0}^{N_y}
\beta_{nm}
P_n(\xi)P_m''(\eta).
\end{equation*}
The mixed second-order derivative is
\begin{equation*}
\frac{\partial^2\hat{\phi}}{\partial x\,\partial y}
=
\frac{\partial\xi}{\partial x}
\frac{\partial\eta}{\partial y}
\frac{\partial^2\hat{\phi}}{\partial\xi\,\partial\eta}
=
\kappa_x\kappa_y
\sum_{n=0}^{N_x}\sum_{m=0}^{N_y}
\beta_{nm}
P_n'(\xi)P_m'(\eta).
\end{equation*}

\section{Implementation of the Neural Solvers}
\label{sec:sec4}

This section presents the numerical implementations of the PINN and the MELM framework. The PINN implementation is first described, followed by the formulation and implementation of the MELM framework.

{

\subsection{{PINNs Implementation}}

The PINN solvers are implemented following the formulation and training framework described in Section~\ref{subsec:PINNs}. Benchmark-specific choices of network architecture, collocation, optimization, loss weighting, and advanced training strategies are reported in the implementation tables of this section.

\subsubsection{Steady-state PINN}

For steady-state problems, the solution is represented by a spatial neural network $\hat{\phi}_{\theta}(\mathbf{x})$, where $\mathbf{x}$ denotes the spatial coordinates. Since the problem is time-independent, the training loss contains only the governing-equation and Dirichlet BC contributions defined in Eqs.~\eqref{eq:pinn_L_pde} and~\eqref{eq:pinn_L_bc}. Interior and boundary collocation points are sampled from their corresponding domains, and the spatial derivatives required by the convection--diffusion operator are evaluated using automatic differentiation.

The network architecture and training configuration are specified for each steady-state benchmark in Tables~\ref{tab:1_parameters}--\ref{tab:3_parameters}, including the network dimensions, activation function, collocation and batch sizes, optimization parameters, loss weights, and the advanced training strategies activated for each case.

\subsubsection{Continuous-Time PINN}

For transient problems, the CT-PINN employs a single neural network $\hat{\phi}_{\theta}(\mathbf{x},t)$ to represent the solution over the entire space--time domain. Time is therefore treated as an additional network input rather than being discretized during the neural-network solution process. The governing equation, Dirichlet BCs, and initial condition are enforced through the corresponding loss terms defined in Section~\ref{subsec:PINNs}, while the required spatial and temporal derivatives are evaluated using automatic differentiation.

For the transient benchmarks considered in this work, the network architecture, collocation and batch sizes, optimization parameters, loss weights, and activated training strategies are specified in Tables~\ref{tab:4_parameters} and~\ref{tab:5_parameters}. The advanced training strategies are selected according to the requirements of each benchmark and are not intrinsic to the CT formulation itself.

\subsubsection{Discrete‑Time PINN}
\label{subsec:pinn_dt}

The DT-PINN formulation follows the general time-stepping framework described in~\cite{Wang2023} and is adapted in the present work to the BDF temporal discretization introduced in Section~\ref{sec:temporal_discretization}. In contrast to the CT formulation, time is not provided as an input to the neural network. Instead, a spatial network $\hat{\phi}_{\theta}^{(n+1)}(\mathbf{x})$ is trained at each discrete time level, with the temporal evolution incorporated into the governing-equation residual through the solutions computed at the preceding time levels. The same BDF formulation and time-step size are used for the DT-PINN and DT-MELM formulations, allowing the two methods to be compared under the same temporal discretization while differing in their spatial solution strategy.

For the first time level, the BDF1 scheme is used, leading to the residual
\begin{equation}
R_{\theta}^{(1)}
=
\frac{\hat{\phi}_{\theta}^{(1)}-\phi^{(0)}}{\Delta t}
+\mathbf{a}\cdot\nabla\hat{\phi}_{\theta}^{(1)}
-k\nabla^{2}\hat{\phi}_{\theta}^{(1)}
-f^{(1)},
\label{eq:dt_pinn_bdf1}
\end{equation}
where $\phi^{(0)}$ denotes the prescribed initial condition. For subsequent time levels ($n\ge1$), the BDF2 scheme is employed, giving
\begin{equation}
R_{\theta}^{(n+1)}
=
\frac{
3\hat{\phi}_{\theta}^{(n+1)}
-4\phi^{(n)}
+\phi^{(n-1)}
}{2\Delta t}
+\mathbf{a}\cdot\nabla\hat{\phi}_{\theta}^{(n+1)}
-k\nabla^{2}\hat{\phi}_{\theta}^{(n+1)}
-f^{(n+1)}.
\label{eq:dt_pinn_bdf2}
\end{equation}

At each time level, the physics-informed loss is constructed from the PDE and Dirichlet BC contributions,
\begin{equation*}
\mathcal{L}^{(n+1)}(\theta)
=
\lambda_{\mathrm{PDE}}\,\mathcal{L}_{\mathrm{PDE}}^{(n+1)}(\theta)
+
\lambda_{\mathrm{BC}}\,\mathcal{L}_{\mathrm{BC}}^{(n+1)}(\theta),
\end{equation*}
where $\mathcal{L}_{\mathrm{PDE}}^{(n+1)}$ and $\mathcal{L}_{\mathrm{BC}}^{(n+1)}$ follow the definitions in Eqs.~\eqref{eq:pinn_L_pde} and~\eqref{eq:pinn_L_bc}, respectively, with the corresponding time-discrete residual defined by Eq.~\eqref{eq:dt_pinn_bdf1} for the first time level and by Eq.~\eqref{eq:dt_pinn_bdf2} thereafter. No separate initial-condition loss is introduced, since the initial condition enters the first BDF1 step directly and the subsequent temporal evolution is propagated through the previously computed solution fields. After convergence, the solution at time $t_{n+1}$ is stored and used in the residual construction at the next time level.

At each time level, the network parameters are initialized independently using the Glorot scheme. A warm-start strategy based on reusing the converged parameters from the preceding time level was also investigated, but no consistent improvement in convergence or final accuracy was observed; therefore, independent initialization was adopted for the reported results. The sequential training procedure is summarized in Algorithm~\ref{alg:dt_pinn}.

}

\begin{algorithm}[H] 
\caption{Sequential training of the DT-PINN using BDF2 with BDF1 initialization} 
\SetAlgoLined 
\KwIn{Domain $\Omega$, final time $T$, time step $\Delta t$, diffusivity $k$, convection velocity $\mathbf{a}$, source term $f$, BC $g$, and initial condition $\phi^{(0)}$.} 
\KwIn{Spatial neural network architecture $\hat{\phi}_{\theta}(\mathbf{x})$.} 
\KwResult{Sequentially computed solutions $\{\phi^{(n)}(\mathbf{x})\}_{n=1}^{N_t}$.} Set $N_t=T/\Delta t$\; 
\For{$n=0$ \KwTo $N_t-1$}{ Sample boundary points and interior collocation points\; Initialize $\theta^{(n+1)}$ independently using the Glorot scheme\; 
\eIf{$n=0$}{ Use the BDF1 residual in 
Eq.~\eqref{eq:dt_pinn_bdf1}\; }
{ Use the BDF2 residual in Eq.~\eqref{eq:dt_pinn_bdf2}\; } 
\For{$\text{epoch}=1$ 
\KwTo $\text{max\_epochs}$}{ Evaluate the selected time-discrete residual at the interior collocation points\; 
Evaluate the PDE and boundary losses\; Compute $\mathcal{L}^{(n+1)} = \lambda_{\mathrm{PDE}}\mathcal{L}_{\mathrm{PDE}}^{(n+1)} + \lambda_{\mathrm{BC}}\mathcal{L}_{\mathrm{BC}}^{(n+1)}$\; 
Update $\theta^{(n+1)}$ using the prescribed optimizer and learning-rate schedule\; } 
Store $\hat{\phi}_{\theta^{(n+1)}}(\mathbf{x})$ as $\phi^{(n+1)}(\mathbf{x})$\; } 
\label{alg:dt_pinn} 
\end{algorithm}

{

Thus, the DT-PINN formulation consists of a sequence of spatial PINN optimization problems coupled through the BDF time integration. This construction allows the same temporal discretization to be used across the DT-PINN and DT-MELM formulations, facilitating a direct comparison of their respective spatial solution strategies.

}

\subsection{MELM Implementation}

The MELM framework is implemented using both
CT and DT formulations. In the CT formulation, the temporal variable is treated as an additional independent variable, and the governing equation is enforced over the entire space--time domain by constructing a single global system of equations. In contrast, the DT formulation discretizes the temporal derivative using the BDF scheme described in Section~\ref{sec:temporal_discretization}, after which a sequence of spatial problems is solved successively at each time level. For steady-state problems, the temporal variable and its associated basis functions are omitted, reducing the formulation to the spatial domain only.

\subsubsection{Continuous-Time Formulation}
For the transient convection--diffusion equation
\eqref{eq:cde} subject to the boundary and initial conditions
\eqref{eq:icbc}, let

\begin{equation*}
\tilde{\mathbf{P}}(x,y,t)
=
\begin{bmatrix}
\tilde P_0(x)\tilde P_0(y)\tilde P_0(t) &
\tilde P_0(x)\tilde P_0(y)\tilde P_1(t) &
\cdots &
\tilde P_{N_x}(x)\tilde P_{N_y}(y)\tilde P_{N_\tau}(t)
\end{bmatrix},
\end{equation*}
which defines the tensor-product basis vector, where
$N_x$, $N_y$, and $N_\tau$ denote the maximum polynomial degrees in the
$x$-, $y$-, and $t$-directions, respectively. The basis functions are
arranged in lexicographic order. The MELM approximation is
given by

\begin{equation*}
\hat{\phi}(x,y,t)
=
\tilde{\mathbf{P}}(x,y,t)\,
\boldsymbol{\beta},
\end{equation*}
where
\begin{equation*}
\boldsymbol{\beta}
=
\begin{bmatrix}
\beta_{000} &
\beta_{001} &
\cdots &
\beta_{N_xN_yN_\tau}
\end{bmatrix}^{T}
\end{equation*}
is the vector of unknown output weights corresponding to the basis
functions.
Using the derivative expressions derived in
Section~\ref{sec:derivatives}, the required derivatives are
\begin{equation*}
\frac{\partial\hat{\phi}}{\partial t}
=
\tilde{\mathbf{P}}_t
\boldsymbol{\beta},
\qquad \frac{\partial\hat{\phi}}{\partial x}
=
\tilde{\mathbf{P}}_x
\boldsymbol{\beta}, \qquad \frac{\partial\hat{\phi}}{\partial y}
=
\tilde{\mathbf{P}}_y
\boldsymbol{\beta},
\end{equation*}
\begin{equation*}
\nabla^2\hat{\phi}
=
\left(
\tilde{\mathbf{P}}_{xx}
+
\tilde{\mathbf{P}}_{yy}
\right)
\boldsymbol{\beta}.
\end{equation*}
Here,
\(\tilde{\mathbf P}_x\),
\(\tilde{\mathbf P}_y\),
\(\tilde{\mathbf P}_t\),
\(\tilde{\mathbf P}_{xx}\),
and
\(\tilde{\mathbf P}_{yy}\)
denote the row vectors obtained by differentiating the corresponding basis functions with respect to the indicated variables.
Substituting the above approximation into
Eq.~\eqref{eq:cde} gives
\begin{equation*}
\left(
\tilde{\mathbf{P}}_t
+a_x\tilde{\mathbf{P}}_x
+a_y\tilde{\mathbf{P}}_y
-k\left(
\tilde{\mathbf{P}}_{xx}
+\tilde{\mathbf{P}}_{yy}
\right)
\right)
\boldsymbol{\beta}
=
f(x,y,t).
\end{equation*}
Enforcing the resulting equation at the interior collocation points yields
\begin{equation*}
\mathbf{A}_{\mathrm{int}}
\boldsymbol{\beta}
=
\mathbf{f}_{\mathrm{int}},
\end{equation*}
where $\mathbf{A}_{\mathrm{int}}$ is obtained by evaluating the differential operator acting on the basis functions at the interior collocation points, and $\mathbf{f}_{\mathrm{int}}$ is the vector of source-term values evaluated at the corresponding collocation points.
Enforcing the prescribed boundary and initial conditions at the corresponding collocation points yields
\begin{equation*}
\mathbf{A}_{\mathrm{bc}}\boldsymbol{\beta}=\mathbf{g}
\quad \text{and} \quad
\mathbf{A}_{\mathrm{ic}}\boldsymbol{\beta}=\boldsymbol{\phi}_0,
\end{equation*}
where $\mathbf{A}_{\mathrm{bc}}$ and $\mathbf{A}_{\mathrm{ic}}$ are the coefficient matrices corresponding to the boundary and initial conditions, respectively, while $\mathbf{g}$ and $\boldsymbol{\phi}_0$ are the vectors of prescribed boundary and initial condition values evaluated at the corresponding collocation points.
Introducing the positive weighting parameters $\lambda_1$ and $\lambda_2$ for the boundary and initial conditions, respectively, the interior, boundary, and initial condition systems are assembled into the weighted global linear system
\begin{equation}
\mathbf{H}\boldsymbol{\beta}
=
\mathbf{T},
\label{eq:global_system}
\end{equation}
where
\begin{equation*}
\mathbf{H}=
\begin{bmatrix}
\mathbf{A}_{\mathrm{int}}\\[2mm]
\sqrt{\lambda_1}\,\mathbf{A}_{\mathrm{bc}}\\[2mm]
\sqrt{\lambda_2}\,\mathbf{A}_{\mathrm{ic}}
\end{bmatrix},
\qquad
\mathbf{T}=
\begin{bmatrix}
\mathbf{f}_{\mathrm{int}}\\[2mm]
\sqrt{\lambda_1}\,\mathbf{g}\\[2mm]
\sqrt{\lambda_2}\,\boldsymbol{\phi}_0
\end{bmatrix}.
\end{equation*}
Since Eq.~\eqref{eq:global_system} is generally overdetermined,
the output weight vector $\boldsymbol{\beta}$ is obtained by solving
the least-squares problem
\begin{equation*}
\min_{\boldsymbol{\beta}}
\left\|
\mathbf{H}\boldsymbol{\beta}
-
\mathbf{T}
\right\|_2^2,
\end{equation*}
using a QR-based least-squares solution. The computed output weight
vector $\boldsymbol{\beta}$ is subsequently substituted into the MELM
approximation to evaluate the solution at any point in the space--time
domain.

\subsubsection{Discrete-Time Formulation}

Using the BDF time discretization introduced in
Section~\ref{sec:temporal_discretization}, the transient
convection--diffusion equation is transformed into a sequence of
time-discrete boundary value problems. The first time step is computed
using the BDF1 scheme, while the BDF2 scheme is employed for all
subsequent time steps. At time level $t_{n+1}$, the MELM approximation is given by
\begin{equation*}
\hat{\phi}^{(n+1)}(x,y)
=
\tilde{\mathbf{P}}(x,y)
\boldsymbol{\beta}^{(n+1)},
\qquad n=0,1,2,\ldots,N_t-1,
\end{equation*}
where $N_t=T/\Delta t$ denotes the total number of time steps, with $T$
being the final time and $\Delta t$ the time-step size. Here,
\begin{equation*}
\tilde{\mathbf{P}}(x,y)=
\begin{bmatrix}
\tilde P_0(x)\tilde P_0(y) &
\tilde P_0(x)\tilde P_1(y) &
\cdots &
\tilde P_{N_x}(x)\tilde P_{N_y}(y)
\end{bmatrix},
\end{equation*}
which defines the row vector of tensor-product basis functions arranged in lexicographic order, where $N_x$ and $N_y$ denote the maximum polynomial degrees in the $x$- and $y$-directions, respectively. The vector
\(
\boldsymbol{\beta}^{(n+1)}
\)
is the column vector of unknown output weights at time level $t_{n+1}$.
Using the expressions derived in Section~\ref{sec:derivatives}, the corresponding spatial derivatives are given by

\[
\frac{\partial\hat{\phi}^{(n+1)}}{\partial x}
=
\tilde{\mathbf{P}}_x
\boldsymbol{\beta}^{(n+1)},
\qquad 
\frac{\partial\hat{\phi}^{(n+1)}}{\partial y}
=
\tilde{\mathbf{P}}_y
\boldsymbol{\beta}^{(n+1)},
\]
and
\[
\nabla^2\hat{\phi}^{(n+1)}
=
\left(
\tilde{\mathbf{P}}_{xx}
+
\tilde{\mathbf{P}}_{yy}
\right)
\boldsymbol{\beta}^{(n+1)}.
\]
The prescribed initial condition is denoted by
\(
\hat{\phi}^{(0)}=\phi^{(0)}.
\)
At the first time step ($n=0$), substituting the BDF1 approximation
into Eq.~\eqref{eq:cde}, multiplying by $\Delta t$, and rearranging
gives
\begin{equation*}
\hat{\phi}^{(1)}
+\Delta t\,a_x\frac{\partial\hat{\phi}^{(1)}}{\partial x}
+\Delta t\,a_y\frac{\partial\hat{\phi}^{(1)}}{\partial y}
-\Delta t\,k\nabla^2\hat{\phi}^{(1)}
=
\phi^{(0)}
+\Delta t\,f^{(1)}.
\end{equation*}
Substituting the MELM approximation into the resulting equation yields
\begin{equation*}
\left[
\tilde{\mathbf{P}}
+\Delta t\,a_x\tilde{\mathbf{P}}_x
+\Delta t\,a_y\tilde{\mathbf{P}}_y
-\Delta t\,k
\left(
\tilde{\mathbf{P}}_{xx}
+
\tilde{\mathbf{P}}_{yy}
\right)
\right]
\boldsymbol{\beta}^{(1)}
=
\phi^{(0)}
+\Delta t\,f^{(1)}.
\end{equation*}
Enforcing the resulting equation at the interior collocation points together with the prescribed BCs at the boundary collocation points yields the linear system
\begin{equation}
\label{eq:melm_bdf1}
\mathbf{A}_{\mathrm{BDF1}}
\boldsymbol{\beta}^{(1)}
=
\mathbf{b}_{\mathrm{BDF1}},
\end{equation}
where
\(\mathbf{A}_{\mathrm{BDF1}}\)
is obtained by evaluating the governing equation and the prescribed BCs at the corresponding collocation points, while
\(\mathbf b_{\mathrm{BDF1}}\)
contains the contributions from the initial condition, the source term, and the prescribed BCs. Since the BDF1 system is overdetermined, the unknown output weight vector $\boldsymbol{\beta}^{(1)}$ is obtained using a QR-based least-squares solution:
\begin{equation*}
\min_{\boldsymbol{\beta}^{(1)}}
\left\|
\mathbf{A}_{\mathrm{BDF1}}
\boldsymbol{\beta}^{(1)}
-
\mathbf{b}_{\mathrm{BDF1}}
\right\|_2^2.
\end{equation*}
For subsequent time steps ($n=1,2,\ldots,N_t-1$), substituting the
BDF2 approximation into Eq.~\eqref{eq:cde}, multiplying by
$2\Delta t$, and rearranging gives
\begin{equation*}
3\hat{\phi}^{(n+1)}
+2\Delta t\,a_x\frac{\partial\hat{\phi}^{(n+1)}}{\partial x}
+2\Delta t\,a_y\frac{\partial\hat{\phi}^{(n+1)}}{\partial y}
-2\Delta t\,k\nabla^2\hat{\phi}^{(n+1)}
=
4\hat{\phi}^{(n)}
-\hat{\phi}^{(n-1)}
+2\Delta t\,f^{(n+1)}.
\end{equation*}
Substituting the MELM approximation into the resulting equation yields
\begin{equation*}
\left[
3\tilde{\mathbf{P}}
+2\Delta t\,a_x\tilde{\mathbf{P}}_x
+2\Delta t\,a_y\tilde{\mathbf{P}}_y
-2\Delta t\,k
\left(
\tilde{\mathbf{P}}_{xx}
+
\tilde{\mathbf{P}}_{yy}
\right)
\right]
\boldsymbol{\beta}^{(n+1)}
=
4\hat{\phi}^{(n)}-\hat{\phi}^{(n-1)}
+2\Delta t\,f^{(n+1)}.
\end{equation*}
Enforcing the resulting equation at the interior collocation points together with the prescribed BCs at the boundary collocation points yields the linear system
\begin{equation}
\label{eq:melm_bdf2}
\mathbf{A}_{\mathrm{BDF2}}
\boldsymbol{\beta}^{(n+1)}
=
\mathbf{b}_{\mathrm{BDF2}},
\qquad
n=1,2,\ldots,N_t-1,
\end{equation}
where $\mathbf{A}_{\mathrm{BDF2}}$ is obtained by evaluating the governing equation and the prescribed BCs at the corresponding collocation points, while $\mathbf{b}_{\mathrm{BDF2}}$ contains the contributions from the solutions at the two previous time levels, $\hat{\phi}^{(n)}$ and $\hat{\phi}^{(n-1)}$, the source term, and the prescribed BCs. Since the BDF2 system is overdetermined, the unknown output weight vector $\boldsymbol{\beta}^{(n+1)}$ is obtained by solving the least-squares problem
\begin{equation*}
\min_{\boldsymbol{\beta}^{(n+1)}}
\left\|
\mathbf{A}_{\mathrm{BDF2}}
\boldsymbol{\beta}^{(n+1)}
-
\mathbf{b}_{\mathrm{BDF2}}
\right\|_2^2,
\qquad
n=1,2,\ldots,N_t-1.
\end{equation*}
Since the coefficient matrix
$\mathbf{A}_{\mathrm{BDF2}}$
is identical for all BDF2 time steps, its QR factorization is computed
only once. At each subsequent time step, only the right-hand side vector
$\mathbf{b}_{\mathrm{BDF2}}$ is updated, while the precomputed QR
factors are reused to solve the corresponding least-squares problem.
This avoids repeated matrix factorizations and significantly reduces the
overall computational cost. The computed output weight vector
$\boldsymbol{\beta}^{(n+1)}$
defines the MELM approximation at time level
$t_{n+1}$.

In addition to the above overdetermined formulation, a square-system
formulation using LU decomposition is also considered. In this formulation, the number of collocation equations is chosen to be equal to the number of basis functions, resulting in a square coefficient matrix. For the BDF1 step, the linear system is solved using a direct dense solver to obtain $\boldsymbol{\beta}^{(1)}$. For the subsequent BDF2 time steps, the LU factorization of $\mathbf{A}_{\mathrm{BDF2}}$ is computed once and reused, while only the right-hand side is updated.

The choice between the square-system LU formulation and the
overdetermined QR formulation is guided by the numerical requirements
of the problem, particularly the conditioning and numerical stability
of the resulting linear systems and the accuracy of the computed
solution. Accordingly, the formulation providing the most satisfactory
numerical performance is adopted for each problem. The sequential solution procedure for the DT-MELM is summarized in Algorithm~\ref{alg:dt_melm}.

\begin{algorithm}[H]
\caption{Sequential solution of the DT-MELM using BDF2 with BDF1 initialization}
\SetAlgoLined
\KwIn{Domain $\Omega$, final time $T$, time step $\Delta t$, diffusivity $k$, convection velocity $\mathbf{a}$, source term $f$, BC $g$, and initial condition $\phi^{(0)}$}
\KwIn{Basis functions.}

\KwResult{Sequentially computed solutions $\{\hat{\phi}^{(n)}(\mathbf{x})\}_{n=1}^{N_t}$.}
Set $N_t=T/\Delta t$\;
Construct $\tilde{\mathbf{P}}(\mathbf{x})$ and the required spatial
derivative matrices\;
Construct $\mathbf{A}_{\mathrm{BDF2}}$ given by Eq.~\eqref{eq:melm_bdf2}
and factorize it using LU decomposition if the system is square or QR
decomposition if the system is overdetermined\;
\For{$n=0$ \KwTo $N_t-1$}{

\eIf{$n=0$}{
Construct the BDF1 linear system given by Eq.~\eqref{eq:melm_bdf1}\;
Solve for $\boldsymbol{\beta}^{(1)}$ using a direct dense solver if the
system is square, or a QR-based least-squares solution if the system is overdetermined\;
}{

Update $\mathbf{b}_{\mathrm{BDF2}}$ using the solutions
$\hat{\phi}^{(n)}$ and $\hat{\phi}^{(n-1)}$ at the two preceding time levels
and solve for $\boldsymbol{\beta}^{(n+1)}$ using the precomputed factorization\;
}

Store
$\hat{\phi}^{(n+1)}(\mathbf{x})
=\tilde{\mathbf{P}}(\mathbf{x})\boldsymbol{\beta}^{(n+1)}$\;
}
\label{alg:dt_melm}
\end{algorithm}

{ \subsection{Computational Environment and Implementation Details}

The computational environments used for the PINN, FEM, and MELM simulations are described in this subsection, including the software frameworks, numerical solvers, and computational resources. The main hardware configurations are summarized in Table~\ref{tab:comp_res}.

\subsubsection{PINN implementation}

The PINN models are implemented in JAX~\cite{jax2018github} using the JAXPI framework~\cite{Wang2023}. The source code is publicly available through the JAXPI repository. The experiments employing the SOAP optimizer use a customized version of the \texttt{pirate} branch of JAXPI. Although this branch also provides PirateNets-based architectures, only the standard MLP and Modified MLP configurations are used in the present study.

All PINN simulations are performed using a single GPU. The one-dimensional problems are executed on an NVIDIA GeForce 940MX, while the two-dimensional problems are executed on an NVIDIA T1000. The corresponding computational resources are summarized in
Table~\ref{tab:comp_res}.

\subsubsection{Finite-element reference solver}

The FEM reference solutions are computed using a finite-element solver. For the linear convection--diffusion problems considered here, the finite-element discretization leads to a global linear system for the scalar unknown $\phi$. The resulting systems are  solved using the BiCGSTAB iterative method~\cite{Vorst1992}, implemented through the PETSc parallel linear algebra library~\cite{petsc-efficient}, with an additive Schwarz method (ASM) preconditioner employing local incomplete LU (ILU) factorizations. The relative tolerance of the Krylov solver is set to $10^{-8}$ for steady problems and $10^{-15}$ for transient problems. The FEM simulations are executed using parallel CPU resources.

\subsubsection{MELM implementation}

The MELM computations are performed on a local workstation and, for selected problems, on the Aqua HPC system. All MELM computations are performed using CPU resources. The corresponding processors and memory allocations are summarized in Table~\ref{tab:comp_res}. The source codes associated with the present study are available on GitHub at \href{https://github.com/GadadeNavnath/Neural-solvers-for-CDEs} {this repository}.

\begin{table}[H]
\centering
\caption{Computational resources used for the FEM, PINN,
and MELM simulations.}
\label{tab:comp_res}
\renewcommand{\arraystretch}{1.3}
\setlength{\tabcolsep}{8pt}

\begin{tabular}{
p{0.23\textwidth}
p{0.36\textwidth}
p{0.30\textwidth}
}
\toprule
\textbf{Method} &
\textbf{Processor / GPU} &
\textbf{Memory} \\
\midrule

FEM &
Intel Core i9-13900 &
32~GB RAM \\

PINN (1D) &
NVIDIA GeForce 940MX &
2~GB VRAM \\

PINN (2D) &
NVIDIA T1000 &
4~GB VRAM \\

MELM (Local) &
Intel Core i5-7500 &
8~GB RAM \\

MELM (Aqua HPC) &
4 CPU cores &
16~GB RAM (Problems~\ref{case_2}
and~\ref{case_3})
\newline
190~GB RAM (Problem~\ref{case_5}, CT) \\

\bottomrule
\end{tabular}
\end{table}

For all methods, the reported computational times correspond to the total wall-clock time required to obtain the numerical solution under the respective computational configurations. Since the simulations were performed using different hardware resources and computational architectures, these timings should be interpreted as indicative measures of computational cost rather than hardware-normalized performance comparisons.

}

\section{Numerical Results and Discussion}
\label{sec:sec5}

This section presents the numerical results for the benchmark convection--diffusion problems considered in this study. The results include the collocation-point distributions, numerical solution comparisons, error analysis, and computational performance of the PINN and MELM approaches. Where no exact solution is available, the converged FEM solution is used as the numerical reference for assessing the PINN and MELM approximations.

The collocation strategies used by the different numerical approaches follow their respective formulations described in the preceding sections. For the PINN, collocation points are sampled randomly within the corresponding spatial, boundary, and temporal domains. For the CT-MELM formulation, Chebyshev points are used along the spatial directions $(x,y)$ and the temporal direction $t$ for transient problems, whereas for steady problems they are used only in the spatial directions. For the DT-MELM formulation, they are used only in the spatial directions, with the temporal evolution discretized using the uniform time step $\Delta t$ introduced in Section~\ref{sec:sec4}. The specific numbers of collocation points and other problem-dependent parameters are reported in the corresponding implementation tables.

The pointwise absolute error is defined as
\begin{equation*}
\text{Absolute error}
=
\left|
\hat{\phi}(\mathbf{x},t)
-
\phi(\mathbf{x},t)
\right|,
\end{equation*}
where $\hat{\phi}$ denotes the numerical approximation obtained with the
PINN or MELM, and $\phi$ denotes the corresponding exact or FEM reference
solution.

The relative $L_2$ error is defined as
\begin{equation*}
\text{Relative } L_2\text{ error}
=
\frac{
\left\|
\hat{\phi}-\phi
\right\|_2
}{
\left\|
\phi
\right\|_2
}.
\end{equation*}

For the transient problems, the same temporal sampling is also used when comparing the CT and DT formulations. Specifically, for a given problem, the CT solution is evaluated at the discrete time levels defined by the time step $\Delta t$ used in the corresponding DT formulation. Thus, for Problem~\ref{case_4}, the CT solution is sampled at the time levels associated with $\Delta t=0.01$, while for Problem~\ref{case_5} it is sampled using $\Delta t=0.001$. The resulting error measures are therefore evaluated over a common space--time set for the corresponding CT and DT formulations, with the FEM mesh defining the spatial evaluation points and the DT time levels defining the temporal evaluation points.

\begin{problem}[One-Dimensional Steady Convection–Diffusion Problem]\label{case_1}
Consider the one-dimensional steady convection--diffusion equation
\begin{equation*}
a\frac{d\phi}{dx} - k\frac{d^{2}\phi}{dx^{2}}=0, \qquad x\in(0,1),   
\end{equation*}
subject to the Dirichlet BCs  \(
\phi(0)=1,
\:
\phi(1)=0.
\)
The problem parameters are
\(
a=1,
\:
k=2\times10^{-4},
\:
L=1,
\)
yielding the Péclet number
\begin{equation*}
Pe=\frac{aL}{k}=5{,}000.
\end{equation*}
The corresponding exact solution is
\begin{equation*}
\phi(x)
=
\frac{1-e^{Pe(x-1)}}
     {1-e^{-Pe}}.
\end{equation*}
  
\end{problem}

The FEM reference solution is computed using a one-dimensional mesh with $5{,}000$ elements, with mesh refinement concentrated near the right boundary where the sharp boundary layer develops. The element sizes decrease from approximately $5.05\times10^{-4}$ in the coarser region to $3.42\times10^{-6}$ within the refined region, providing sufficient spatial resolution of the high-$Pe$ boundary layer.

Table~\ref{tab:1_parameters} summarizes the implementation parameters used for the PINN and MELM methods. For the MELM, $375$ basis functions are used (corresponding to polynomial degrees $0$--$374$). Figure~\ref{fig:1} compares the exact solution with the numerical solutions obtained using the PINN and MELM approaches. This high-$Pe$ benchmark is primarily used to assess the ability of the two numerical representations to resolve the sharp boundary layer generated near the outflow boundary. As shown in Fig.~\ref{fig:1_comp}, both methods reproduce the exact solution over the computational domain, while the enlarged view in Fig.~\ref{fig:1_comp_zoom} shows that the boundary layer near $x=1$ is captured without visible nonphysical oscillations. Figures~\ref{fig:1_error_pinn} and \ref{fig:1_error_melm} present the corresponding pointwise absolute errors, showing that the largest deviations are associated with the region of strongest solution variation. The relative $L_2$ errors and computational times are reported in Table~\ref{tab:comparison}.

\begin{table}[H]
\centering
\caption{Implementation parameters used for PINN and MELM for Problem~\ref{case_1}.} 

\label{tab:1_parameters}

\begin{minipage}[t]{0.48\textwidth}
\centering
\textbf{PINN}\\[2mm]
\begin{tabular}{ll}
\toprule
\textbf{Parameter} & \textbf{Value} \\
\midrule
\multicolumn{2}{l}{\textit{Architecture}} \\
Network type & Modified MLP \\
Layers & 6 \\
Neurons & 128 \\
Activation & GELU \\
Fourier feature scale & 10.0 \\
RWF & $\mu=0.5,\ \sigma=0.1$ \\
\midrule
\multicolumn{2}{l}{\textit{Collocation and training}} \\
Batch size (residual) & 1{,}024 \\
Learning rate & 0.001 \\
Training iterations & 200{,}000 \\
$\lambda_i$ (BC, PDE) & $(10.0,\ 1.0)$ \\
PRNG seed & 123 \\
\midrule
\multicolumn{2}{l}{\textit{Optimization}} \\
Optimizer & SOAP \\
$\beta_1$ & 0.95 \\
LR schedule & Step decay \\
LR decay rate & 0.9 \\
LR decay steps & 2{,}000 \\
\midrule
\multicolumn{2}{l}{\textit{Loss weighting}} \\
Weighting scheme & NTK \\
\bottomrule
\end{tabular}
\end{minipage}
\hfill
\begin{minipage}[t]{0.45\textwidth}
\centering
\textbf{MELM}\\[2mm]
\begin{tabular}{ll}
\toprule
\textbf{Parameter} & \textbf{Value} \\
\midrule
\multicolumn{2}{l}{\textit{Architecture}} \\
Network type & MELM (Legendre) \\
Layers & 1 \\
No. of Basis functions & 375 \\
Basis & Shifted Legendre \\
\midrule
\multicolumn{2}{l}{\textit{Collocation and solution}} \\
Collocation points & 502 \\
$\lambda_i$ (BC) & $ 1.0 $ \\
\midrule
\multicolumn{2}{l}{\textit{Linear solver}} \\
System size & $502\times375$ \\
Method & QR-based least-squares \\
\bottomrule
\end{tabular}
\end{minipage}

\end{table}

\begin{figure}[H]
\centering

\begin{subfigure}{0.40\textwidth}
    \centering
    \includegraphics[width=\linewidth]{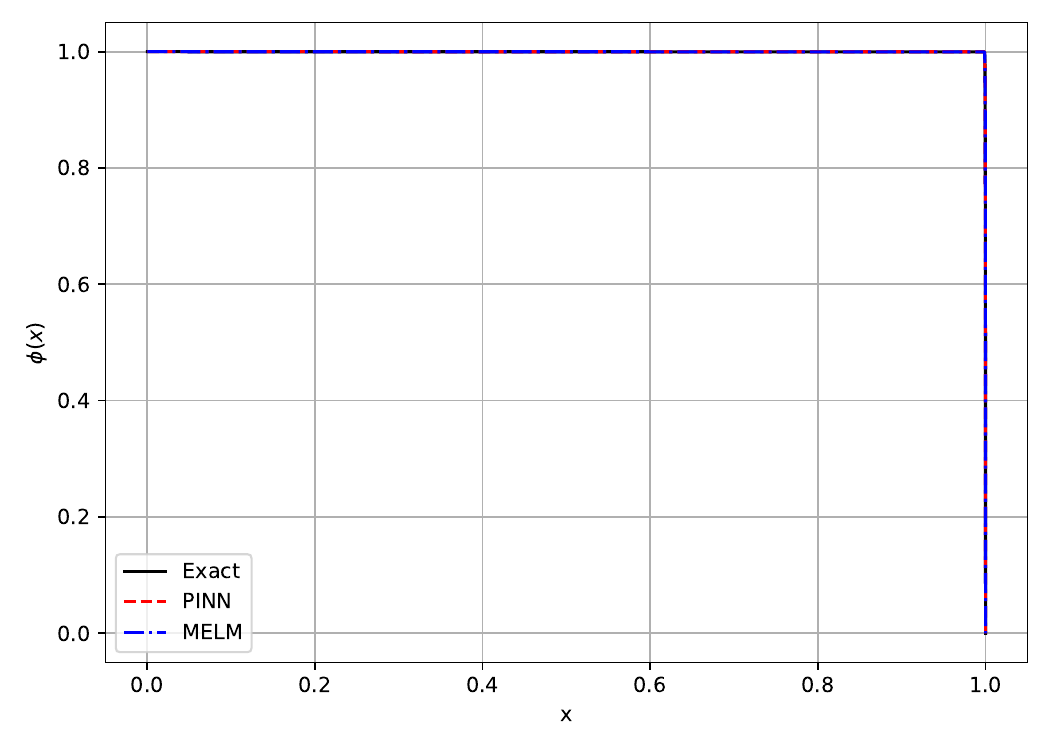}
    \caption{Comparison of the solutions.}
    \label{fig:1_comp}
\end{subfigure}
\hfill
\begin{subfigure}{0.40\textwidth}
    \centering
    \includegraphics[width=\linewidth]{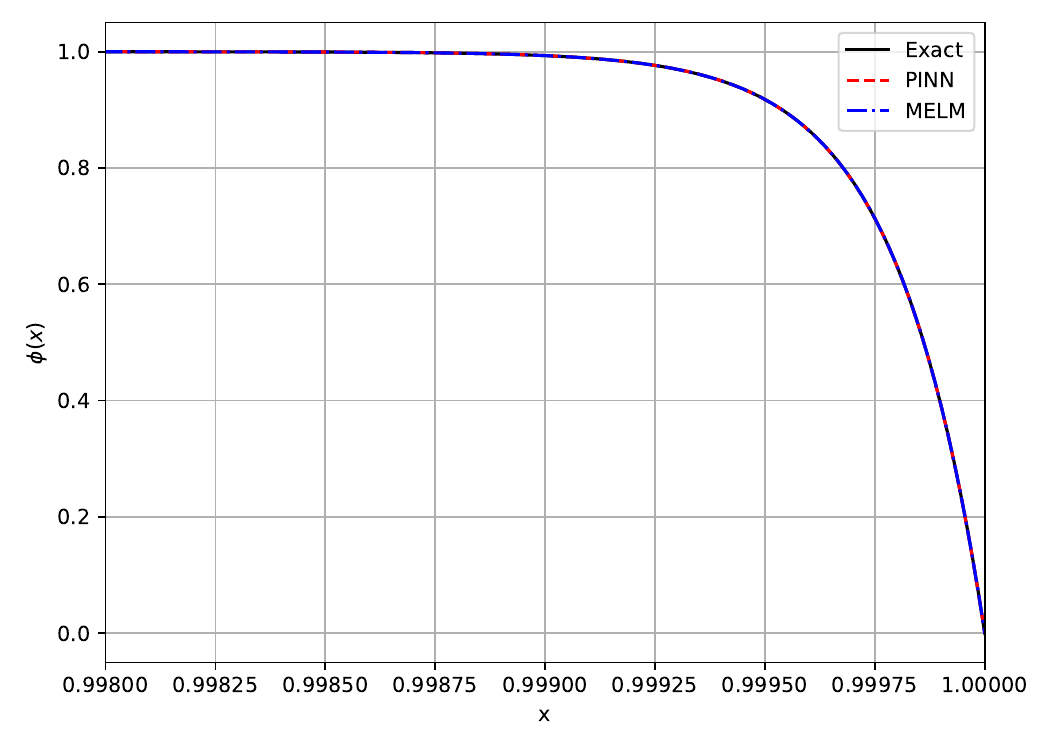}
    \caption{Zoomed-in view near the boundary layer.}
    \label{fig:1_comp_zoom}
\end{subfigure}

\vspace{0.35cm}

\begin{subfigure}{0.40\textwidth}
    \centering
    \includegraphics[width=\linewidth]{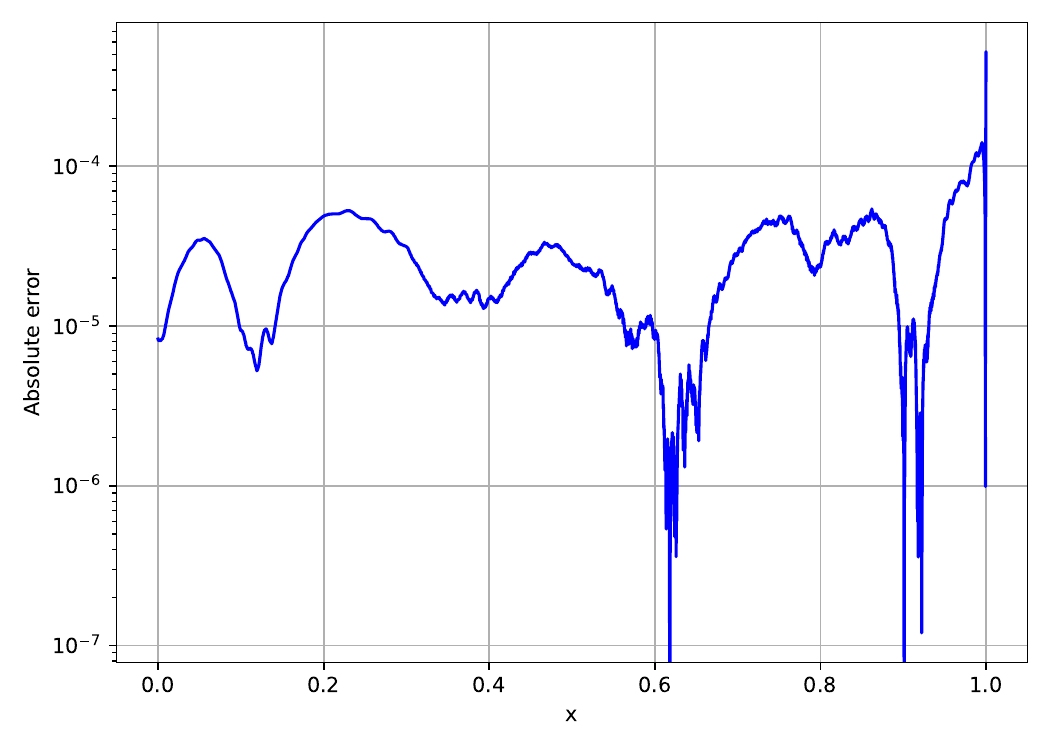}
    \caption{PINN}
    \label{fig:1_error_pinn}
\end{subfigure}
\hfill
\begin{subfigure}{0.40\textwidth}
    \centering
    \includegraphics[width=\linewidth]{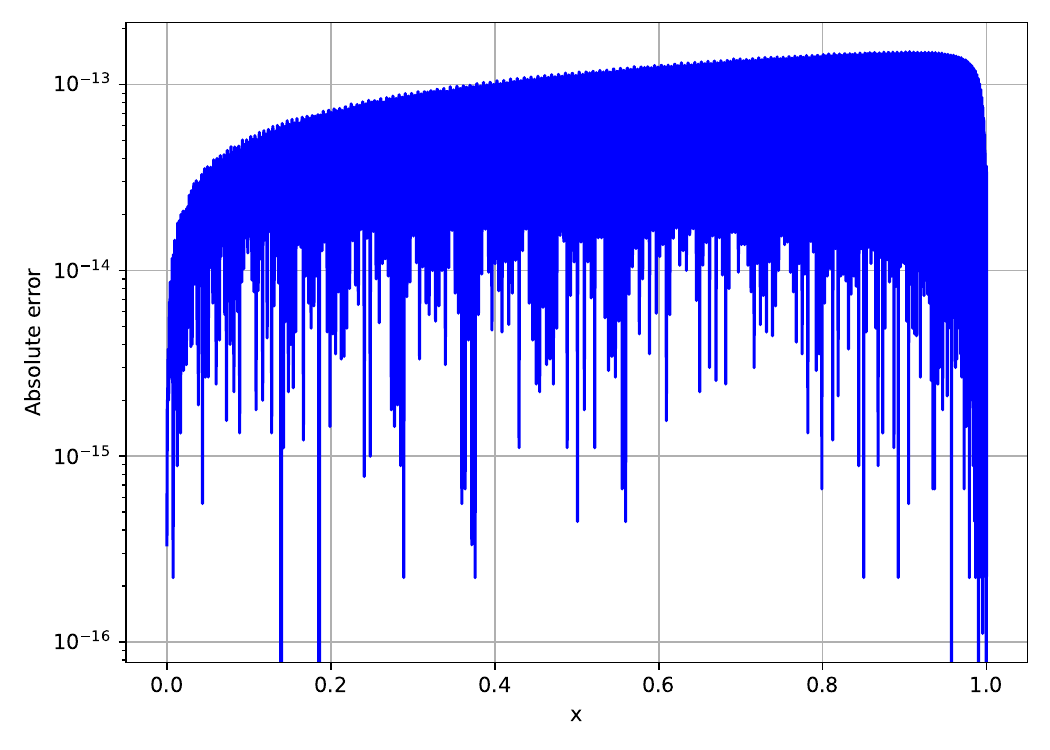}
    \caption{MELM}
    \label{fig:1_error_melm}
\end{subfigure}

\caption{Comparison of the exact, PINN, and MELM solutions for Problem~\ref{case_1}. The first row shows the solution comparison and a zoomed-in view of the boundary layer, while the second row presents the corresponding absolute errors.}
\label{fig:1}

\end{figure}

\begin{problem}[Two-Dimensional Steady Problem (Continuous BCs)]\label{case_2}
    Consider the steady two-dimensional convection--diffusion equation
\begin{equation*}
a_x\frac{\partial\phi}{\partial x}
+
a_y\frac{\partial\phi}{\partial y}
-
k\left(
\frac{\partial^{2}\phi}{\partial x^{2}}
+
\frac{\partial^{2}\phi}{\partial y^{2}}
\right)
=
5,
\qquad
(x,y)\in\Omega:=\left(-\frac{1}{2},\frac{1}{2}\right)
\times \left(-\frac{1}{2},\frac{1}{2}\right),
\end{equation*}
subject to the homogeneous Dirichlet BC
\begin{equation*}
\phi(x,y)=0
\quad
\text{on } \partial\Omega.
\end{equation*}
The problem parameters are
\begin{equation*}
a_x=a_y=\frac{1}{\sqrt{2}},\qquad
U=1,\qquad
k=0.001,\qquad
L=1,
\end{equation*}
yielding the Péclet number
\begin{equation*}
Pe=\frac{UL}{k}=1{,}000.
\end{equation*}
\end{problem}

The FEM reference solution is computed on a structured quadrilateral mesh containing $90{,}000$ elements. The same mesh bias factor is applied along the $x$- and $y$-directions, producing a progressive concentration of elements toward the corner where the sharp solution layer develops. The corresponding element side lengths vary approximately from $0.00583$ in the coarser regions to $0.000626$ near the corner, as shown in Fig.~\ref{fig:P2_2dmesh}. This localized refinement provides increased spatial resolution in the region of the strongest solution gradients for the FEM results.

\begin{figure}[H]
    \centering
    \includegraphics[width=0.5\linewidth]{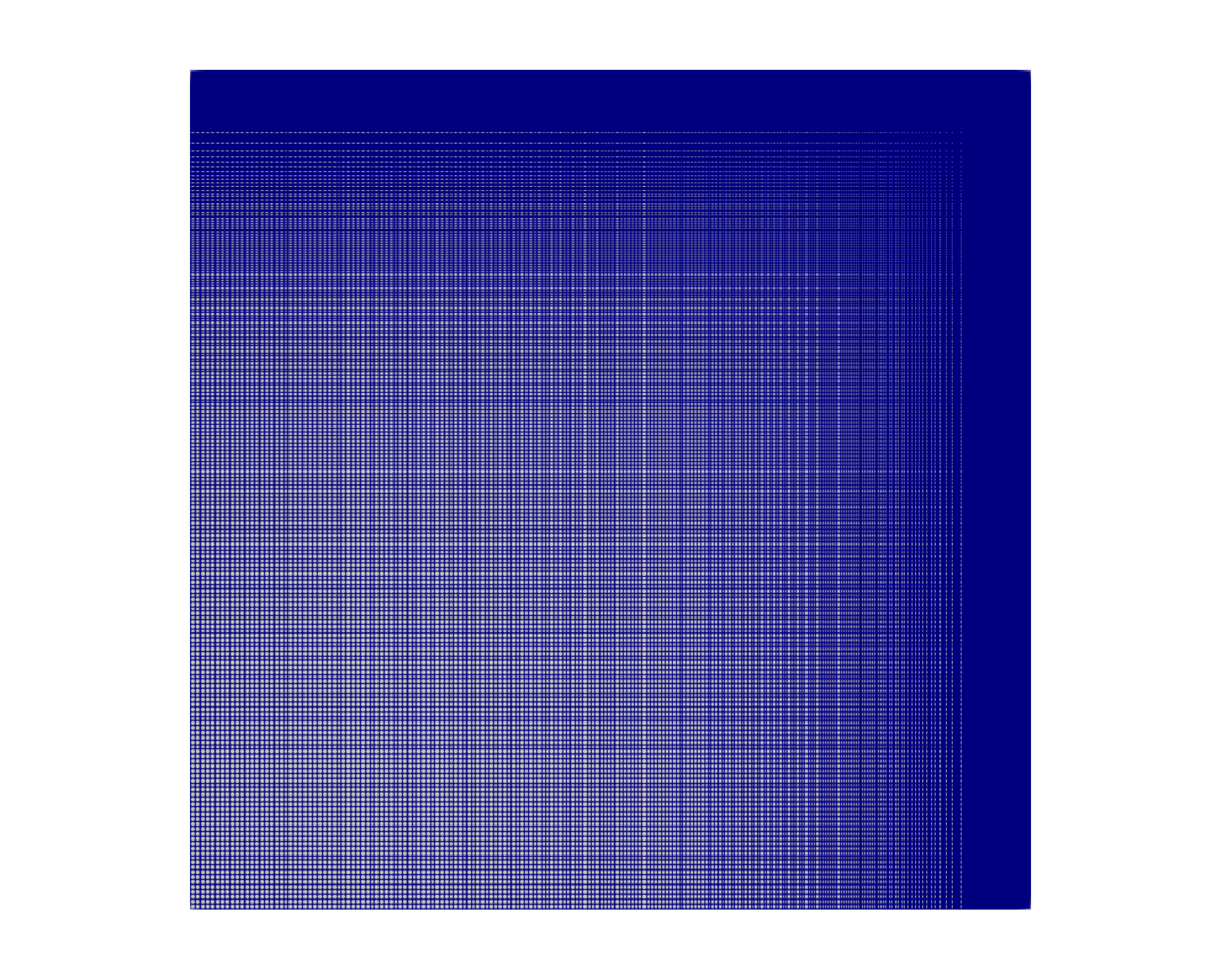}
    \caption{Structured FEM mesh used to compute the reference solution for Problem~\ref{case_2}. Local mesh refinement is introduced near the corner associated with the sharp solution layer.}
    \label{fig:P2_2dmesh}
\end{figure}

Table~\ref{tab:2_parameters} summarizes the implementation parameters used for the PINN and MELM methods in Problem~\ref{case_2}. For the MELM, $110$ basis functions are used along each of the $x$- and $y$-directions, resulting in $110\times110=12{,}100$ tensor-product basis functions. The numerical solutions obtained using the FEM, PINN, and MELM approaches are compared in Fig.~\ref{fig:2_surface}, while the corresponding contour representations are presented in Fig.~\ref{fig:2_contour}. Fig.~\ref{fig:2_error} illustrates the absolute error surfaces of the PINN and MELM solutions with respect to the FEM reference solution, highlighting the spatial distribution of the approximation errors. To provide a detailed local comparison, Fig.~\ref{fig:2_cross_section} presents solution profiles along selected cross-sections of the computational domain, and Fig.~\ref{fig:2_boundary_profile} compares the solution profiles along the domain boundaries.

The numerical results show close agreement between the PINN and MELM approximations and the FEM reference solution, as observed in the contour, cross-sectional, and boundary-profile comparisons. The error surfaces indicate that the largest deviations are localized in regions with steep solution gradients. The relative $L_2$ errors and computational times are reported in Table~\ref{tab:comparison}.

\begin{table}[H]
\centering
\caption{Implementation parameters used for PINN and MELM for Problem~\ref{case_2}.}
\label{tab:2_parameters}

\begin{minipage}[t]{0.48\textwidth}
\centering
\textbf{PINN}\\[2mm]
\begin{tabular}{ll}
\toprule
\textbf{Parameter} & \textbf{Value} \\
\midrule
\multicolumn{2}{l}{\textit{Architecture}} \\
Network type & Modified MLP \\
Layers & 6 \\
Neurons & 128 \\
Activation & tanh \\
Fourier feature scale & 10.0 \\
RWF & $\mu=0.5,\ \sigma=0.1$ \\
\midrule
\multicolumn{2}{l}{\textit{Collocation and training}} \\
Batch size (residual) & 256 \\
Batch size (boundary) & 256 \\
Learning rate & $10^{-3}$ \\
Training iterations & 200{,}000 \\
$\lambda_i$ (BC, PDE) & $(10.0,\ 1.0)$ \\
PRNG seed & 123 \\
\midrule
\multicolumn{2}{l}{\textit{Optimization}} \\
Optimizer & SOAP \\
$\beta_1$ & 0.9 \\
LR schedule & Step decay \\
LR decay rate & 0.9 \\
LR decay steps & 2{,}000 \\
\midrule
\multicolumn{2}{l}{\textit{Loss weighting}} \\
Weighting scheme & NTK \\
\bottomrule
\end{tabular}
\end{minipage}
\hfill
\begin{minipage}[t]{0.45\textwidth}
\centering
\textbf{MELM}\\[2mm]
\begin{tabular}{ll}
\toprule
\textbf{Parameter} & \textbf{Value} \\
\midrule
\multicolumn{2}{l}{\textit{Architecture}} \\
Network type & MELM (Legendre) \\
Layers & 1 \\
No. of Basis functions & $12{,}100$ \\
Basis & Shifted Legendre \\
\midrule
\multicolumn{2}{l}{\textit{Collocation and solution}} \\
Collocation points & $14{,}100$ \\
$\lambda_i$ (BC) & $ 1.0 $ \\
\midrule
\multicolumn{2}{l}{\textit{Linear solver}} \\
System size & $14{,}100\times12{,}100$ \\
Method & QR-based least-squares \\
\bottomrule
\end{tabular}
\end{minipage}

\end{table}

\begin{figure}[H]
    \centering
    \includegraphics[width=1\linewidth]{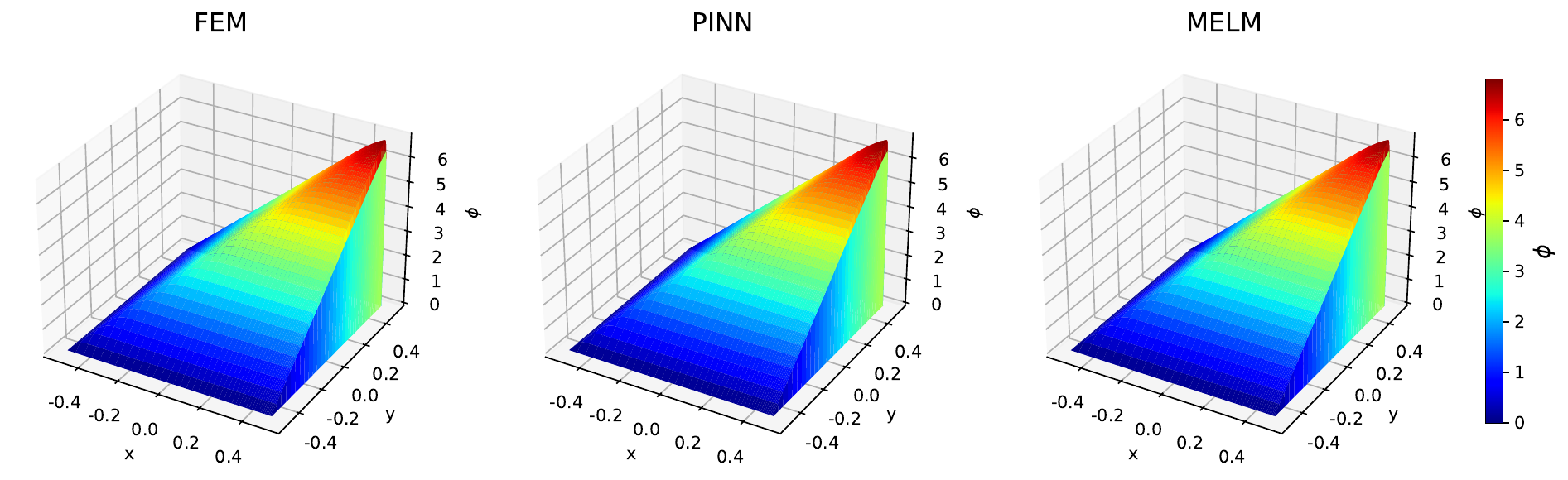}
    \caption{Comparison of the FEM, PINN, and MELM solutions for Problem~\ref{case_2}.}
    \label{fig:2_surface}
\end{figure}

\begin{figure}[H]
    \centering
    \includegraphics[width=1\linewidth]{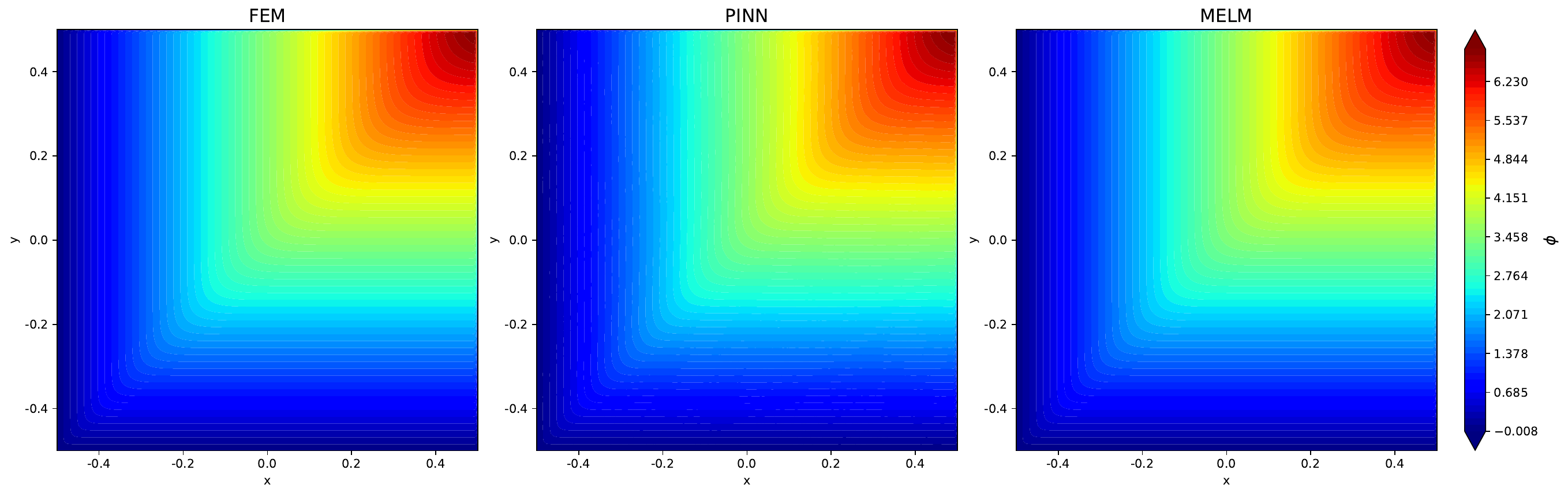}
    \caption{Comparison of the FEM, PINN, and MELM solution contours for Problem~\ref{case_2}.}
    \label{fig:2_contour}
\end{figure}

\begin{figure}[H]
    \centering
    \includegraphics[width=1\linewidth]{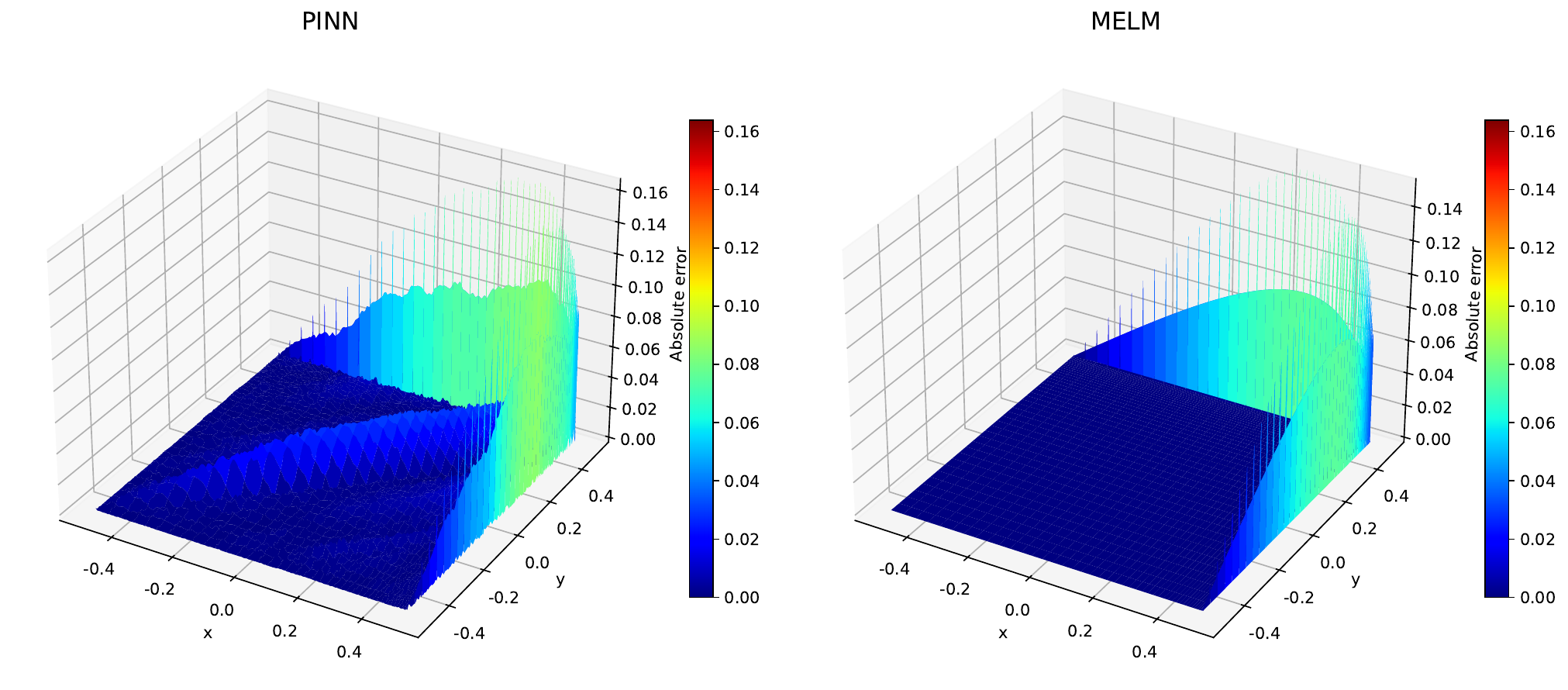}
    \caption{Comparison of the PINN and MELM absolute error surfaces for Problem~\ref{case_2}.}
    \label{fig:2_error}
\end{figure}

\begin{figure}[H]
    \centering
    \includegraphics[width=0.96\linewidth]{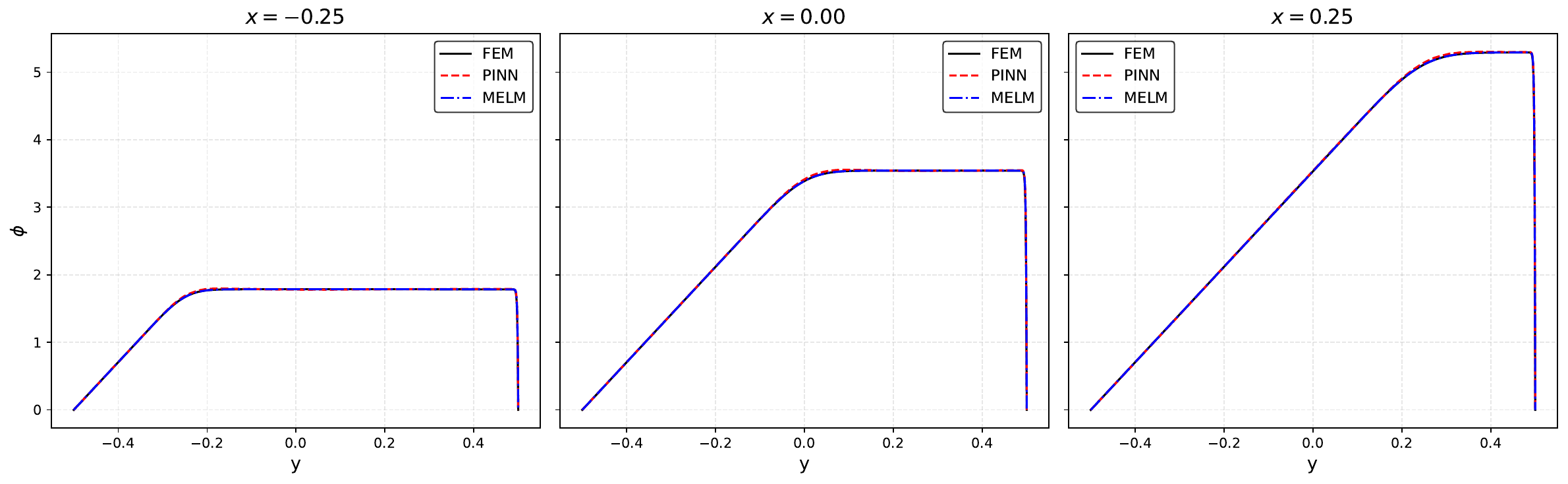}
    \caption{Cross-sectional comparison of the FEM, PINN, and MELM solutions for Problem~\ref{case_2} at selected values of $x$.}
    \label{fig:2_cross_section}
\end{figure}

\begin{figure}[H]
    \centering
    \includegraphics[width=0.9\linewidth]{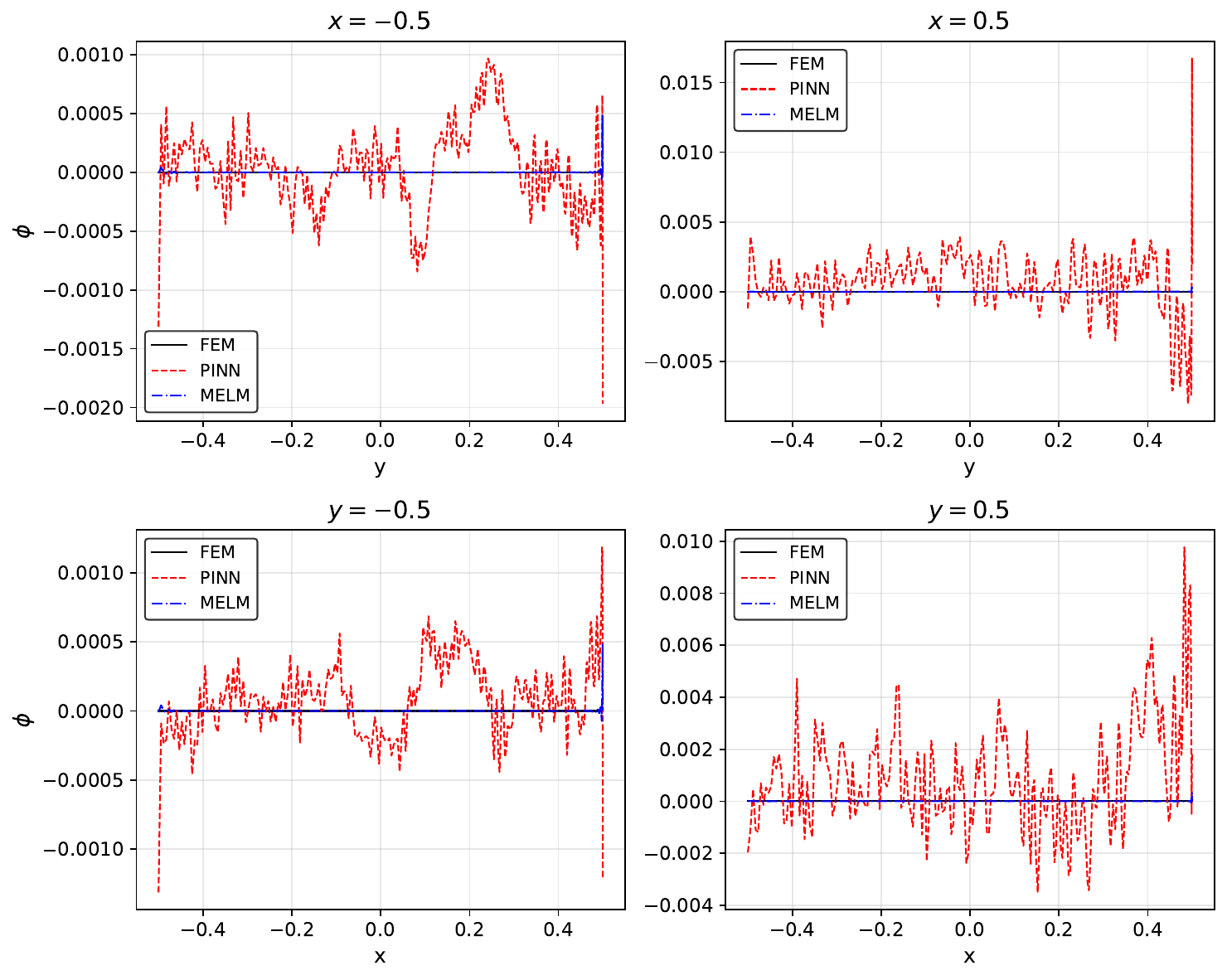}
    \caption{Comparison of the FEM, PINN, and MELM boundary profiles for Problem~\ref{case_2}.}
    \label{fig:2_boundary_profile}
\end{figure}

\begin{problem}[Two-Dimensional Steady Problem (Discontinuous BCs)]\label{case_3}
    Consider the steady two-dimensional convection--diffusion equation
\begin{equation*}
a_x\frac{\partial\phi}{\partial x}
+
a_y\frac{\partial\phi}{\partial y}
-
k\left(
\frac{\partial^{2}\phi}{\partial x^{2}}
+
\frac{\partial^{2}\phi}{\partial y^{2}}
\right)
=
0,
\qquad
(x,y)\in\Omega:=\left(-\frac{1}{2},\frac{1}{2}\right)
\times
\left(-\frac{1}{2},\frac{1}{2}\right),
\end{equation*}
subject to the discontinuous Dirichlet BCs
\begin{equation*}
\phi(x,y)=
\begin{cases}
1, & y=\dfrac{1}{2}, \\[2mm]
1, & x=-\dfrac{1}{2},\quad \dfrac{1}{4}\le y\le\dfrac{1}{2}, \\[2mm]
0, & \text{otherwise on } \partial\Omega.
\end{cases}
\end{equation*}
The problem parameters are
\begin{equation*}
a_x=\frac{1}{\sqrt{5}},
\qquad
a_y=-\frac{2}{\sqrt{5}},
\qquad
U=1,
\qquad
k=0.001,
\qquad
L=1,
\end{equation*}
yielding the Péclet number
\begin{equation*}
Pe=\frac{UL}{k}=1{,}000.
\end{equation*}
\end{problem}

The FEM reference solution is computed using a structured quadrilateral mesh with $90{,}000$ elements. Unlike Problem~\ref{case_2}, where refinement is directed toward the corner, the present mesh is refined most strongly around the boundary discontinuity and progressively coarsens away from this location. In consequence, mesh biasing is applied along both the $x$- and $y$-directions near the discontinuity and the element side lengths vary approximately from $0.00583$ to $0.000626$, with the smallest elements located immediately adjacent to the discontinuity, where the local elements retain an approximately square shape. Away from this region, the elements become progressively anisotropic. The resulting mesh provides enhanced spatial resolution of the strong gradients associated with the discontinuous BC and is used to obtain the FEM reference solution.

\begin{figure}[H]
    \centering
    \includegraphics[width=0.5\linewidth]{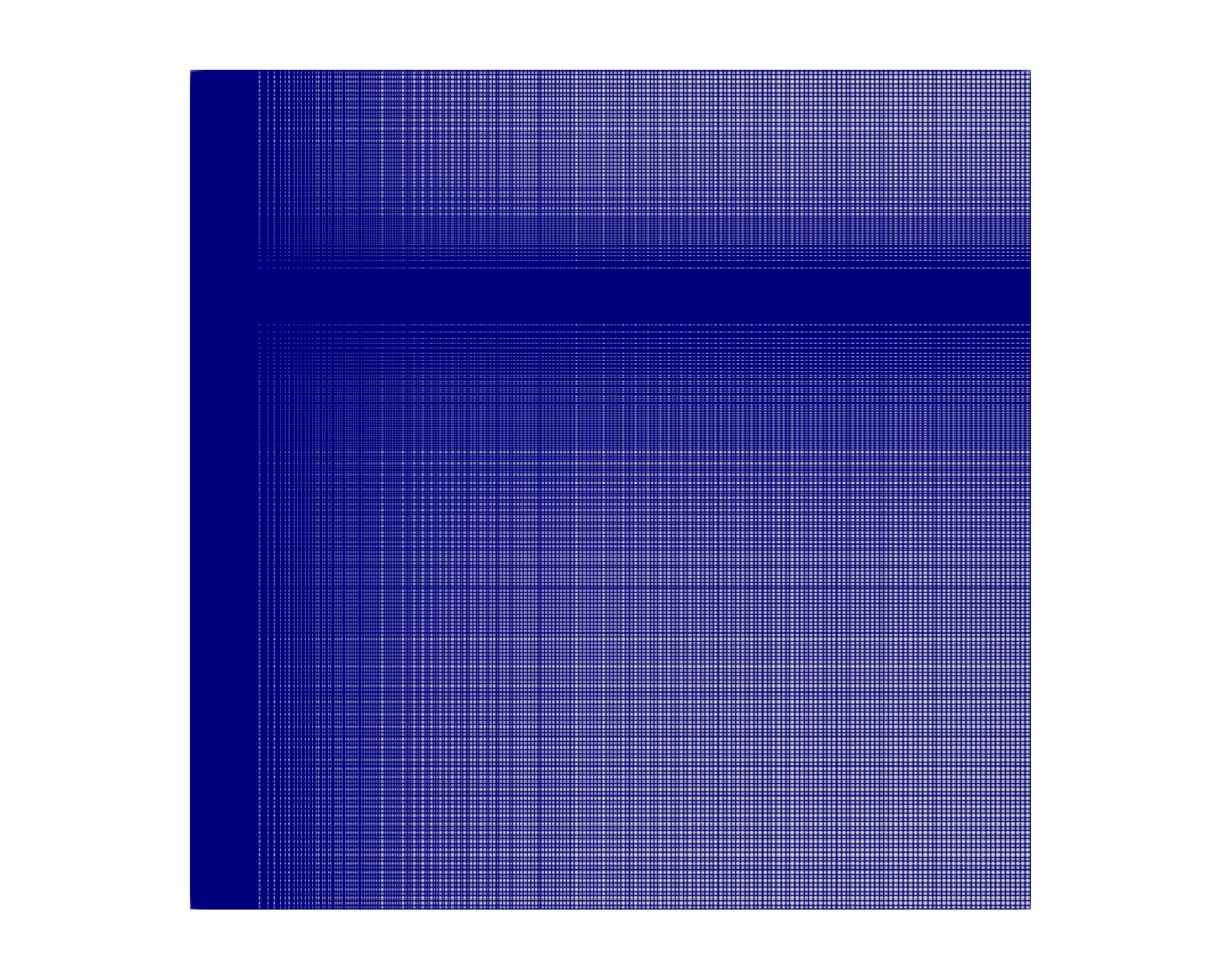}
    \caption{Structured FEM mesh used to compute the reference solution for Problem~\ref{case_3}. Local mesh refinement is introduced near the discontinuity on the left boundary.}
    \label{fig:P3_2dmesh}
\end{figure}

Table~\ref{tab:3_parameters} summarizes the implementation parameters used for the PINN and MELM methods in Problem~\ref{case_3}. For the MELM, $110$ basis functions are used along each of the $x$- and $y$-directions, resulting in $110\times110=12{,}100$ tensor-product basis functions. The numerical solutions obtained using the FEM, PINN, and MELM approaches are compared in Fig.~\ref{fig:3_surface}, while the corresponding contour representations are presented in Fig.~\ref{fig:3_contour}. Fig.~\ref{fig:3_error} illustrates the absolute error surfaces of the PINN and MELM solutions with respect to the FEM reference solution, and Figs.~\ref{fig:3_profiles} and~\ref{fig:3_boundary} compare the solution profiles at selected cross-sections and along the domain boundary, respectively.

The results for $Pe=1{,}000$ illustrate the behavior of the PINN and MELM approximations in the presence of discontinuous boundary data. As shown in Fig.~\ref{fig:3_boundary}, the MELM boundary solution exhibits localized oscillations in the vicinity of the discontinuities, consistent with the Gibbs phenomenon typically associated with high-order spectral approximations. The PINN, in contrast, does not exhibit comparable oscillations but shows larger pointwise deviations near the discontinuity locations. Away from the discontinuities, both approaches provide close agreement with the FEM reference solution. The corresponding error surfaces and cross-sectional profiles further illustrate that the largest deviations are concentrated in regions affected by the sharp gradients generated by the discontinuous BC. Quantitative comparisons in terms of relative $L_2$ error and computational time are reported in Table~\ref{tab:comparison}.

\begin{table}[H]
\centering
\caption{Implementation parameters used for PINN and MELM for Problem~\ref{case_3}.}
\label{tab:3_parameters}

\begin{minipage}[t]{0.48\textwidth}
\centering
\textbf{PINN}\\[2mm]
\begin{tabular}{ll}
\toprule
\textbf{Parameter} & \textbf{Value} \\
\midrule
\multicolumn{2}{l}{\textit{Architecture}} \\
Network type & Modified MLP \\
Layers & 5 \\
Neurons & 256 \\
Activation & GELU \\
Fourier feature scale & 30.0 \\
RWF & $\mu=1.0,\ \sigma=0.1$ \\
\midrule
\multicolumn{2}{l}{\textit{Collocation and training}} \\
Batch size (residual) & 512 \\
Batch size (boundary) & 256 \\
Learning rate & 0.001 \\
Training iterations & 200{,}000 \\
$\lambda_i$ (BC, PDE) & $(20.0,\ 1.0)$ \\
PRNG seed & 123 \\
\midrule
\multicolumn{2}{l}{\textit{Optimization}} \\
Optimizer & SOAP \\
$\beta_1$ & 0.9 \\
LR schedule & Step decay \\
LR decay rate & 0.9 \\
LR decay steps & 2{,}000 \\
\midrule
\multicolumn{2}{l}{\textit{Loss weighting}} \\
Weighting scheme & NTK \\
\bottomrule
\end{tabular}
\end{minipage}
\hfill
\begin{minipage}[t]{0.45\textwidth}
\centering
\textbf{MELM}\\[2mm]
\begin{tabular}{ll}
\toprule
\textbf{Parameter} & \textbf{Value} \\
\midrule
\multicolumn{2}{l}{\textit{Architecture}} \\
Network type & MELM (Legendre) \\
Layers & 1 \\
No. of Basis functions & $12{,}100$ \\
Basis & Shifted Legendre \\
\midrule
\multicolumn{2}{l}{\textit{Collocation and solution}} \\
Collocation points & $14{,}100$ \\
$\lambda_i$ (BC) & $ 1.0 $ \\
\midrule
\multicolumn{2}{l}{\textit{Linear solver}} \\
System size & $14{,}100\times12{,}100$ \\
Method & QR-based least-squares \\
\bottomrule
\end{tabular}
\end{minipage}

\end{table}

\begin{figure}[H]
    \centering
    \includegraphics[width=1\linewidth]{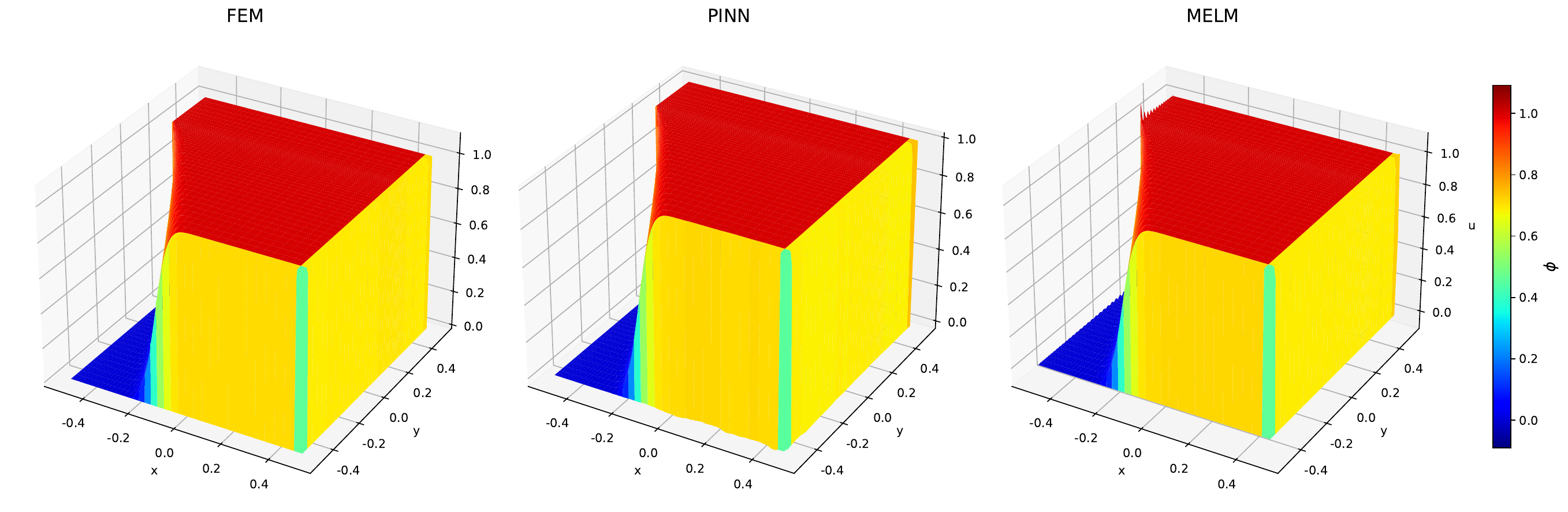}
    \caption{Comparison of the FEM, PINN, and MELM solutions for Problem~\ref{case_3}.}
    \label{fig:3_surface}
\end{figure}

\begin{figure}[H]
    \centering
    \includegraphics[width=1\linewidth]{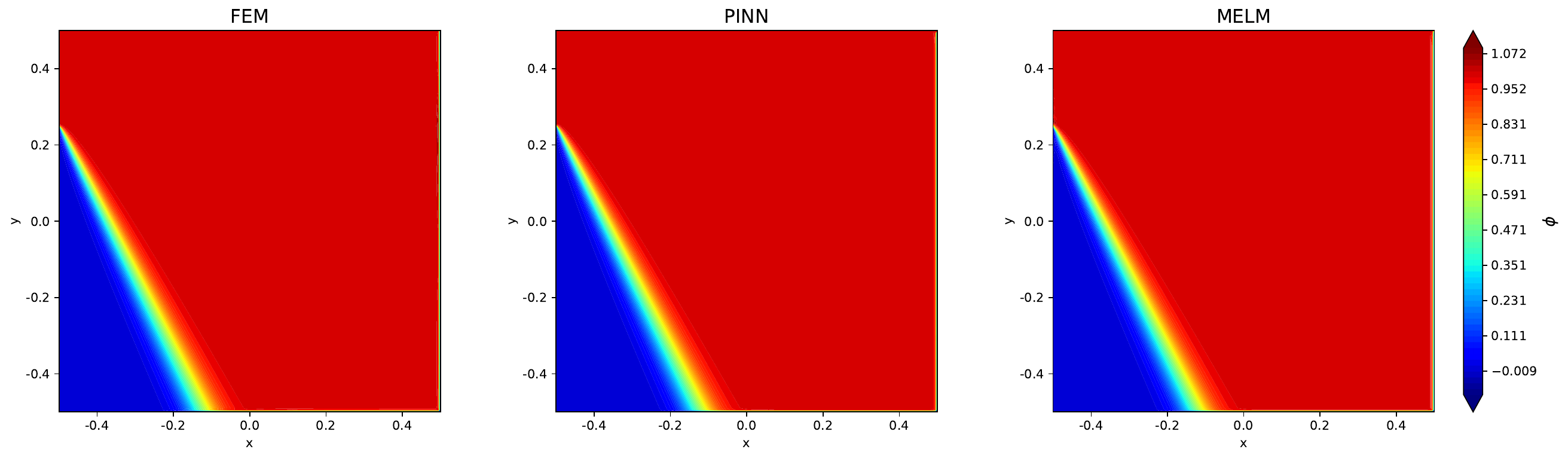}
    \caption{Comparison of the FEM, PINN, and MELM solution contours for Problem~\ref{case_3}.}
   \label{fig:3_contour}
\end{figure}

\begin{figure}[H]
    \centering
    \includegraphics[width=1\linewidth]{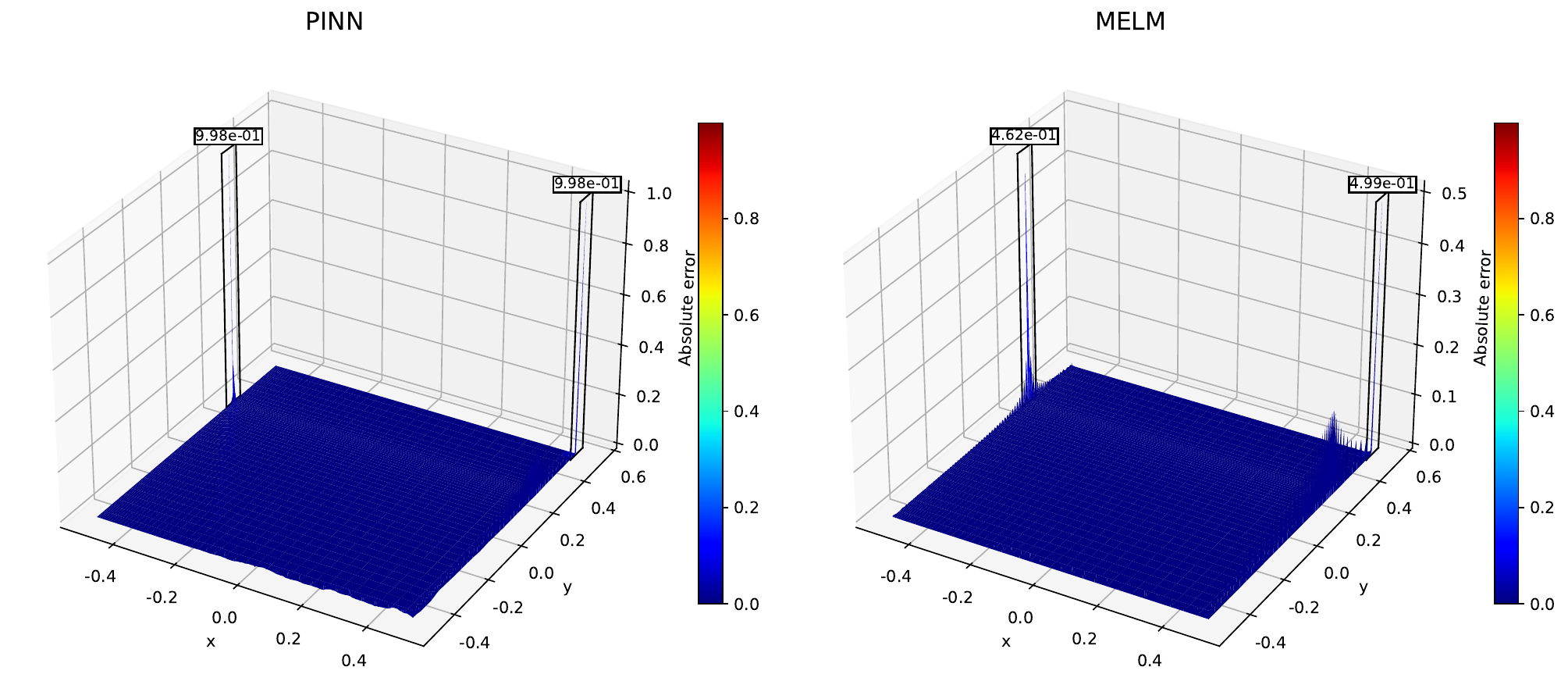}
   \caption{Comparison of the PINN and MELM absolute error surfaces for Problem~\ref{case_3}.}
    \label{fig:3_error}
\end{figure}

\begin{figure}[H]
    \centering
    \includegraphics[width=0.96\linewidth]{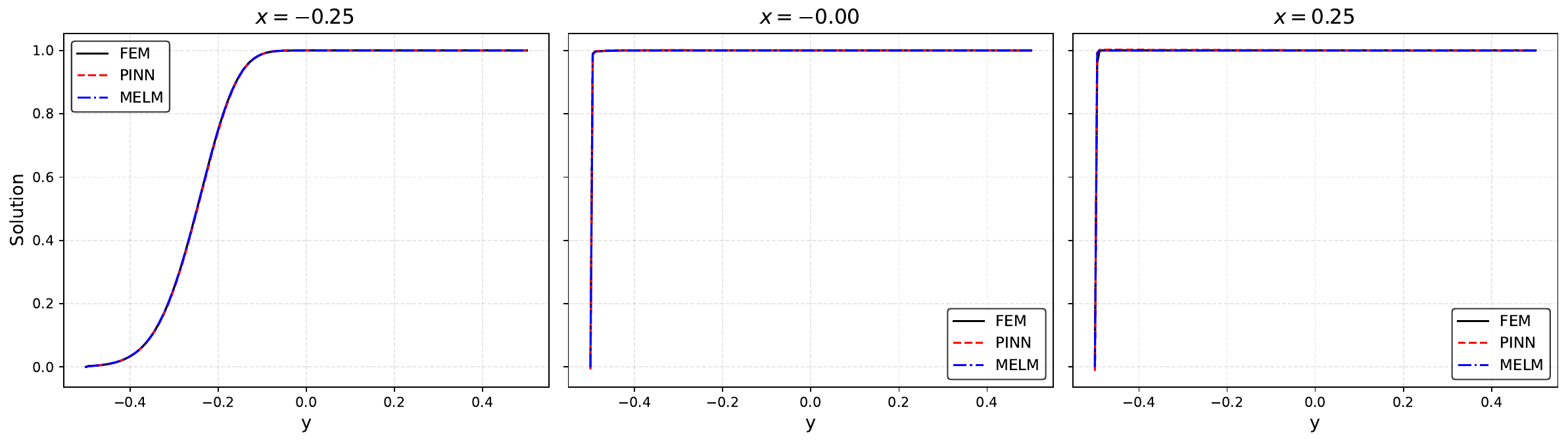}
    \caption{Cross-sectional comparison of the FEM, PINN, and MELM solutions for Problem~\ref{case_3} at selected values of $x$.}
   \label{fig:3_profiles}
\end{figure}

\begin{figure}[H]
    \centering
    \includegraphics[width=0.9\linewidth]{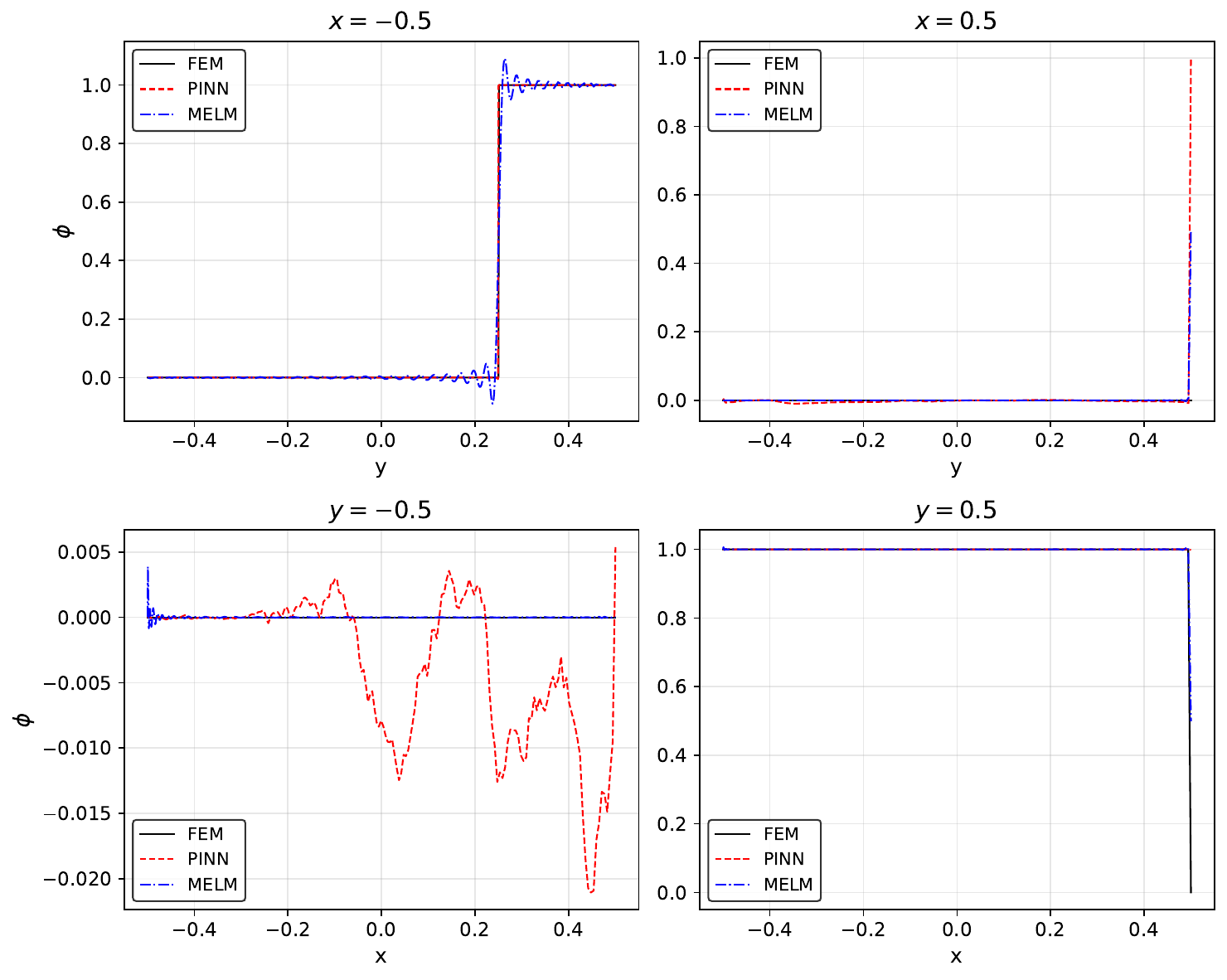}
    \caption{Comparison of the FEM, PINN, and MELM boundary profiles for Problem~\ref{case_3}.}
   \label{fig:3_boundary}
\end{figure}

The differences observed near the boundary discontinuities also illustrate the distinct approximation mechanisms underlying the three numerical approaches. The stabilized FEM solution benefits from a locally refined mesh and residual-based stabilization, which contribute to resolving the sharp gradients generated by the discontinuous boundary data. In contrast, MELM employs a globally supported polynomial representation, leading to localized oscillations consistent with the Gibbs phenomenon. Interestingly, the PINN does not exhibit comparable oscillatory behavior, despite producing larger pointwise deviations near the discontinuities. This suggests that the flexibility of the trainable nonlinear representation may allow smoother approximations of abrupt transitions, although this behavior does not necessarily imply greater local accuracy. The results therefore highlight a distinction between suppressing visible oscillations and accurately resolving sharp solution features. Furthermore, the present comparison does not isolate the respective contributions of FEM stabilization, mesh refinement, and PINN architecture, and should not be interpreted as demonstrating an intrinsic shock-capturing capability of PINNs.

\begin{problem}[One-Dimensional Unsteady Convection-Diffusion Problem]\label{case_4}
    Consider the one-dimensional unsteady convection--diffusion equation
\begin{equation*}
\frac{\partial \phi}{\partial t}
+
a\frac{\partial \phi}{\partial x}
-
k\frac{\partial^2 \phi}{\partial x^2}
=
0,
\qquad
(x,t)\in\Omega\times(0,T],
\end{equation*}
defined on the domain
\(
\Omega=(0,L),
\)
where \(a\) is the convection velocity and \(k\) is the diffusion coefficient. The exact solution is given by
\begin{equation*}
\phi(x,t)
=
\frac{1}{\sqrt{4t+1}}
\exp\!\left(
-\frac{(x-1-at)^2}
{k(4t+1)}
\right).
\end{equation*}
The initial and Dirichlet BCs are obtained directly from the exact solution.
The problem parameters are
\begin{equation*}
U=1,
\qquad
a=1,
\qquad
k=0.1,
\qquad
L=9,
\qquad
T=5,
\end{equation*}
yielding the Péclet number
\begin{equation*}
Pe=\frac{UL}{k}=90.
\end{equation*}
\end{problem}

The FEM reference solution is computed using a one-dimensional uniform mesh with $5{,}000$ elements and $5{,}001$ nodes, corresponding to a uniform element size of $\Delta x=0.0018$. The same spatial discretization is used as the reference grid for evaluating the transient numerical solutions.

Table~\ref{tab:4_parameters} summarizes the implementation parameters used for the PINN and MELM methods for both the CT and DT formulations.

\begin{table}[H]
\centering
\caption{Parameters used for the CT-PINN, DT-PINN, CT-MELM, and DT-MELM models for Problem~\ref{case_4}.}
\label{tab:4_parameters}

\begin{minipage}[t]{0.48\textwidth}
\centering
\textbf{PINN}\\[2mm]
\begin{tabularx}{\linewidth}{@{}>{\raggedright\arraybackslash}X
                              >{\centering\arraybackslash}X
                              >{\centering\arraybackslash}X@{}}
\toprule
\textbf{Parameter} & \textbf{CT} & \textbf{DT} \\
\midrule
$T$ & 5 & 5 \\
$\Delta t$ & --- & 0.01 \\
\midrule
\multicolumn{3}{l}{\textit{Architecture}} \\
Network type & Modified MLP & Modified MLP \\
Layers & 4 & 3 \\
Neurons & 128 & 64 \\
Activation & tanh & GELU \\
Fourier feature scale & 1.0 & 15.0 \\
RWF & \makecell{$\mu=1.0,$\\$\sigma=0.1$}
    & \makecell{$\mu=0.5,$\\$\sigma=0.1$} \\
\midrule
\multicolumn{3}{l}{\textit{Collocation and training}} \\
Batch size (residual) & 128 & 256 \\
Batch size (initial) & 32 & --- \\
Learning rate & 0.001 & 0.001 \\
Training iterations & 180{,}000 & 3{,}000 \\
$\lambda_i$ (BC, PDE, IC) & $(2.0,\ 1.0,\ 2.0)$ & $(1.0,\ 2.0,\ -)$ \\
PRNG seed & 42 & 123 \\
\midrule
\multicolumn{3}{l}{\textit{Optimization}} \\
Optimizer & SOAP & SOAP \\
$\beta_1$ & 0.9 & 0.9 \\
LR schedule & Step decay & Step decay \\
LR decay rate & 0.9 & 0.9 \\
LR decay steps & 2{,}000 & 2{,}000 \\
\midrule
\multicolumn{3}{l}{\textit{Loss weighting}} \\
Weighting scheme & NTK & NTK \\
Causal tolerance & 1.0 & --- \\
Number of chunks & 32 & --- \\
\bottomrule
\end{tabularx}
\end{minipage}
\hfill
\begin{minipage}[t]{0.45\textwidth}
\centering
\textbf{MELM}\\[2mm]
\begin{tabularx}{\linewidth}{@{}>{\raggedright\arraybackslash}X
                              >{\centering\arraybackslash}X
                              >{\centering\arraybackslash}X@{}}
\toprule
\textbf{Parameter} & \textbf{CT} & \textbf{DT} \\
\midrule
$T$ & 5 & 5 \\
$\Delta t$ & --- & 0.01 \\
\midrule
\multicolumn{3}{l}{\textit{Architecture}} \\
Network type & MELM (Legendre) & MELM (Legendre) \\
Layers & 1 & 1 \\
No. of Basis functions & 8{,}075 & 60 \\
Basis & Shifted Legendre & Shifted Legendre \\
\midrule
\multicolumn{3}{l}{\textit{Collocation and solution}} \\
Collocation points & 9{,}384 & 60 \\
$\lambda_i$ (BC,\ IC) & $(1.0,\ 1.0)$ & $(1.0,\ -)$ \\
\midrule
\multicolumn{3}{l}{\textit{Linear solver}} \\
System size & $9{,}384\times8{,}075$ & $60\times60$ \\
Method & QR-based least-squares & LU-based direct solver \\
\bottomrule
\end{tabularx}
\end{minipage}

\end{table}

\subsubsection*{Case 1: Continuous-Time Results}
The number of basis functions used for the MELM is 95 along the
$x$-direction and 85 along the $t$-direction, yielding a total of
$95\times85=8{,}075$ tensor-product basis functions.
Fig.~\ref{fig:problem4_comparison} compares the exact solution with
the PINN and MELM solutions using contour plots.
Fig.~\ref{fig:profile} compares the exact and numerical solution
profiles at selected time instants. Both the PINN and MELM closely follow
the exact solution as the profile propagates in the positive
$x$-direction and gradually decreases in amplitude with time. The close agreement between the numerical and exact profiles shows that both methods reproduce the main spatial and temporal features of the solution. Fig.~\ref{fig:prob4_error_comparison} compares the
absolute error distributions of the PINN and MELM solutions. For both
methods, the maximum pointwise errors are concentrated in the region
of the sharp solution gradient, where the solution changes rapidly over
a small spatial interval.

\begin{figure}[!htbp]
\centering

\begin{subfigure}{0.32\textwidth}
    \centering
    \includegraphics[width=\linewidth]{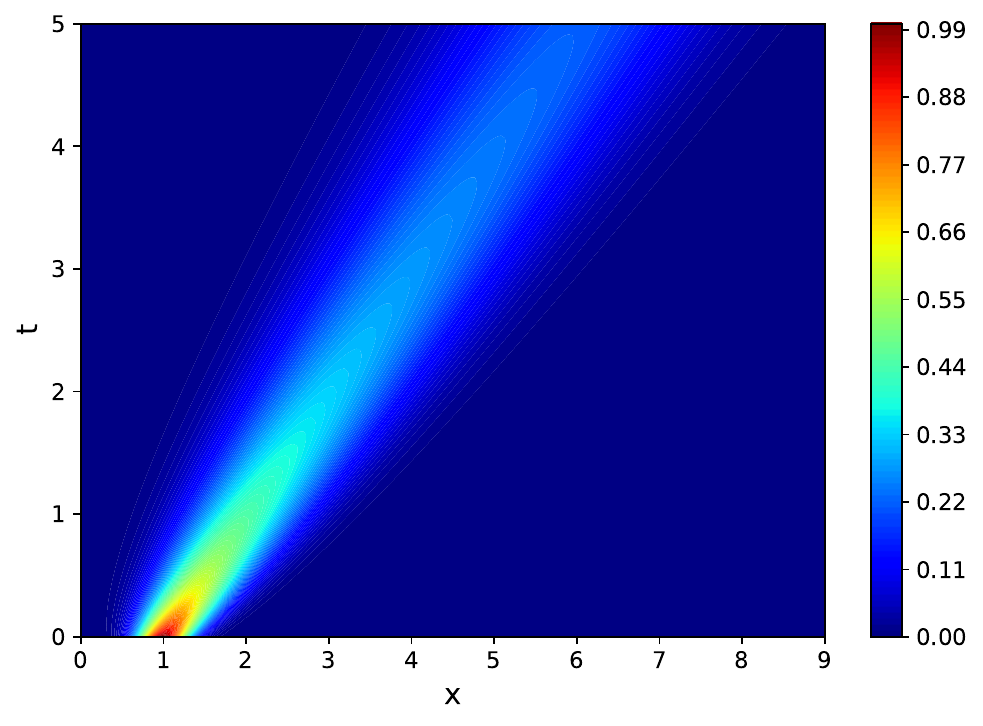}
    \caption{Exact}
    \label{fig:problem4_exact}
\end{subfigure}
\hfill
\begin{subfigure}{0.32\textwidth}
    \centering
    \includegraphics[width=\linewidth]{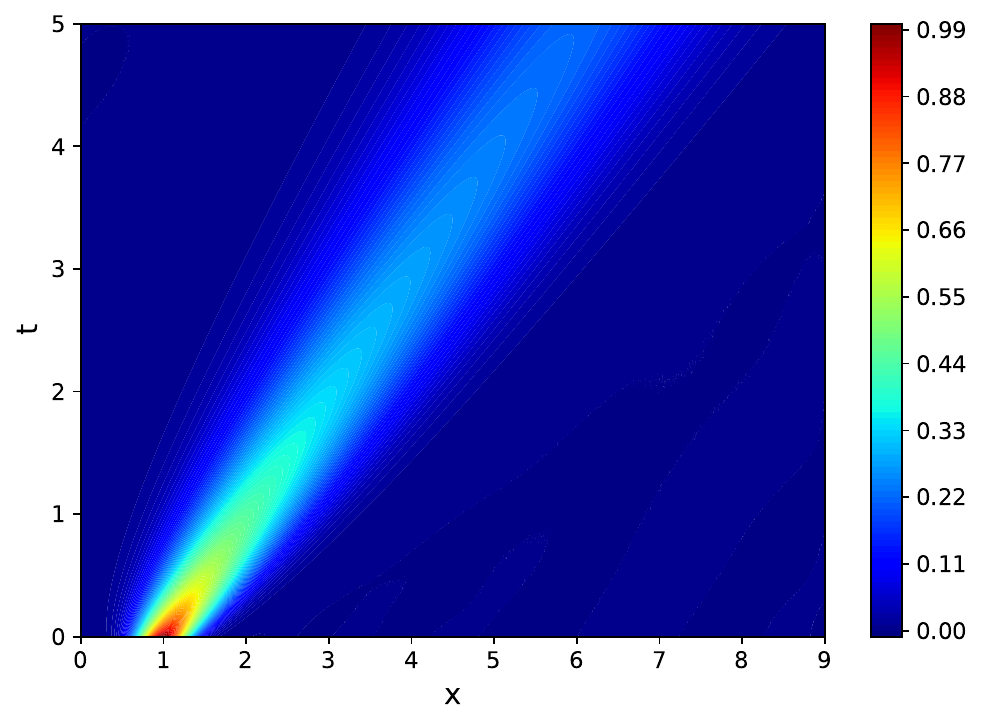}
    \caption{PINN}
    \label{fig:problem4_pinn}
\end{subfigure}
\hfill
\begin{subfigure}{0.32\textwidth}
    \centering
    \includegraphics[width=\linewidth]{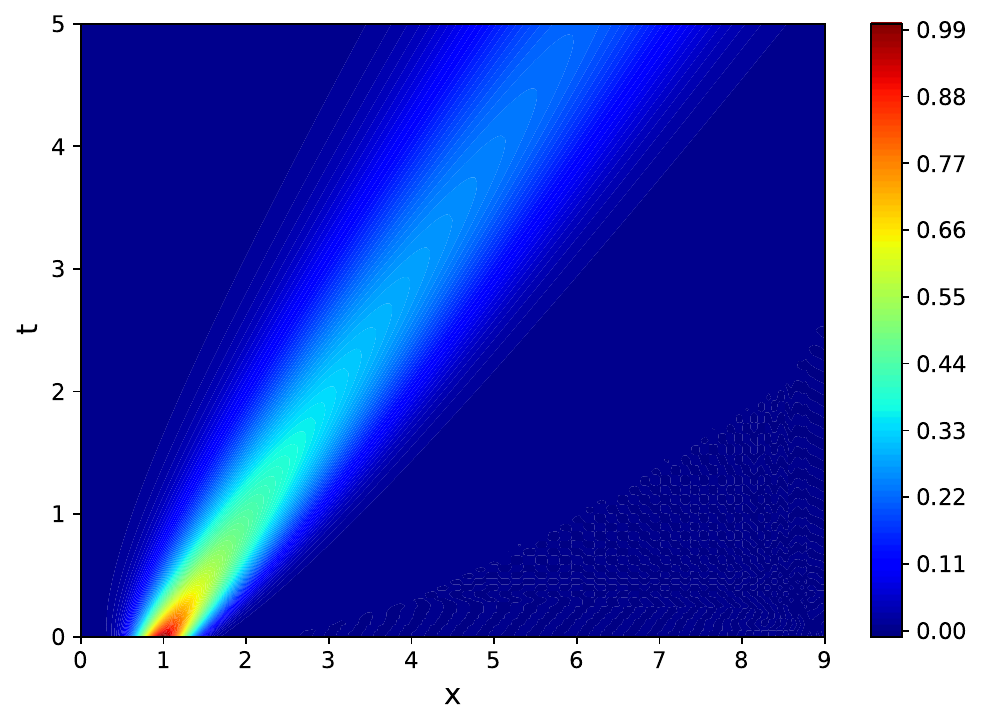}
    \caption{MELM}
    \label{fig:problem4_melm}
\end{subfigure}

\caption{Comparison of the exact, PINN, and MELM solutions for Problem~\ref{case_4}.}
\label{fig:problem4_comparison}
\end{figure}

\begin{figure}[!htbp]
    \centering
    \begin{subfigure}[b]{0.49\textwidth}
        \centering
        \includegraphics[width=\linewidth]{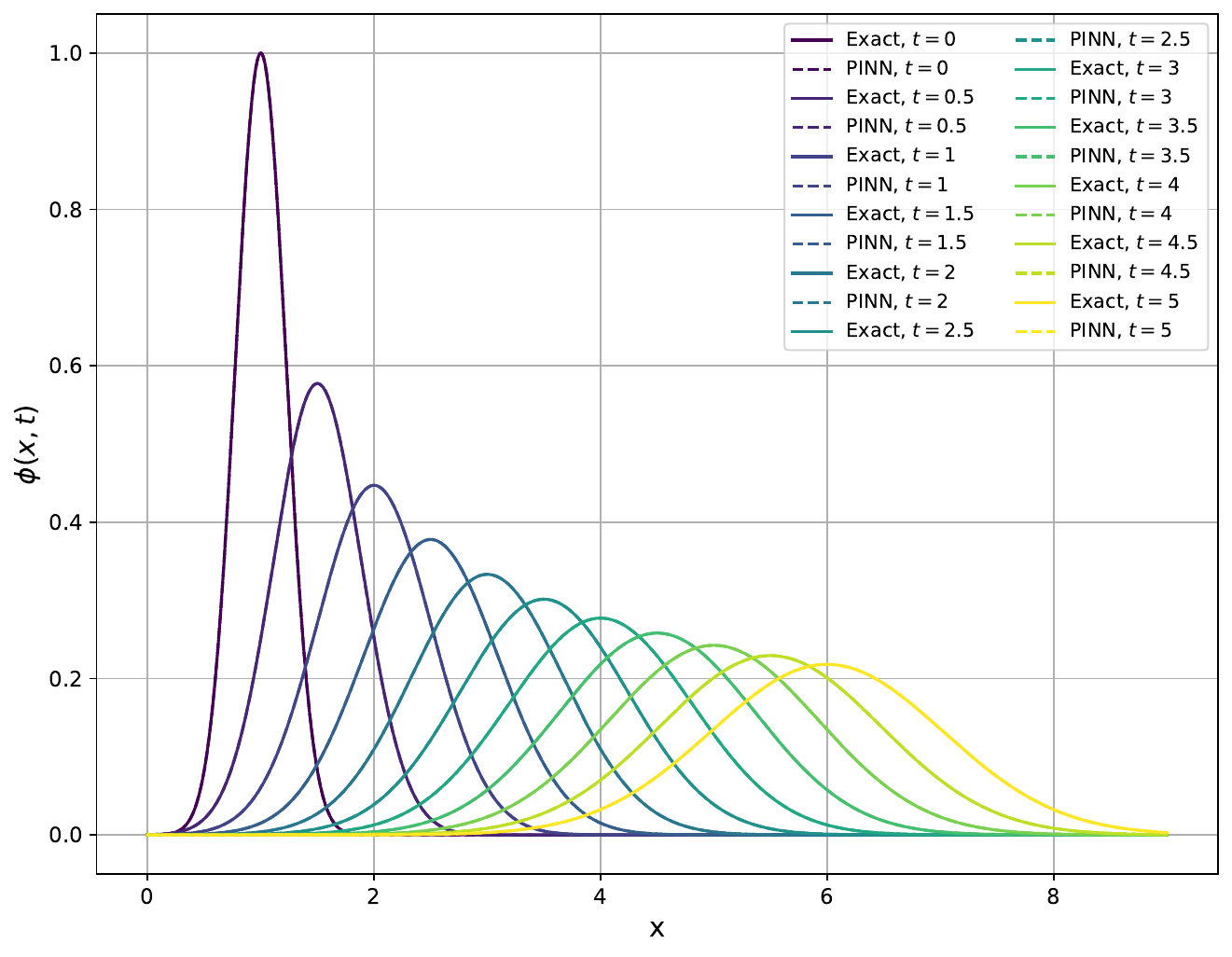}
        \caption{PINN}
        \label{fig:profile_pinn}
    \end{subfigure}
    \hfill
    \begin{subfigure}[b]{0.49\textwidth}
        \centering
        \includegraphics[width=\linewidth]{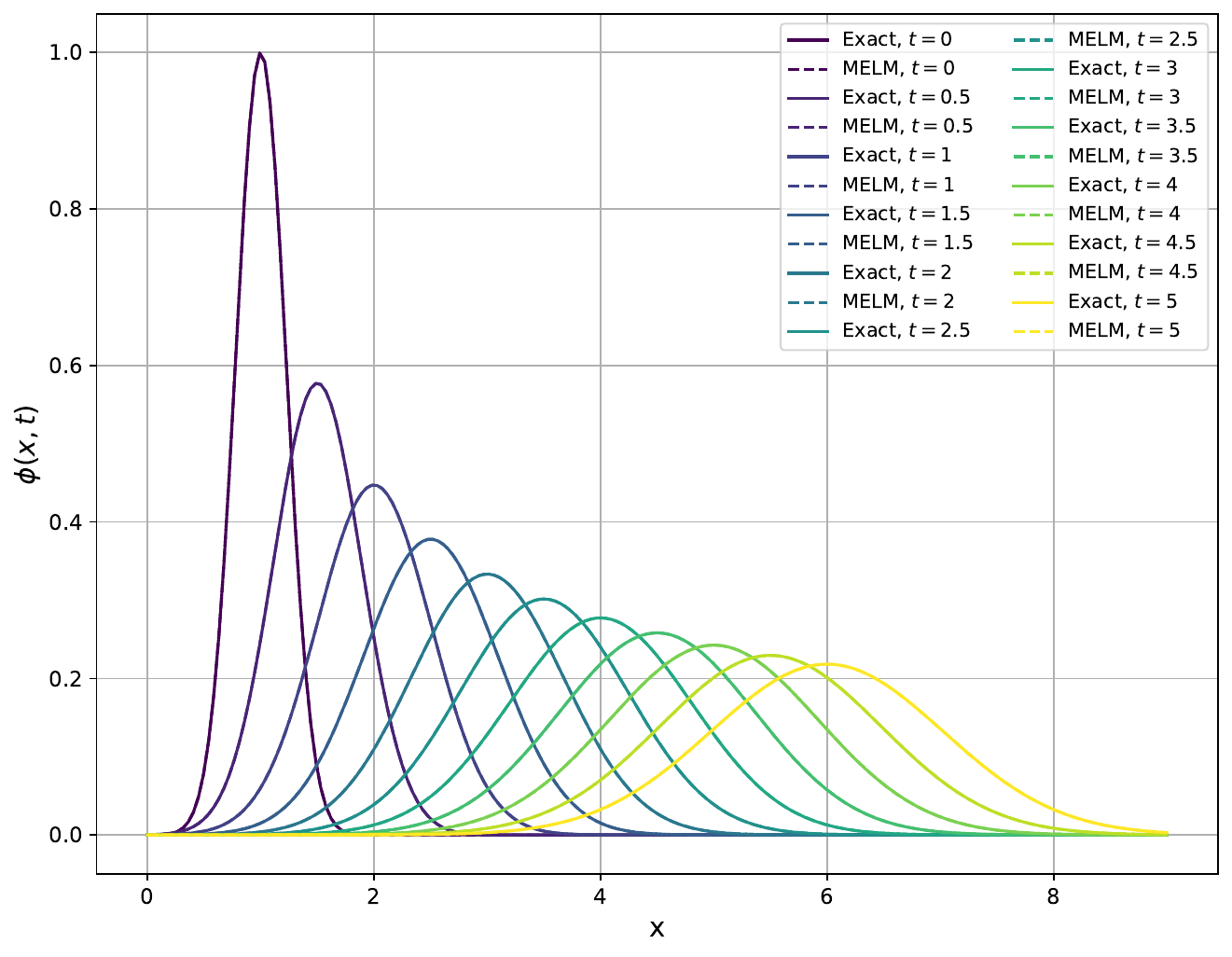}
        \caption{MELM}
        \label{fig:profile_melm}
    \end{subfigure}

    \caption{Exact and numerical solution profiles at selected time instants for Problem \ref{case_4}.}
    \label{fig:profile}
\end{figure}

\begin{figure}[H]
    \centering

    \begin{subfigure}{0.48\textwidth}
        \centering
        \includegraphics[width=\linewidth]{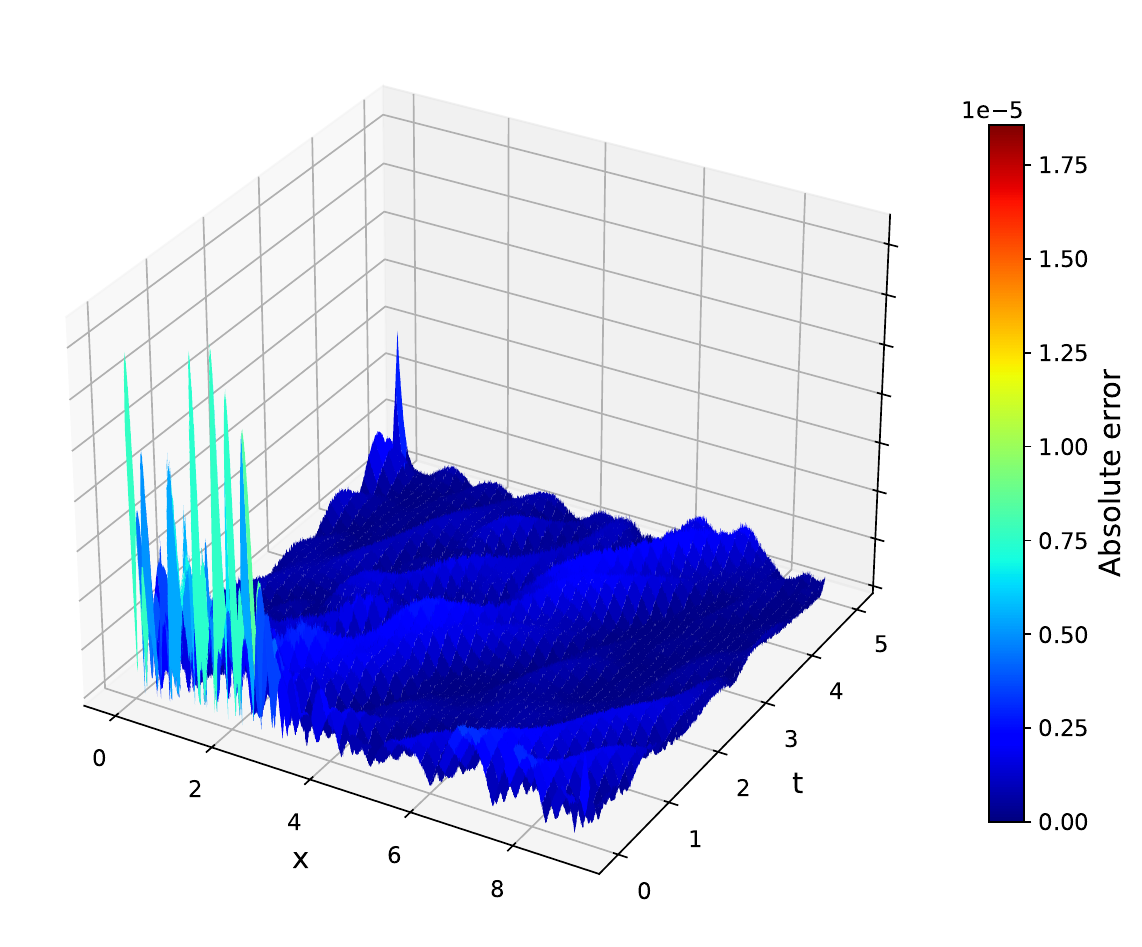}
        \caption{PINN}
        \label{fig:prob4_pinn_error}
    \end{subfigure}
    \hfill
    \begin{subfigure}{0.48\textwidth}
        \centering
        \includegraphics[width=\linewidth]{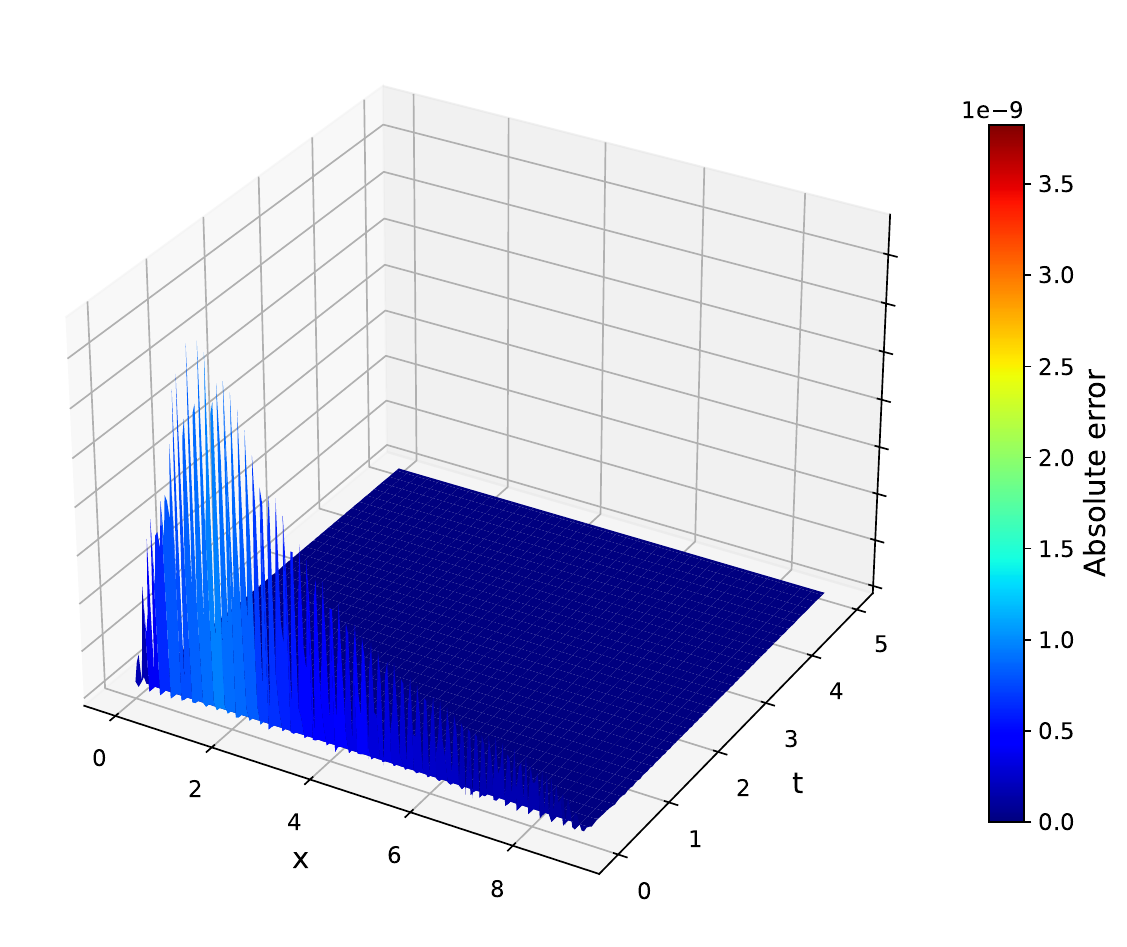}
        \caption{MELM}
        \label{fig:prob4_melm_error}
    \end{subfigure}

    \caption{Comparison of the absolute error distributions for the CT formulation for Problem \ref{case_4}.}
    \label{fig:prob4_error_comparison}
\end{figure}

\subsubsection*{Case 2: Discrete-Time Results}

For both methods, the time step is taken as $\Delta t=0.01$. For the
MELM, 60 basis functions are employed along the $x$-direction, and the
number of collocation points is chosen to match the number of basis
functions, resulting in a square system matrix. The resulting linear
system is then solved using LU decomposition.
Fig.~\ref{fig:profile_dt_1} compares the exact and numerical solution
profiles at selected time instants. The profiles illustrate the temporal evolution obtained with both methods, while the corresponding error distributions reveal differences in their approximation accuracy.

\begin{figure}[!htbp]
    \centering

    \begin{subfigure}[b]{0.49\textwidth}
        \centering
        \includegraphics[width=\linewidth]{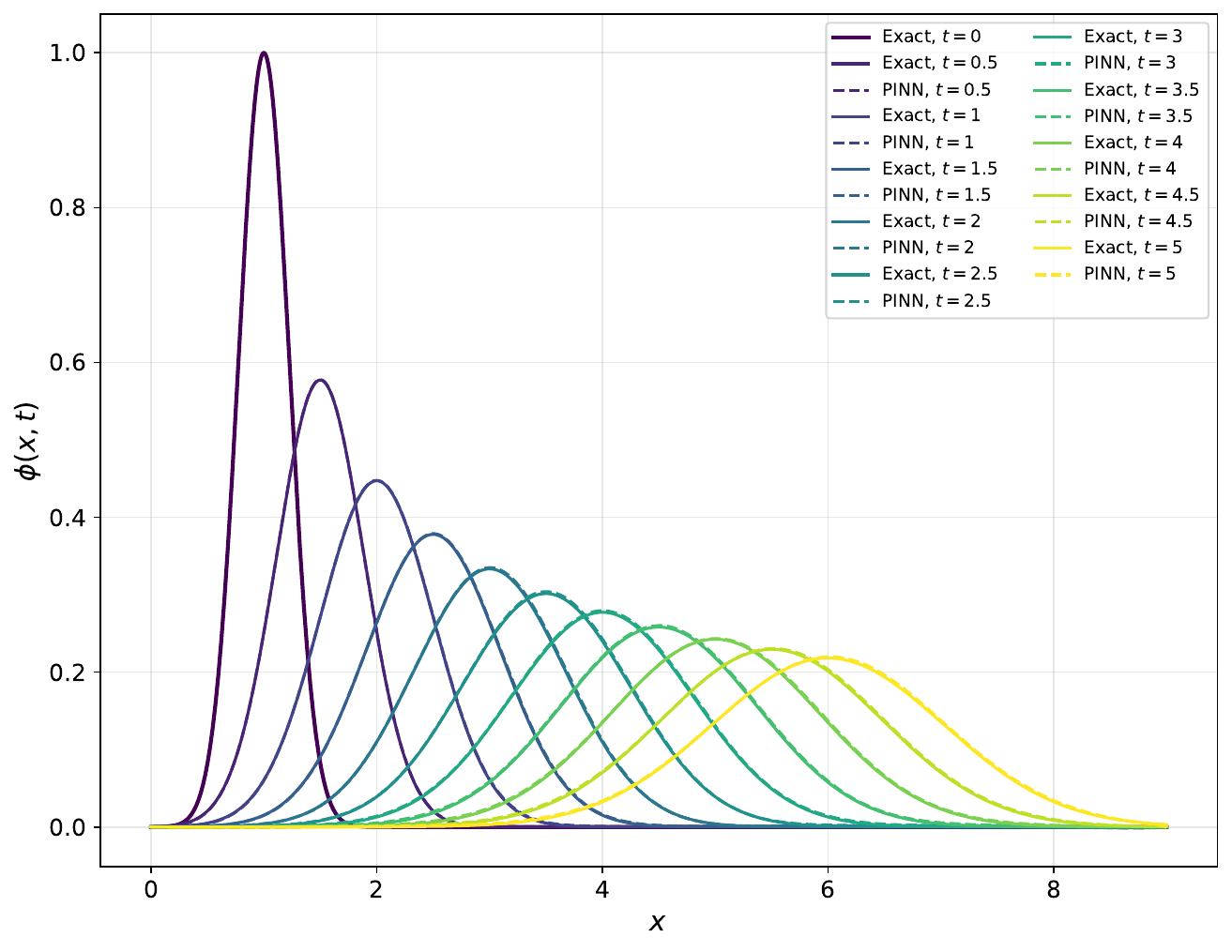}
        \caption{PINN}
        \label{fig:profile_dt_pinn}
    \end{subfigure}
    \hfill
    \begin{subfigure}[b]{0.49\textwidth}
        \centering
        \includegraphics[width=\linewidth]{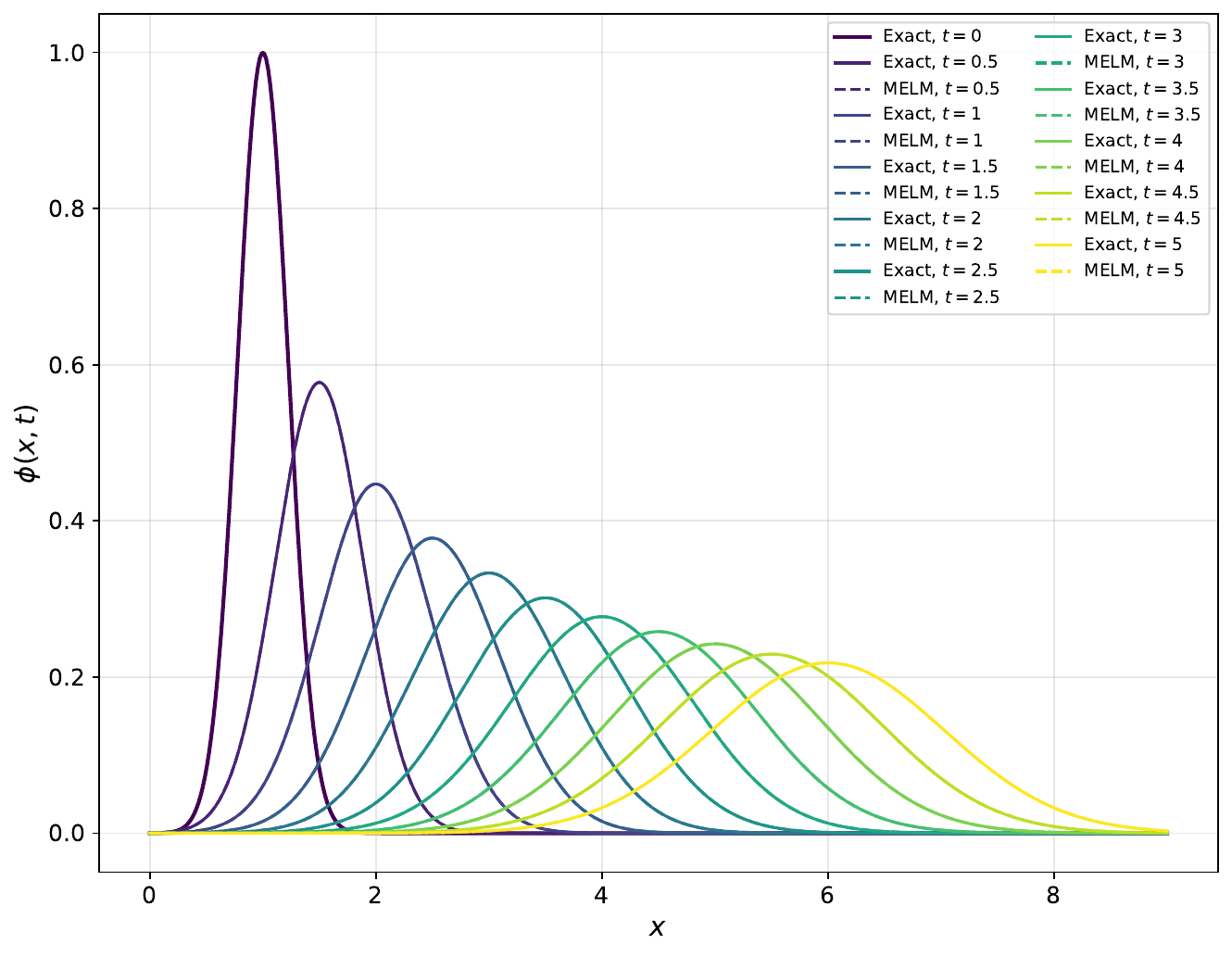}
        \caption{MELM}
        \label{fig:profile_dt_melm}
    \end{subfigure}

     \caption{Exact and numerical solution profiles at selected time instants for Problem \ref{case_4}.}
    \label{fig:profile_dt_1}
\end{figure}

Fig.~\ref{fig:prob4_error_comparison_1} compares the absolute error
distributions of the PINN and MELM solutions. The error plots provide a
visual assessment of the spatial accuracy achieved by the two methods
over the computational domain.
\begin{figure}[!htbp]
    \centering

    \begin{subfigure}{0.48\textwidth}
        \centering
        \includegraphics[width=\linewidth]{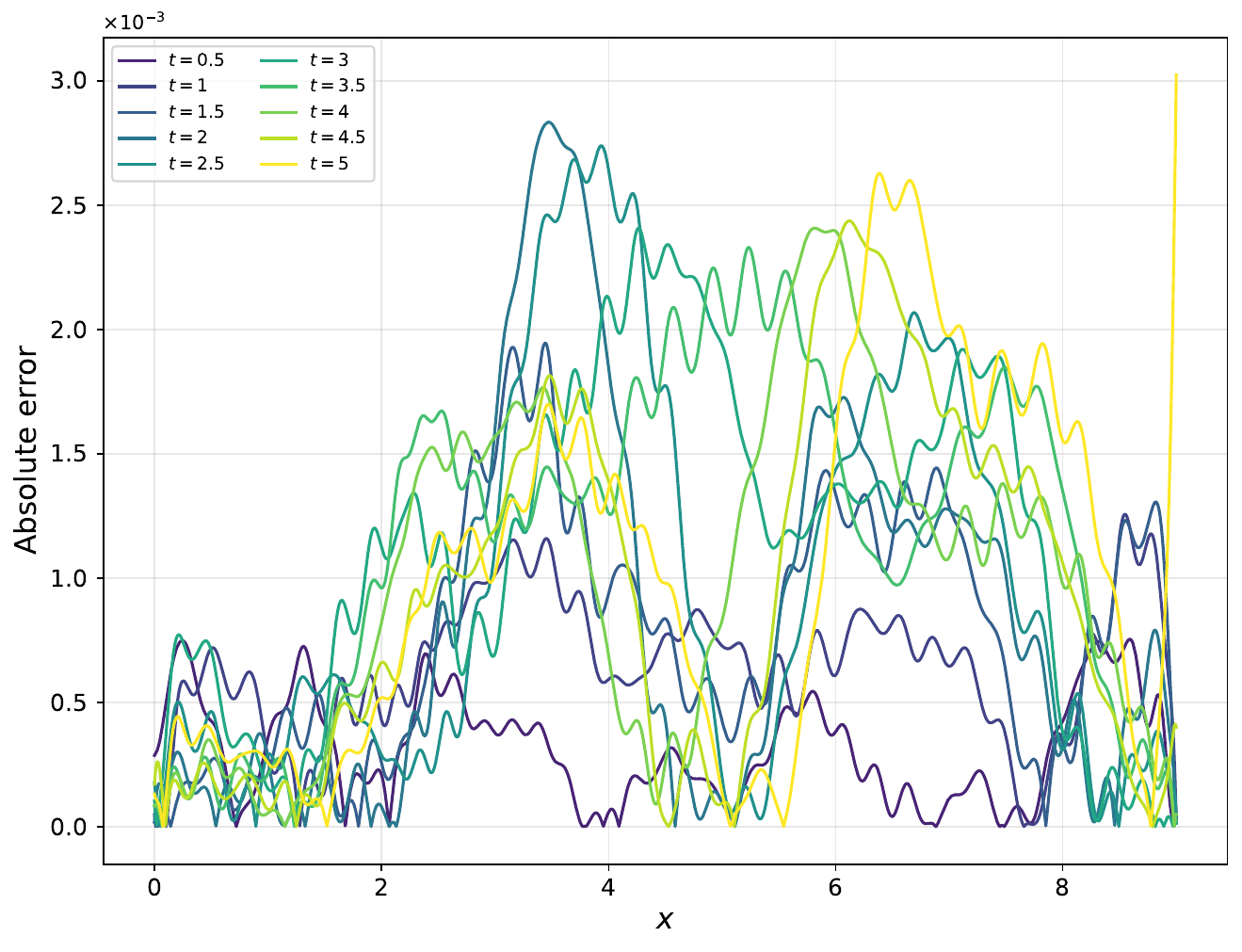}
        \caption{PINN}
        \label{fig:prob4_dt_pinn_error}
    \end{subfigure}
    \hfill
    \begin{subfigure}{0.48\textwidth}
        \centering
        \includegraphics[width=\linewidth]{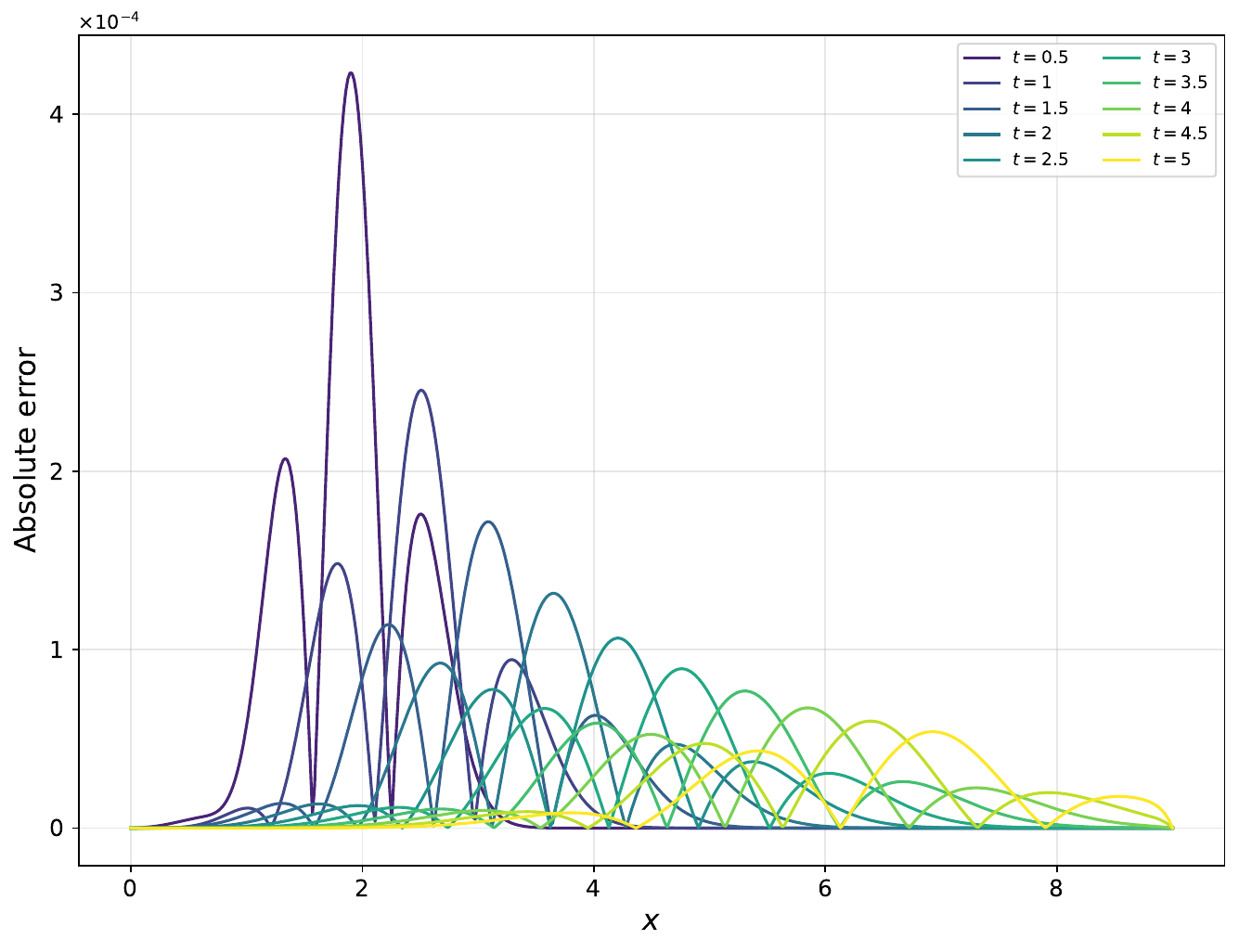}
        \caption{MELM}
        \label{fig:prob4_dt_melm_error}
    \end{subfigure}

    \caption{Comparison of the absolute error distributions for the DT formulation for Problem \ref{case_4}.}
    \label{fig:prob4_error_comparison_1}
\end{figure}

Fig.~\ref{fig:P4_error_vs_time} compares the temporal evolution of
the maximum absolute error and relative $L_2$ error for the PINN and
MELM. For the PINN in Fig.~\ref{fig:P4_linf_error}, the relative
$L_2$ error increases progressively with time and reaches the order of
$10^{-2}$ at later times. The maximum absolute error initially
decreases, followed by an increase and noticeable fluctuations, and remains of the order of $10^{-3}$ during the later stages of the simulation.

In contrast, the MELM in Fig.~\ref{fig:P4_relL2_error} shows a
consistent decrease in both error measures with time. The maximum
absolute error decreases to approximately the order of $10^{-5}$,
while the relative $L_2$ error decreases to approximately the order of
$10^{-4}$ at $t=5$. Thus, the MELM achieves substantially lower errors and exhibits more stable temporal error behavior than the PINN for the considered problem.

\begin{figure}[htbp]
    \centering

    \begin{subfigure}{0.48\textwidth}
        \centering
        \includegraphics[width=\textwidth]{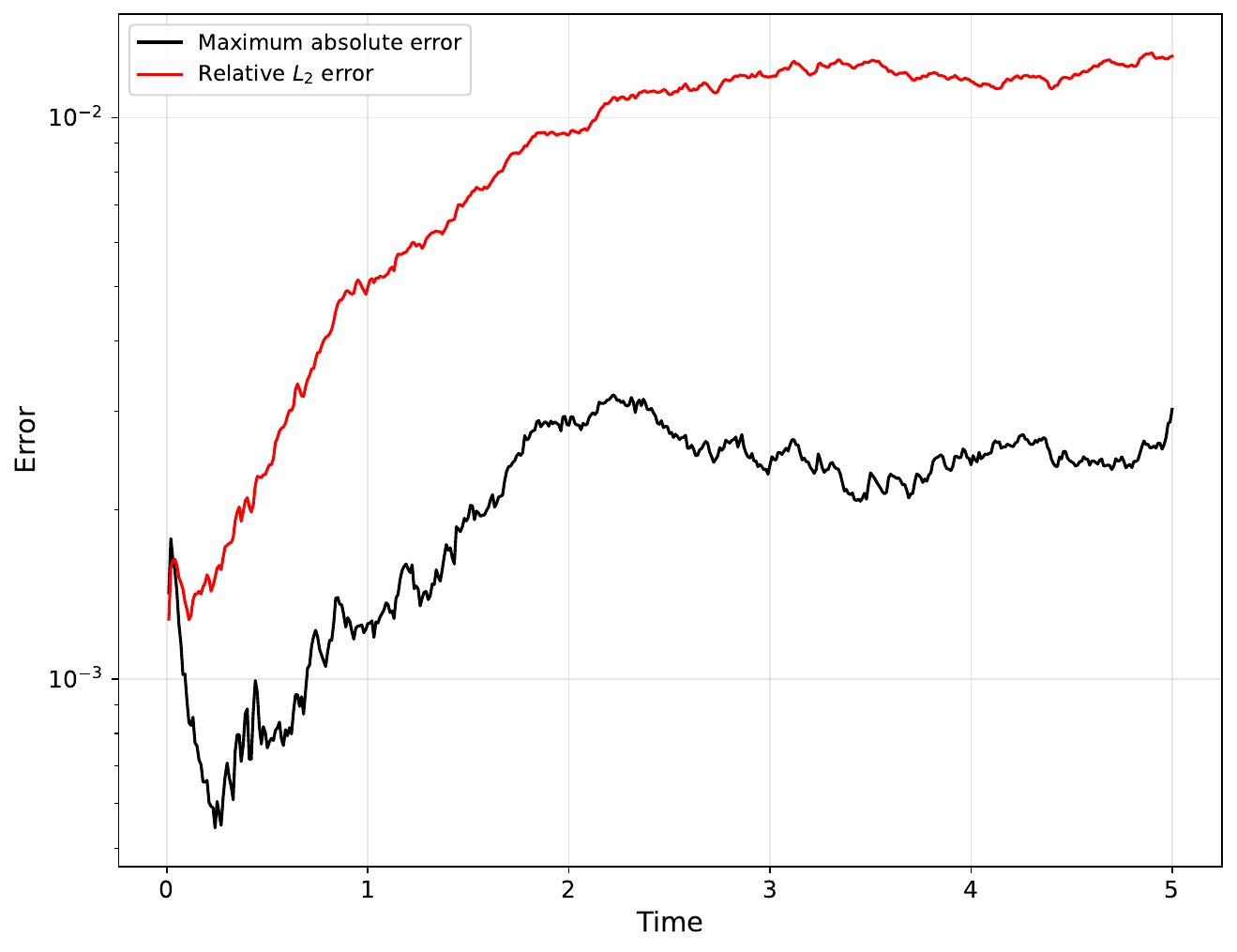}
        \caption{PINN}
        \label{fig:P4_linf_error}
    \end{subfigure}
    \hfill
    \begin{subfigure}{0.48\textwidth}
        \centering
        \includegraphics[width=\textwidth]{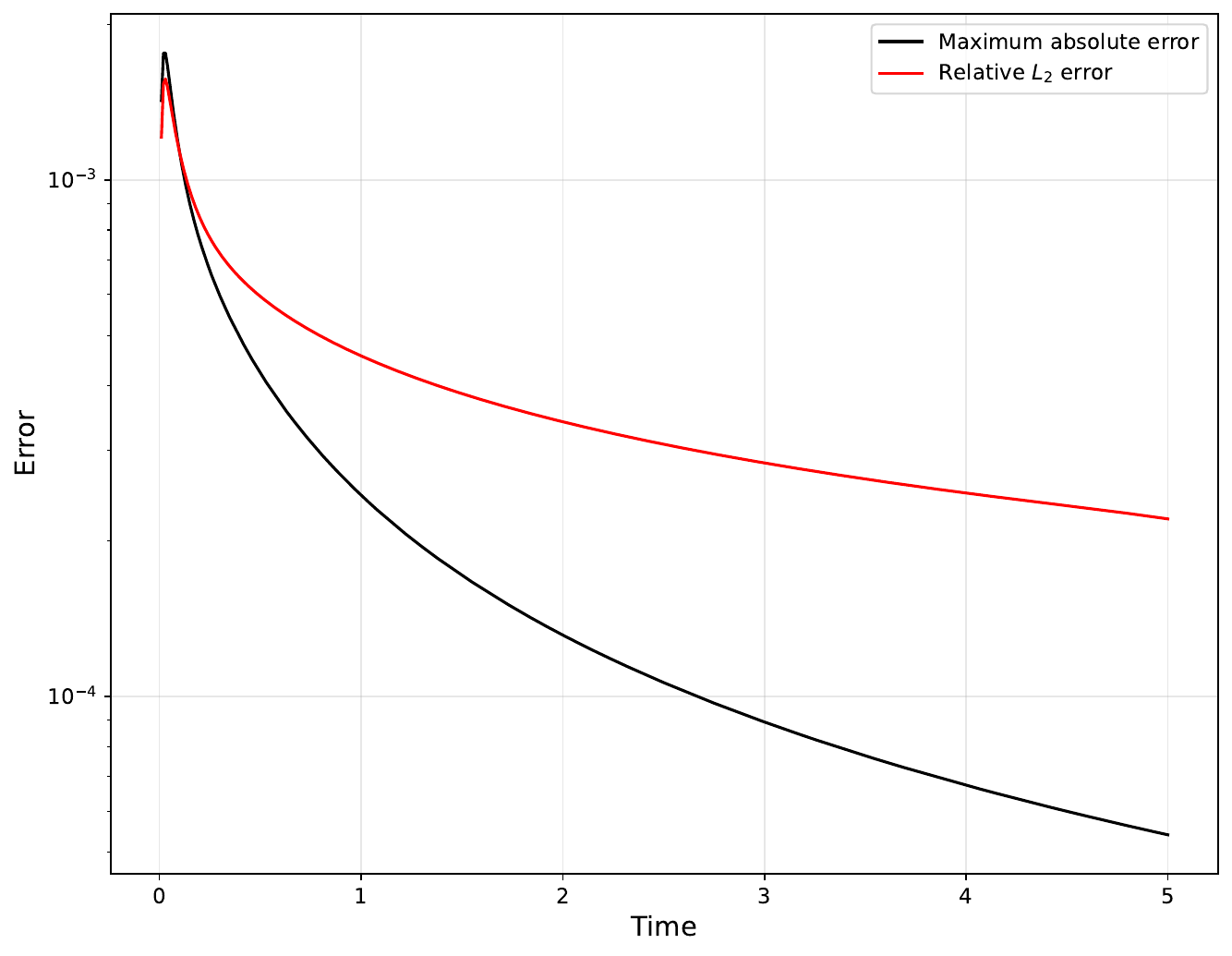}
        \caption{MELM}
        \label{fig:P4_relL2_error}
    \end{subfigure}

    \caption{Evolution of the maximum absolute error and relative $L_2$ error with time for Problem \ref{case_4}.}
    \label{fig:P4_error_vs_time}
\end{figure}

\begin{problem}[Two-Dimensional Unsteady Problem]\label{case_5}
  
Consider the two-dimensional unsteady convection--diffusion equation
\begin{equation*}
\frac{\partial \phi}{\partial t}
+
a_x\frac{\partial \phi}{\partial x}
+
a_y\frac{\partial \phi}{\partial y}
-
k\left(
\frac{\partial^2 \phi}{\partial x^2}
+
\frac{\partial^2 \phi}{\partial y^2}
\right)
=
0,
\qquad
(x,y,t)\in\Omega\times(0,T],
\end{equation*}
defined on the square domain
\(
\Omega=(0,2)\times(0,2).
\)
The exact solution is given by
\begin{equation*}
\phi(x,y,t)
=
\frac{1}{1+4t}
\exp\!\left(
-\frac{
(x-a_x t-0.5)^2+(y-a_y t-0.5)^2
}
{k(1+4t)}
\right).
\end{equation*}
The initial and Dirichlet BCs are obtained directly from the exact solution.
The problem parameters are
\begin{equation*}
U=0.8,
\qquad
a_x=0.8,
\qquad
a_y=0.8,
\qquad
k=0.01,
\qquad
L=2,
\qquad
T=1,
\end{equation*}
yielding the Péclet number
\begin{equation*}
Pe=\frac{UL}{k}=160,
\end{equation*}
where \(U=a_x=a_y=0.8\) is taken as the characteristic velocity.
\end{problem}

The FEM reference solution is computed on a uniform structured quadrilateral mesh containing $90{,}000$ elements, corresponding to a $300\times300$ element discretization. The uniform element side length is $\Delta x=\Delta y\approx 0.00667$, and no local mesh biasing is employed. This mesh provides the common spatial grid used to evaluate the CT and DT numerical solutions against the FEM reference.
\begin{figure}[H]
    \centering
    \includegraphics[width=0.5\linewidth]{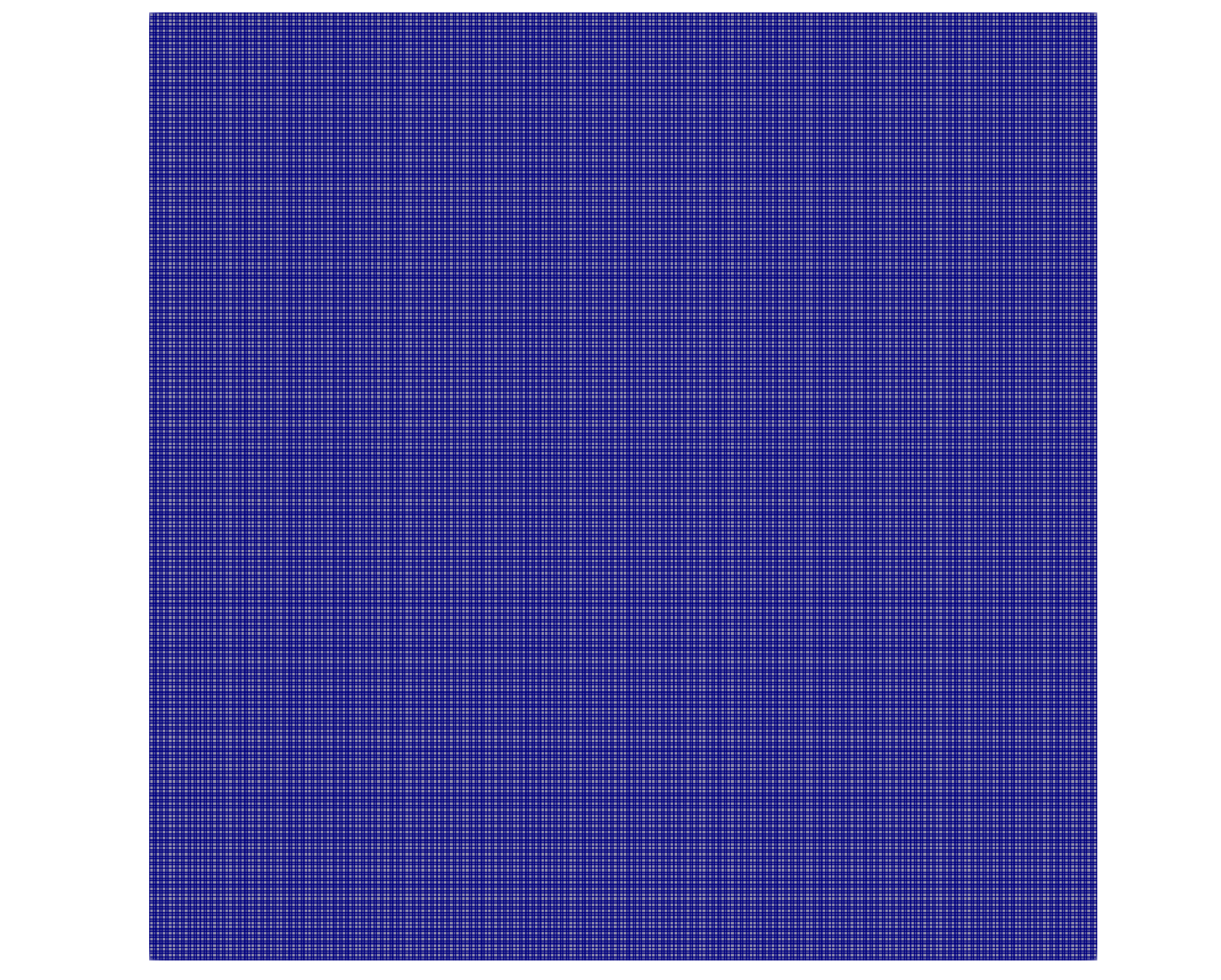}
    \caption{Structured FEM mesh used to compute the reference solution for Problem~\ref{case_5}.}
    \label{fig:P5_2dmesh}
\end{figure}

Table~\ref{tab:5_parameters} summarizes the implementation parameters used for the PINN and MELM methods for both
the CT and DT formulations.

\begin{table}[H]
\centering
\caption{Parameters used for the CT-PINN, DT-PINN, CT-MELM, and DT-MELM models for Problem~\ref{case_5}.}
\label{tab:5_parameters}

\begin{minipage}[t]{0.48\textwidth}
\centering
\textbf{PINN}\\[2mm]
\begin{tabularx}{\linewidth}{@{}>{\raggedright\arraybackslash}X
                              >{\centering\arraybackslash}X
                              >{\centering\arraybackslash}X@{}}
\toprule
\textbf{Parameter} & \textbf{CT} & \textbf{DT} \\
\midrule
$T$ & 1 & 1 \\
$\Delta t$ & --- & 0.001 \\
\midrule
\multicolumn{3}{l}{\textit{Architecture}} \\
Network type & Modified MLP & MLP \\
Layers & 5 & 3 \\
Neurons & 128 & 64 \\
Activation & GELU & GELU \\
Fourier feature scale & 10.0 & 7.0 \\
RWF & \makecell{$\mu=1.0,$\\$\sigma=0.1$}
    & \makecell{$\mu=1.0,$\\$\sigma=0.1$} \\
\midrule
\multicolumn{3}{l}{\textit{Collocation and training}} \\
Collocation points & --- & --- \\
Batch size (residual) & 256 & 1{,}024 \\
Batch size (boundary) & 256 & 256 \\
Batch size (initial) & 256 & --- \\
Learning rate & 0.001 & 0.001 \\
Training iterations & 150{,}000 & 6{,}000 \\
$\lambda_i$ (BC, PDE, IC) & $(1.0,\ 10.0,\ 1.0)$ & $(1.0,\ 7.0,\ -)$ \\
PRNG seed & 123 & 123 \\
\midrule
\multicolumn{3}{l}{\textit{Optimization}} \\
Optimizer & SOAP & SOAP \\
$\beta_1$ & 0.9 & 0.9 \\
LR schedule & Step decay & Step decay \\
LR decay rate & 0.9 & 0.9 \\
LR decay steps & 2{,}000 & 2{,}000 \\
\midrule
\multicolumn{3}{l}{\textit{Loss weighting}} \\
Weighting scheme & NTK & NTK \\
Causal tolerance & 1.0 & --- \\
Number of chunks & 16 & --- \\
\bottomrule
\end{tabularx}
\end{minipage}
\hfill
\begin{minipage}[t]{0.45\textwidth}
\centering
\textbf{MELM}\\[2mm]
\begin{tabularx}{\linewidth}{@{}>{\raggedright\arraybackslash}X
                              >{\centering\arraybackslash}X
                              >{\centering\arraybackslash}X@{}}
\toprule
\textbf{Parameter} & \textbf{CT} & \textbf{DT} \\
\midrule
$T$ & 1 & 1 \\
$\Delta t$ & --- & 0.001 \\
\midrule
\multicolumn{3}{l}{\textit{Architecture}} \\
Network type & MELM (Legendre) & MELM (Legendre) \\
Layers & 1 & 1 \\
No. of Basis functions & 52{,}992 & 2{,}500 \\
Basis & Shifted Legendre & Shifted Legendre \\
\midrule
\multicolumn{3}{l}{\textit{Collocation and solution}} \\
Collocation points & 63{,}105 & 2{,}700 \\
$\lambda_i$ (BC,\ IC) & $(1.0,\ 1.0)$ & $(1.0,\ -)$ \\
\midrule
\multicolumn{3}{l}{\textit{Linear solver}} \\
System size & $63{,}105\times52{,}992$ & $2{,}700\times2{,}500$ \\
Method & Direct least-squares & QR-based least-squares \\
\bottomrule
\end{tabularx}
\end{minipage}

\end{table}

\subsubsection*{Case 1: Continuous-Time Results}

The number of basis functions used for the MELM is $48$ along the
$x$-direction, $48$ along the $y$-direction, and $23$ along the
$t$-direction, yielding a total of
$48 \times 48 \times 23 = 52{,}992$ tensor-product basis functions.

Fig.~\ref{fig:ct5_pinn} presents the CT-PINN results at different
time instances, $t=0.20$, $0.40$, $0.60$, $0.80$, and $1.00$. The computed solutions closely follow the exact solution, reproducing the localized peak and its temporal evolution over the considered time levels. The corresponding error distributions indicate good agreement throughout
the computational domain.

Fig.~\ref{fig:ct5_melm} presents the corresponding results obtained
using the CT-MELM. The MELM also reproduces the localized peak and its temporal evolution, showing close agreement with the exact solution at the considered time levels. Compared with the CT-PINN results,
the PINN and the MELM exhibit similar solution profiles and error
distributions, enabling a direct visual comparison of their
approximation accuracy at the considered time levels.
\begin{figure}[!htbp]
    \centering
    \includegraphics[width=1\linewidth]{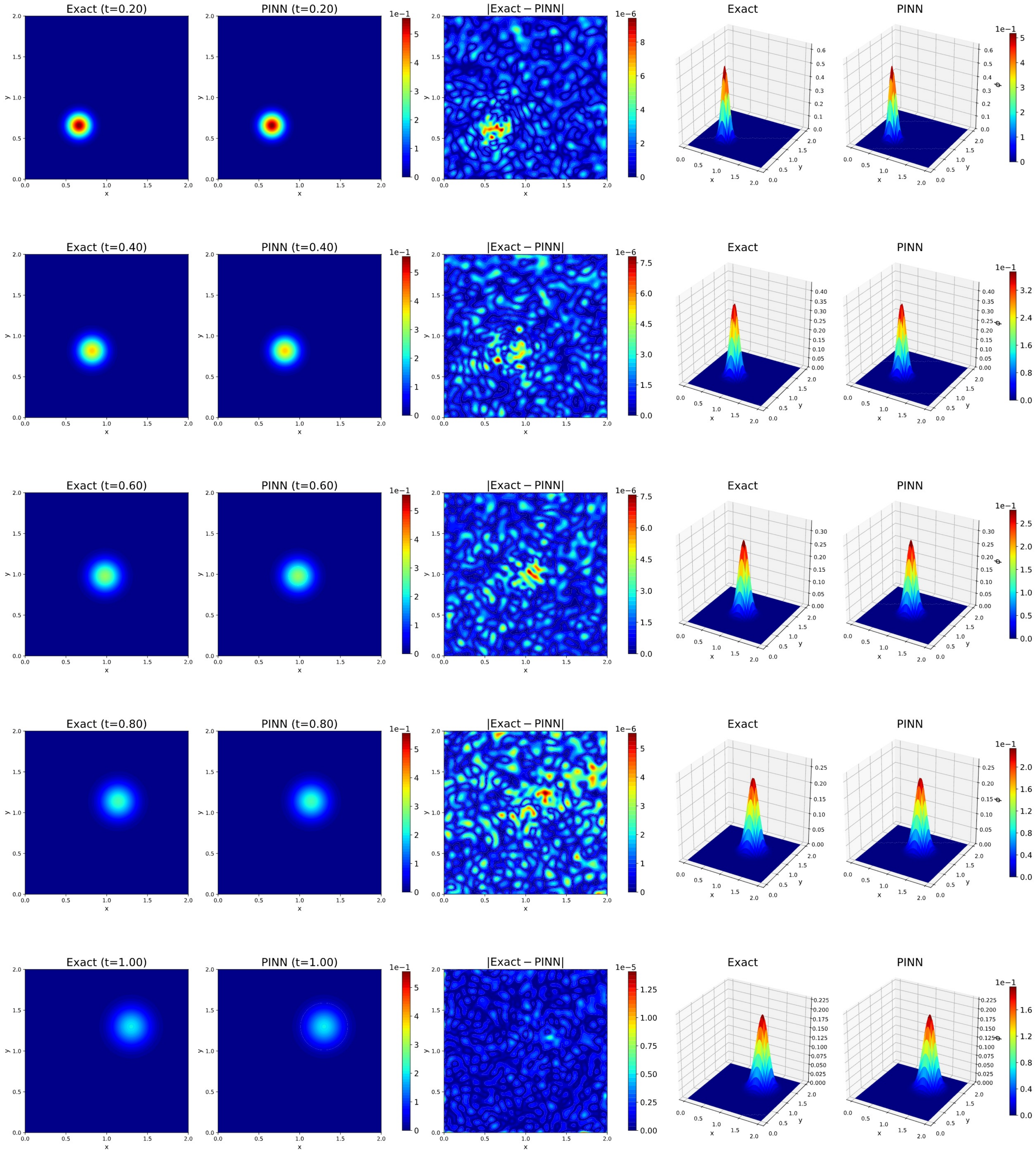}
    \caption{Problem \ref{case_5}: CT-PINN solution at $t=0.2$, $0.4$, $0.6$, $0.8$, and $1.0$ (top to bottom), compared with the exact solution.}
    \label{fig:ct5_pinn}
\end{figure}

\begin{figure}[!htbp]
    \centering
    \includegraphics[width=1\linewidth]{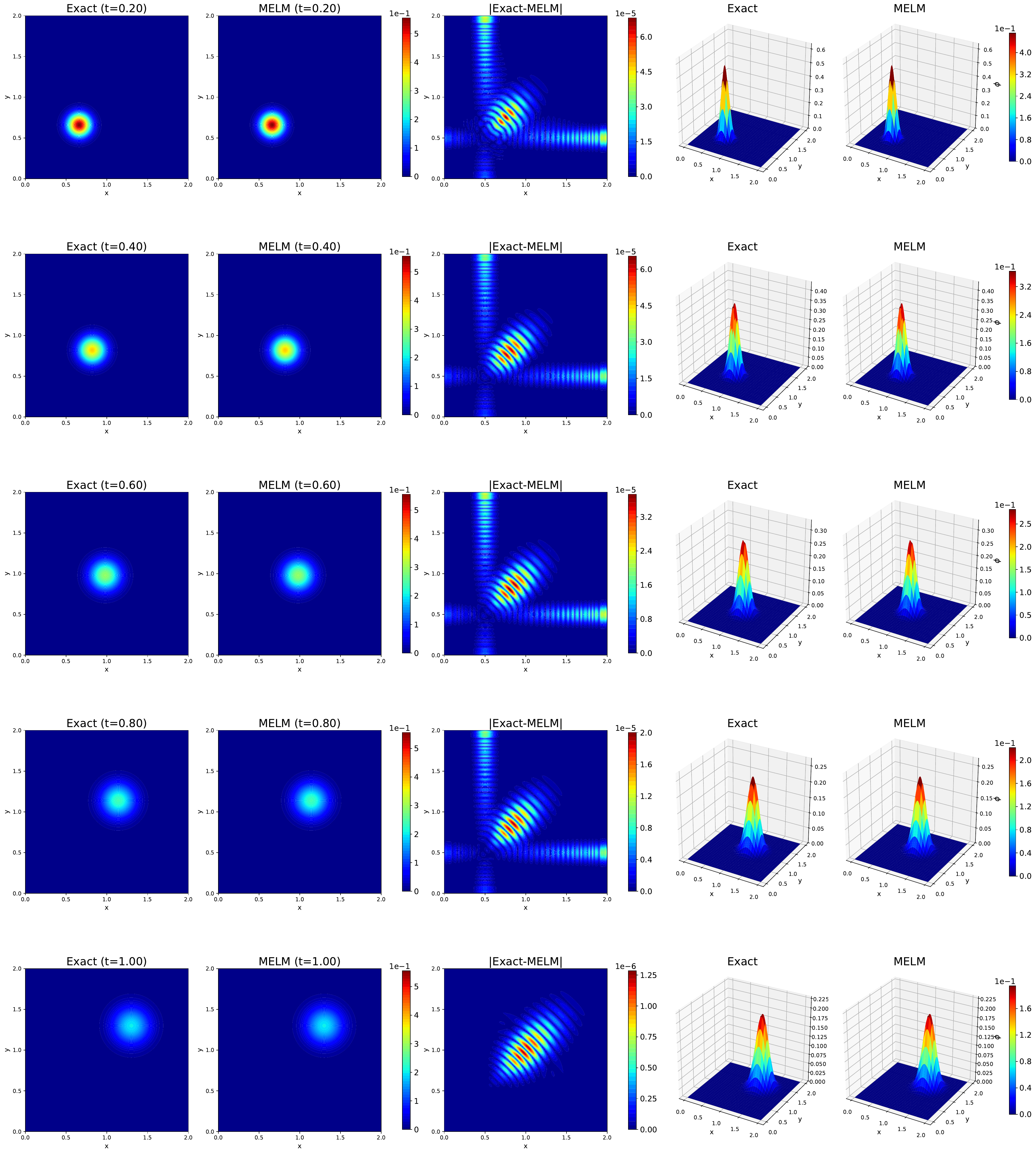}
    \caption{Problem \ref{case_5}: CT-MELM solution at $t=0.2$, $0.4$, $0.6$, $0.8$, and $1.0$ (top to bottom), compared with the exact solution.}
    \label{fig:ct5_melm}
\end{figure}

\subsubsection*{Case 2: Discrete-Time Results}

For both methods, the time step is taken as $\Delta t=0.001$. For the
MELM, $50$ basis functions are used along both the $x$- and
$y$-directions, resulting in $50\times50=2{,}500$ basis functions. A total
of $2{,}700$ collocation points are used, leading to an overdetermined
system. For the first time step, BDF1 is used to obtain the solution,
and the resulting overdetermined system is solved in the least-squares
sense. Subsequently, BDF2 is used for the remaining time steps. The
BDF2 coefficient matrix remains fixed across these time steps, while
only the right-hand side is updated. Its QR factorization is therefore
computed once and reused at each time step, thereby improving
computational efficiency.

Fig.~\ref{fig:dt5_pinn} presents the DT-PINN results for Problem~\ref{case_5} at
$t=0.20$, $0.40$, $0.60$, $0.80$, and $1.00$. The DT-PINN solution reproduces the main features of the localized peak and its temporal evolution, while the corresponding error contours show the spatial distribution of the approximation error. The three-dimensional surfaces further illustrate the evolution of the computed solution over the considered time levels.

Fig.~\ref{fig:dt5_melm} presents the DT-MELM results for Problem~\ref{case_5} at
$t=0.20$, $0.40$, $0.60$, $0.80$, and $1.00$. The MELM solution closely
captures the localized peak and its temporal evolution, showing good
agreement with the exact solution. The error contours and three-dimensional
surface plots further confirm the accuracy of the MELM approximation.

\begin{figure}
    \centering
    \includegraphics[width=1\linewidth]{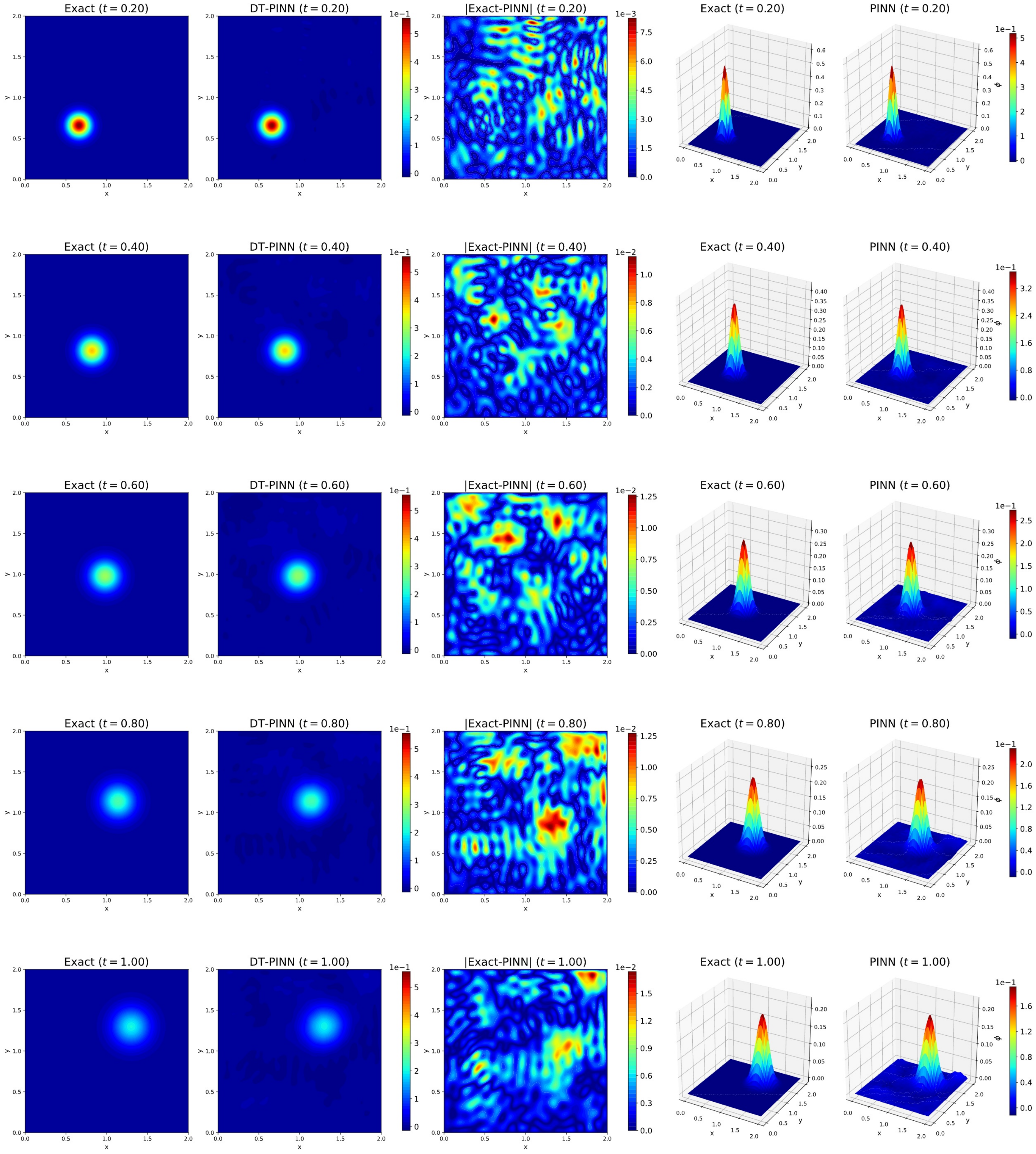}
    \caption{Problem 5: DT-PINN-BDF2 solution at $t=0.2$, $0.4$, $0.6$, $0.8$, and $1.0$ (top to bottom), compared with the exact solution.}
    \label{fig:dt5_pinn}
\end{figure}

\begin{figure}[!htbp]
    \centering
    \includegraphics[width=1\linewidth]{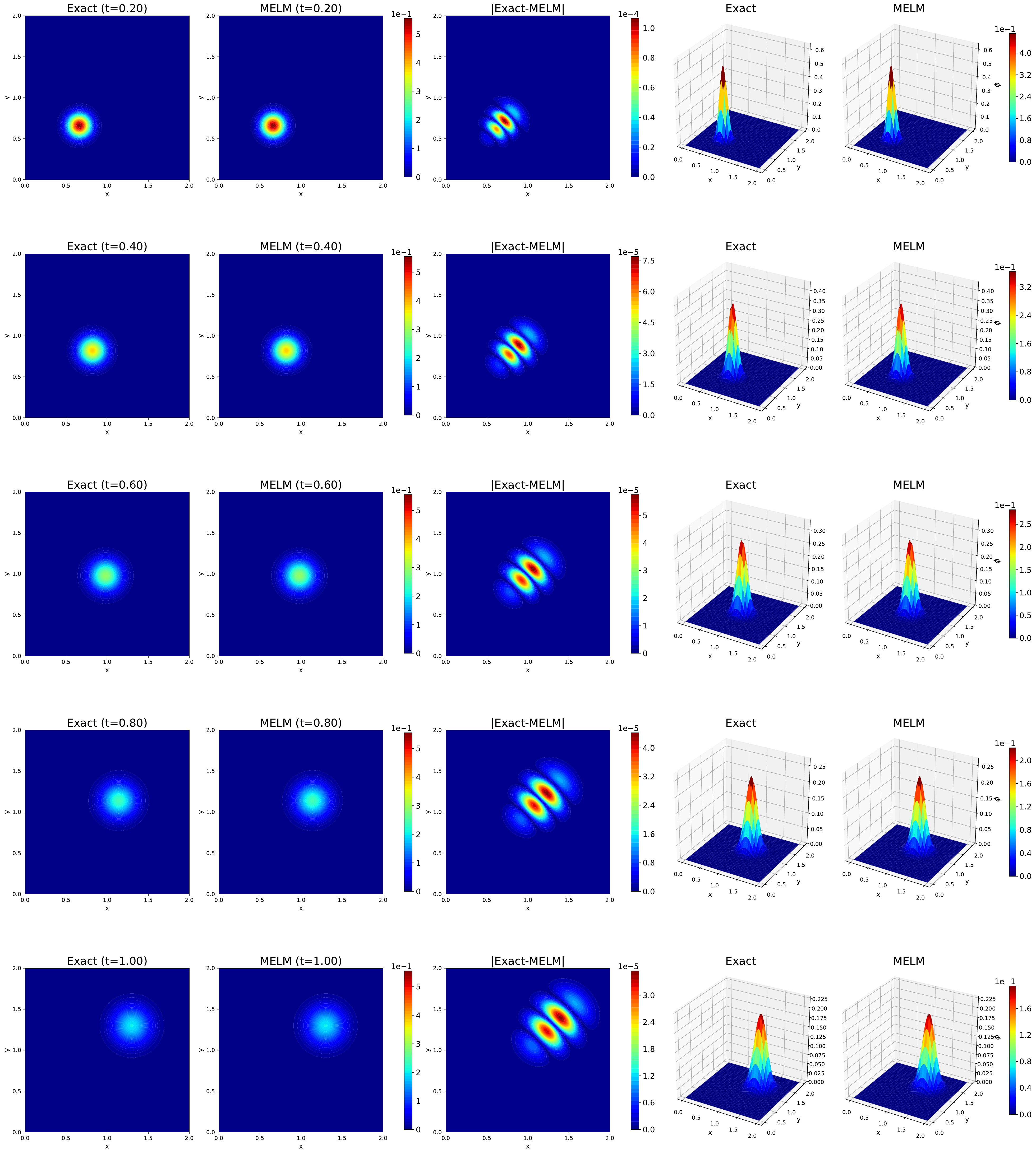}
    \caption{Problem 5: DT-MELM-BDF2 solution at $t=0.2$, $0.4$, $0.6$, $0.8$, and $1.0$ (top to bottom), compared with the exact solution.}
    \label{fig:dt5_melm}
\end{figure}

Fig.~\ref{fig:PINN_vs_MELM_error} presents the temporal
evolution of the relative $L_2$ and maximum absolute 
errors for the PINN and MELM. For the PINN, the relative $L_2$ error
increases rapidly during the initial stage and continues to grow
gradually, while the maximum absolute error remains comparatively small
with some fluctuations. In contrast, the MELM exhibits substantially
smaller errors, with the relative $L_2$ error remaining on the order of
$10^{-4}$ and the maximum absolute error decreasing steadily over time.
These results indicate that the MELM provides more accurate solutions
with more stable error behavior than the PINN throughout the simulation.
\begin{figure}[!htbp]
    \centering

    \begin{subfigure}{0.48\textwidth}
        \centering
        \includegraphics[width=\textwidth]{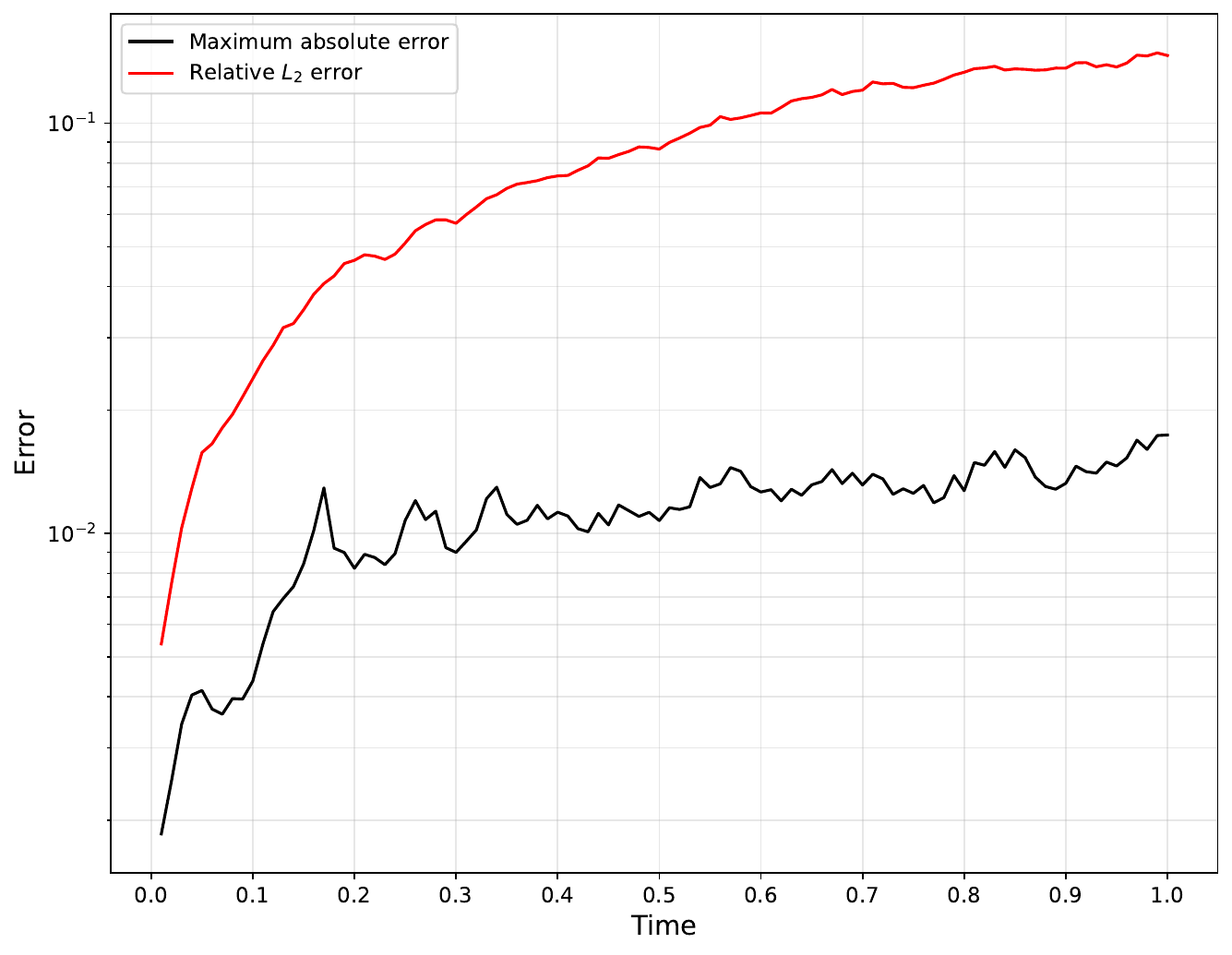}
        \caption{PINN}
        \label{fig:PINN_error_time}
    \end{subfigure}
    \hfill
    \begin{subfigure}{0.48\textwidth}
        \centering
        \includegraphics[width=\textwidth]{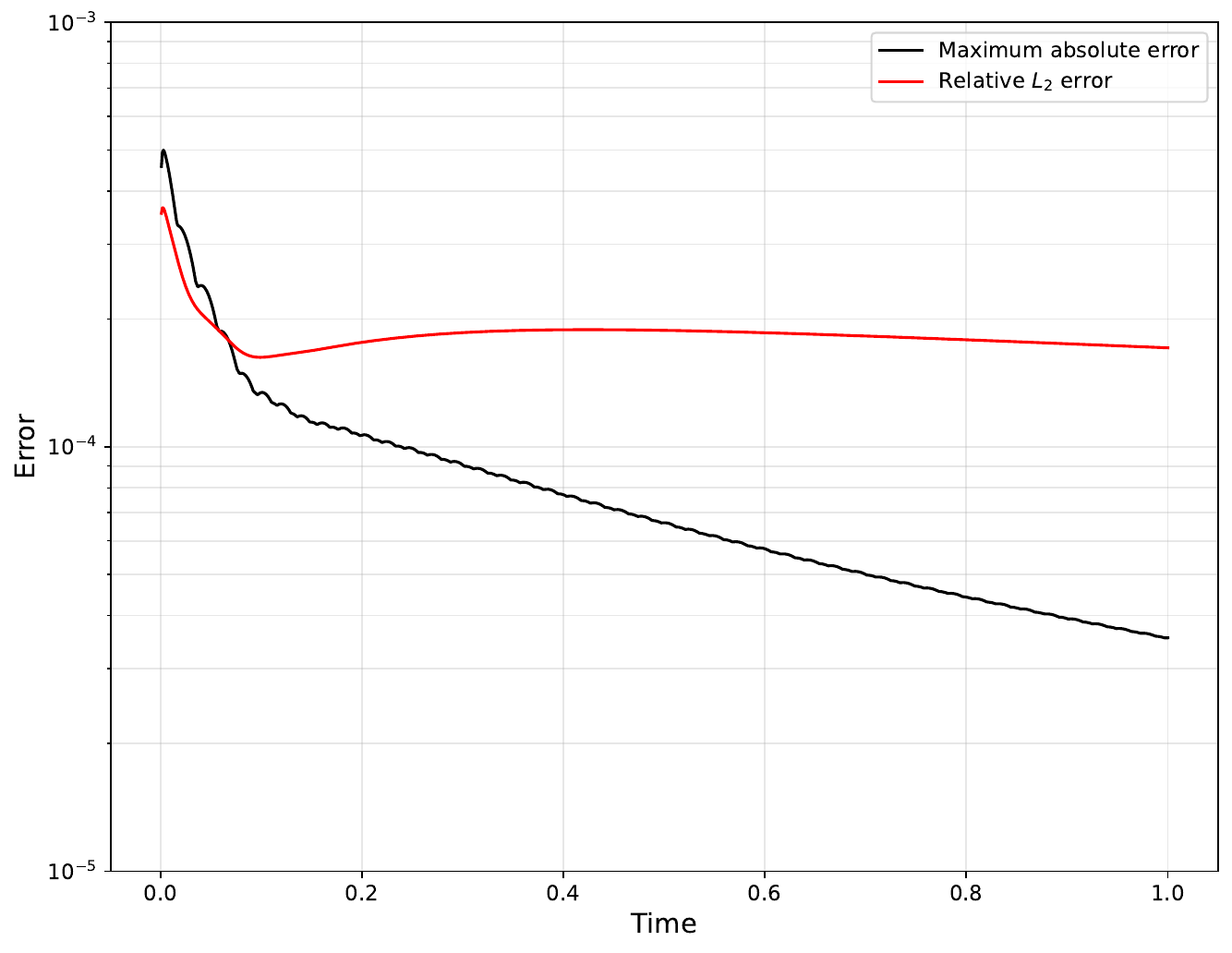}
        \caption{MELM}
        \label{fig:MELM_error_time}
    \end{subfigure}

   \caption{Evolution of the maximum absolute error and relative $L_2$ error with time for Problem \ref{case_5}.}
    \label{fig:PINN_vs_MELM_error}
\end{figure}

The quantitative comparison shows that, for the present two-dimensional space–time problem, the CT-PINN attains a lower relative $L_2$ error and a shorter computational time than the CT-MELM. This behavior differs from the trend observed in the preceding benchmarks and indicates that the global CT-MELM formulation does not scale indefinitely with the size of the space–time basis. As shown below, the DT-MELM formulation does not exhibit this limitation, because the temporal structure of the BDF discretization allows the coefficient matrix to be factorized once and reused across time levels.

\subsection*{Overall Performance Comparison}

Table~\ref{tab:comparison} summarizes the numerical accuracy and computational efficiency of the PINN, MELM, and FEM methods across all benchmark problems. The comparison includes both steady-state and transient convection--diffusion equations, with CT and DT formulations reported separately for the transient problems. For each problem, the relative $L_2$ error and total solution time are presented. The reported computational times were obtained on different computing systems and should therefore be interpreted as indicative measures of computational efficiency; the computational environments are detailed in Section~\ref{sec:sec4}, while the reported error measures are directly comparable due to the common spatial and, for transient problems, temporal evaluation grids described in this section.  It should be noted that the FEM calculations employ first-order Lagrange elements, whereas MELM uses a global high-order Legendre polynomial basis. Therefore, the reported comparisons reflect the performance of the implemented numerical approaches rather than a controlled assessment at equivalent approximation orders.

\renewcommand{\arraystretch}{1.5}
\begin{table}[!htbp]
\centering
\caption{Comparison of PINN, MELM, and FEM for all benchmark problems.}
\label{tab:comparison}
\begin{tabular}{llcccccc}
\hline
\textbf{Problem} & \textbf{Case} &
\multicolumn{3}{c}{\textbf{Relative $L_2$ Error}} &
\multicolumn{3}{c}{\textbf{Solution Time (s)}} \\
\cline{3-5}
\cline{6-8}
& & \textbf{PINN} & \textbf{MELM} & \textbf{FEM}
  & \textbf{PINN} & \textbf{MELM} & \textbf{FEM} \\
\hline

Problem \ref{case_1} & --
& $8.31\times10^{-5}$ & $8.22\times10^{-14}$ & $2.34\times10^{-5}$
& $5{,}730$ & $0.62$ & $0.13$ \\

Problem \ref{case_2} & --
& $6.44\times10^{-3}$ & $4.96\times10^{-3}$ & --
& $3{,}650$ & $481$ & $1.02^{\dagger}$ \\

Problem \ref{case_3} & --
& $1.03\times10^{-2}$ & $9.86\times10^{-3}$ & --
& $9{,}171$ & $464$ & $0.92^{\dagger}$ \\

Problem \ref{case_4} & CT
& $9.90\times10^{-6}$ & $1.14\times10^{-9}$ & --
& $1{,}619$ & $113$ & -- \\

& DT
& $8.65\times10^{-3}$ & $5.40\times10^{-4}$ & $7.55\times10^{-4}$
& $81{,}913$ & $0.26$ & 29.85 \\

Problem \ref{case_5} & CT
& $3.90\times10^{-5}$ & $2.77\times10^{-4}$ & --
& $3{,}814$ & $40{,}284$ & -- \\

& DT
& $8.06\times10^{-2}$ & $1.91\times10^{-4}$ & $2.93\times10^{-3}$
& $231{,}451$ & $22.59$ &  $1{,}860$ \\
\hline
\end{tabular}
\\
\footnotesize{$\dagger$ In Problems~\ref{case_2} and~\ref{case_3} no analytical
solution is available; the FEM solution is used as the numerical reference for
the PINN and MELM approximations, and the reported time corresponds to the cost
of computing that reference. In the remaining cases, FEM is included as an
independent classical solution for comparison.}
\end{table}
\renewcommand{\arraystretch}{1}

\subsubsection*{Overall Observations}

\begin{itemize}

\item \textbf{Higher accuracy with lower computational cost:} For the benchmark problems considered, the MELM generally achieves lower relative $L_2$ errors than the PINN while requiring substantially lower wall-clock time. The accuracy advantage is particularly pronounced for the DT formulations of Problems~\ref{case_4} and~\ref{case_5}. These results demonstrate that the proposed MELM can provide high solution accuracy at substantially lower computational cost than the PINN for the linear convection--diffusion problems considered.

\item 
\textbf{Behavior near discontinuities:}
Problem~\ref{case_3} highlights differences in local approximation behavior. MELM exhibits Gibbs-type oscillations, whereas the PINN produces a smoother solution with larger pointwise deviations near the discontinuities. This illustrates that the absence of visible oscillations does not necessarily imply higher local accuracy.

\item \textbf{Advantage of DT-MELM for time-dependent problems:} For the time-dependent problems considered, DT-MELM provides a substantial computational advantage over CT-MELM. This advantage becomes particularly pronounced for Problem~\ref{case_5}, where the larger space--time system leads to a marked increase in the computational cost of the CT formulation.

\item \textbf{Matrix reuse in DT-MELM:} For the linear time-dependent problems considered, the coefficient matrix remains fixed across discrete time levels, while only the right-hand side changes. Hence, the matrix can be assembled and factorized once and reused at subsequent time levels, which contributes substantially to the computational efficiency of DT-MELM. For nonlinear problems, this property is generally lost because the resulting system depends on the current solution and may require an updated or iterative solution procedure at each time level.

\item \textbf{Scalability:} Problem~\ref{case_5} highlights the different computational scaling of the CT and DT MELM formulations for a larger space--time problem. The CT-MELM requires $40{,}284$~s, whereas the DT-MELM requires only $22.59$~s, corresponding to an approximately $1.78\times10^{3}$-fold reduction in computational time. This result illustrates the computational benefit of the time-marching formulation when the spatial operator is independent of time.

\item \textbf{Performance of Problem~\ref{case_5} in the CT formulation:} Problem~\ref{case_5} in the CT formulation represents a clear exception to the general accuracy--efficiency trend observed for MELM. In this case, the CT-PINN achieves a lower relative $L_2$ error of $3.90\times10^{-5}$ compared with $2.77\times10^{-4}$ for CT-MELM, while also requiring substantially less computational time ($3{,}814$~s compared with $40{,}284$~s). This result indicates that the computational advantage of MELM is not uniform across all formulations and problem sizes, in contrast with the DT formulation discussed above.

\end{itemize}

\section{Conclusion}
\label{sec:sec6}

This study presents a Legendre-based MELM framework for solving convection--diffusion equations under challenging convection-dominated conditions. The proposed framework is investigated in both CT and DT formulations and systematically compared with PINNs for steady and unsteady problems with relatively high Péclet numbers and both continuous and discontinuous BCs. The main conclusions of the study are summarized below.

First, the proposed MELM framework exploits the linearity of the governing equation across the steady, CT, and DT formulations. For the steady and CT formulations, this allows the unknown output weights to be obtained through a single least-squares solve rather than iterative training. The deterministic construction of the MELM also reduces the
dependence on training-related hyperparameters and initialization.
Moreover, the Legendre-based spectral approximation effectively
captures steep gradients and sharp boundary layers across all the
problems considered.

Second, the numerical results demonstrate that the MELM generally achieves competitive or higher accuracy than PINNs across the benchmark problems considered, while requiring substantially lower computational time in most cases. The computational advantage is particularly pronounced for the DT formulations, especially for Problem~\ref{case_5}. However, the CT formulation of Problem~\ref{case_5} represents a clear exception, where the CT-PINN achieves both lower relative $L_2$ error and shorter computational time than the CT-MELM. These results highlight the potential of MELM as an efficient alternative to iterative PINN-based solvers for linear convection--diffusion problems, while also indicating that its computational advantage depends on the problem formulation and size.

Third, for the linear time-dependent problems, the matrix-reuse capability of DT-MELM provides the most significant computational advantage of the proposed framework. Since the coefficient matrix remains unchanged across successive BDF2 time levels while only the right-hand side varies, the matrix factorization can be computed once and reused at subsequent time levels. This substantially reduces the computational cost compared with repeatedly factorizing the corresponding systems, with the benefit becoming particularly pronounced for large space--time problems.

Fourth, the matrix-reuse advantage identified in this work is specific to the linear time-dependent setting. For nonlinear time-dependent PDEs, the resulting algebraic system generally depends on the current solution, so the fixed-matrix property exploited by DT-MELM may no longer hold and an updated or iterative solution procedure may be required at each time level. Extending MELM to nonlinear time-dependent PDEs and systematically assessing its computational efficiency relative to PINNs therefore remains an important direction for future research.

Finally, beyond nonlinear extensions, the MELM framework opens the door to using PINNs with more complex activation functions than usual, while offering a natural foundation for future research to explore multi-layer architectures and the introduction of trainable hidden-layer biases.

\section{Declarations}
	
\textbf{Competing Interest:} The authors have no competing interests to declare that are relevant to the content of this article. \\	
	
\textbf{Data Availability:} There is no associated data available for our manuscript.\\

\textbf{Funding:}  The second author wishes to thank WSAI, IIT Madras for their financial support.\\

\textbf{CRediT authorship contribution statement:} \textbf{Gadade Navnath Ankush:} Writing -- original draft, Validation, Software, Methodology, Investigation, Formal analysis. \textbf{Ram Shiromani:} Writing -- original draft, Validation, Methodology, Investigation, Formal analysis. \textbf{Tom\'as Gandara:} Writing -- original draft, Validation, Software, Methodology, Investigation, Formal analysis. \textbf{Ramon Codina:} Writing -- review \& editing, Supervision, Methodology, Investigation, Formal analysis, Conceptualization. \textbf{Sivaram Ambikasaran:} Writing -- review \& editing, Supervision, Formal analysis, Conceptualization.\\

\textbf{Ethical Approval:} Not applicable.

\nocite{*}

\bibliography{references}
\end{document}